\documentclass[onefignum,onetabnum]{siamart251216} % for arXiv version
\usepackage{lipsum}
\usepackage{amsfonts}
\usepackage{graphicx}
\usepackage{epstopdf}
\usepackage{algorithmic}
\ifpdf
  \DeclareGraphicsExtensions{.eps,.pdf,.png,.jpg}
\else
  \DeclareGraphicsExtensions{.eps}
\fi

\newsiamremark{remark}{Remark}
\newsiamremark{hypothesis}{Hypothesis}
\crefname{hypothesis}{Hypothesis}{Hypotheses}
\newsiamthm{claim}{Claim}
\newsiamremark{fact}{Fact}
\crefname{fact}{Fact}{Facts}

\usepackage{amsopn}

\usepackage{hyperref}
\usepackage{graphicx}
\usepackage{subcaption}
\usepackage{multicol}
\usepackage{xcolor}
\usepackage{tikz}
\usepackage{algorithm}
\usepackage{algorithmic}
\usepackage{multirow}
\usepackage{float}
\usepackage{enumitem}
\usepackage{booktabs}
\usepackage{xstring}
\usepackage{mathtools}
\usepackage{tabularx}
\usetikzlibrary{positioning, arrows.meta, calc, shapes.geometric, fit, backgrounds, decorations, decorations.pathreplacing, decorations.markings, decorations.pathmorphing, decorations.text, calligraphy, babel}
\newcommand{\mycomment}[1]{\hfill $\triangleright$ #1}
\definecolor{lowrank}{RGB}{0,200,230}   % cyan-like
\definecolor{fullrank}{RGB}{200,30,30}  % red-like
\definecolor{identblue}{RGB}{40,100,210}
\definecolor{framecol}{gray}{0.15}

\usepackage{standalone}

\newcommand{\sbkt}[1]{\left[#1\right]}
\newcommand{\bkt}[1]{\left(#1\right)}

\newcommand{\floor}[1]{\left\lfloor#1\right\rfloor}

\newcommand{\cbkt}[1]{\left\{#1\right\}}
\newcommand{\norm}[2][]{\left\lVert #2 \right\rVert_{#1}}
\numberwithin{equation}{section}
\newcommand{\N}{\mathbb{N}}

\definecolor{lightblue}{rgb}{0, 0.3, 0.9} % RGB values 

\headers{A fast direct solver based neural network for solving PDEs}{Jashwanth Reddy Kadaru and Vaishnavi Gujjula}

\title{A fast direct solver based neural network for solving PDEs
}

\author{Jashwanth Reddy Kadaru\thanks{ Department of Computer Science \& Engineering, International Institute of Information Technology Bangalore (IIIT-B), India 
  (\email{Jashwanth.Kadaru095@iiitb.ac.in}).}
\and Vaishnavi Gujjula\thanks{Department of Data Science and Artificial Intelligence, International Institute of Information Technology Bangalore (IIIT-B), India 
  (\email{vaishnavi.gujjula@iiitb.ac.in}).}
}

\ifpdf
\hypersetup{
  pdftitle={A neural network for solving PDEs based on the fast direct solver for HODLR matrices},
  pdfauthor={Jashwanth Reddy Kadaru and Vaishnavi Gujjula}
}
\fi

\begin{document}
% \nolinenumbers
\maketitle

% REQUIRED
\begin{abstract}
The matrices arising from large scale $N$-body problems can be efficiently represented using hierarchical matrices, whose key idea is that the admissible off-diagonal sub-matrices can be well approximated by low-rank matrices across a hierarchy of matrix partitions.
HODLR (Hierarchical Off-Diagonal Low-Rank) matrices are a subclass of hierarchical matrices in which all off-diagonal submatrices at every level of a recursive binary partition are low-rank.
In this article, we present a neural network that learns the inverse operation of HODLR matrices based on the fast direct solver for HODLR matrices developed by Ambikasaran and Darve (2013).
% Further, we extend the network to solve non-linear PDEs, 
We further extend the architecture to learn nonlinear solution operators associated with PDEs by replacing some of the linear layers with deep sub-networks.
% containing nonlinear activations.
We demonstrate the performance of the proposed architecture by performing a comprehensive set of experiments that include
(i)~solving a linear problem such as the Fredholm integral equation of the second kind,
(ii)~solving PDEs such as the nonlinear Schr\"odinger equation, Burgers' equation, and the steady-state Darcy's flow equation,
(iii)~generalization study across varying parameter values,
(iv)~comparing the inference time of the proposed network with the run time of a classical numerical solver, and
(v)~comparing the proposed network with some of the existing neural operator learning networks.
% The key idea is to encode each level of the recursive Sherman--Morrison--Woodbury factorization as a residual block composed of locally connected (LC) neural layers, yielding a network whose parameter count scales as $O(Np\log N)$, where $p$ is the off-diagonal rank.
% The proposed architecture shares the same structure of the fast direct solver inverse.
\end{abstract}

% REQUIRED
\begin{keywords}
Hierarchical matrices, HODLR matrices, fast direct solver, locally connected neural network, convolutional neural network, PDE solvers, operator learning
\end{keywords}

% REQUIRED
\begin{MSCcodes}
68T07, 65R20, 65F30
\end{MSCcodes}

\section{Introduction}
Physical problems governed by linear partial differential equations can be formulated using either integral equation (IE) or differential equation approaches. Discretization of the IE formulation typically results in a dense linear system, whose size becomes huge when solving large-scale problems or when high accuracy is required. Consequently, the use of naive linear algebra solvers for such dense linear systems is computationally expensive. Motivated by these challenges, the development of efficient algorithms for solving large-scale linear systems arising from the discretization of integral equations has been an active area of research in computational physics.

Linear systems can generally be solved using either iterative or direct methods. Iterative methods, such as Jacobi iteration, conjugate gradient, GMRES, and MINRES, are widely used for large systems.
However, the performance of iterative methods is highly sensitive to matrix conditioning, and they become inefficient when multiple right-hand sides need to be processed. 
% Example iterative solvers include Jacobi iteration, conjugate gradient, GMRES, and MINRES. 
The computationally intensive step in the iterative solvers is the computation of matrix-vector product in each iteration. To accelerate computation, fast matrix-vector product algorithms are to be used. 

An alternative approach is to use direct solvers. However, conventional direct methods such as LU and QR factorizations are computationally expensive for large dense systems. To accelerate computation, fast direct solvers are to be used.
% have been developed to significantly reduce computational complexity while maintaining high accuracy.

These computational challenges have motivated the development of multiscale frameworks that exploit the inherent structure of physical operators. 
Multiscale methods including multigrid methods~\cite{brandt1977multi}, the Fast Multipole Method (FMM)~\cite{greengard1987fast}, wavelet-based methods \cite{beylkin1991fast} and hierarchical matrix frameworks~\cite{borm2003introduction, grasedyck2003construction} exploit the underlying multiscale structure to reduce computational complexity. 

% The hierarchical matrices introduced by Hackbusch \cite{hackbusch1999sparse,hackbusch2015hierarchical} provide a robust algebraic framework for approximating dense matrices by exploiting the low-rank structure. 
In the hierarchical matrix framework introduced by Hackbusch~\cite{hackbusch1999sparse,hackbusch2015hierarchical}, the computational domain is organized into a hierarchical tree of sub-domains. The resulting system matrix is composed of sub-matrices that represent the interactions between the sub-domains at a given level of hierarchy. Most of these sub-matrices correspond to interactions between far-field sub-domains, while a few correspond to interactions between near-field sub-domains. The sub-matrices corresponding to the far-field interactions exhibit low-rank. 
Furthermore, a near-field interaction at a given level in the hierarchy can be broken down into interactions among their respective subdomains. At a finer level, some of these interactions are far-field interactions, which can be well represented using low-rank approximations, and a few are near-field interactions. 
% Thus, using this framework we can represent most of the off-diagonal blocks of the system matrix using low-rank sub-matrices. 
Consequently, as the hierarchical refinement proceeds deeper into the tree, the number of low-rank sub-matrices increases. This low-rank structure can be leveraged to develop fast direct and iterative solvers.
This is the key idea behind the hierarchical matrix framework. The fast iterative and direct solvers developed for hierarchical matrices can achieve near-linear (quasi-linear) complexity \cite{martinsson2005fast,chandrasekaran2006fast,bebendorf2008hierarchical,ying2004kernel,fong2009blackbox,coulier2017ifmm}. 

% Developing solution methods for nonlinear PDEs is an active area of research, with significant effort directed towards reducing their computational cost.
Developing computationally efficient numerical methods for evaluating non-linear solution operators associated with PDEs is an active area of research.
% , with significant effort directed towards reducing their computational cost.
In general, these methods are problem-specific, and constructing them is challenging due to the complexity of the underlying solution map.

The combination of numerical methods and deep learning has emerged as a promising direction for learning non-linear solution operators associated with PDEs. Deep learning approaches offer powerful and flexible function approximators capable of handling high-dimensional non-linearities \cite{raissi2019physics}. By the universal approximation theorem, a sufficiently wide feedforward neural network with a non-polynomial activation can uniformly approximate any continuous function on a compact domain to arbitrary accuracy under mild conditions \cite{cybenko1989,hornik1989}. However, selecting an appropriate architecture for a given problem remains challenging and is an active area of research.

In recent years, several neural networks have been proposed for learning non-linear solution operators associated with PDEs. 
The Fourier Neural Operator (FNO) of Li et al.~\cite{li2021fourier} parameterizes integral kernel operators through learned Fourier-mode coefficients, enabling efficient resolution-invariant operator learning.
% The Fourier Neural Operator (FNO) of Li et al.~\cite{li2021fourier} learns operators in Fourier space by parameterizing integral kernel layers via their Fourier coefficients, enabling efficient resolution-invariant learning. 
DeepONet~\cite{lu2021deeponet} uses a branch-trunk architecture to approximate infinite-dimensional operators. A theoretical framework for neural operator approximation has been established by Kovachki et al.~\cite{kovachki2023neural}. 
% However, these general-purpose architectures treat the underlying discretization as a black box, ignoring the rich multiscale algebraic structure embedded in the operators.
% Another recent work has sought to encode multiscale structures into model architectures. 
A notable recent advancement in incorporating multiscale structures into model architectures is the multiscale neural network based on $\mathcal{H}$-matrices (from here on abbreviated as MNN) and $\mathcal{H}2$-matrices proposed by Fan et al.~\cite{fan2019multiscale,fan2019multiscale_h2}. In these articles, the authors build (i) a neural network that learns the forward operation of the hierarchical matrix (hierarchical matrix-vector product) and then (ii) extend it to learn non-linear solution operators associated with PDEs.
% by replacing selected linear sub-networks with local deep neural networks.
% By interpreting the hierarchical block structure of $\mathcal{H}$-matrices as a cascade of locally connected neural layers (convolutional layers in the case of translation invariance) with linear activations, the multiscale neural network efficiently learns the forward operation of a matrix (matrix-vector multiplication). It was further extended to solve the nonlinear PDEs by replacing the linear activations with non-linear activations. 
% Related efforts on learning structured matrix factorizations include butterfly factorizations \cite{dao2019butterfly,li2020butterfly}, which target matrices with complementary low-rank structure (i.e., off-diagonal blocks that are low-rank in a Fourier or frequency-domain sense).
% ; these are complementary to the HODLR structure exploited here, as butterfly matrices arise in oscillatory kernels where spatial-domain off-diagonal blocks may not be low-rank. 

In this article, similar to the work by Fan et al.~\cite{fan2019multiscale}, we design a neural network inspired by hierarchical matrices. However, while Fan et al.~\cite{fan2019multiscale} introduced a network to learn the forward operation of a hierarchical matrix, we introduce 
\begin{enumerate}[leftmargin=1.2em]
\item 
    A neural network to learn the inverse operation of a hierarchical matrix. By learning the inverse operator, our model effectively functions as a fast direct solver. 
    While learning the forward map of a matrix (matrix-vector product) is valuable, the primary computational bottleneck in many forward and inverse problems is solving the associated linear system.
    % For a matrix of size $N\times N$, a matrix-vector multiplication costs $\mathcal{O}(N^{2})$ whereas a solving a linear system of equations costs $\mathcal{O}(N^{3})$ when done naively. 
    To the best of our knowledge, this is the first neural network approach structured to invert a hierarchical matrix or, equivalently, to solve the underlying linear system it represents. 
    % Te computational complexity of their approach is $O(N \log N)$. 
    After training, the network takes a right hand side vector as input and, via a forward pass, solves the corresponding linear system.

\item 
    Further, we extend this neural network to learn non-linear solution operators associated with PDEs by replacing selected linear layers with deeper sub-networks.
\end{enumerate}

    The architecture of the proposed neural network is inspired from the fast direct solver developed by Ambikasaran and Darve \cite{ambikasaran2013fast} for a sub-class of hierarchical matrices called HODLR (Hierarchically Off-Diagonal Low-Rank) matrices. HODLR matrices are constructed based on the weak admissibility criterion, in which all interactions other than self-interactions are approximated as low-rank. As a consequence, at each level of the hierarchy, all the off-diagonal blocks of the matrix are approximated as low-rank. Their approach avoids the prohibitive cost of dense matrix inversion by constructing a recursive factorization of the matrix into a sequence of block-diagonal factors. 
    % The inverse is then applied efficiently via the Sherman-Morrison-Woodbury (SMW) identity \cite{woodbury1950}. 
    After constructing the factorization, solutions to new right-hand sides can be obtained at low computational cost, making the method attractive for problems involving multiple right-hand sides.

We evaluate the performance of the proposed neural network through a comprehensive set of experiments that includes:
\begin{enumerate}[leftmargin=1.2em]
    \item Solving the linear system arising from the Fredholm integral equation.
    \item Solving nonlinear PDEs, such as the nonlinear Schr\"odinger equation (NLSE), Burgers' equation, and the linear Darcy flow problem. 
    The proposed network approximates the solutions well to a relative error in the order of $10^{-4}$ to $10^{-2}$.
    \item Solving PDEs across parametric sweeps. For the NLSE 1D problem, the proposed network is observed to generalize well across a range of parameter values. When trained on a dataset spanning multiple parameter values, the network accurately infers solutions corresponding to unseen parameter values.
    
    Further, the total offline cost of our neural network comprises the time required for data generation and training. Although this cost is higher than that of a classical numerical solver, it is effectively amortized in settings involving non-linear operators and extensive parametric sweeps. In such scenarios, classical solvers must be re-run for each parameter instance. In contrast, the proposed neural network learns over a range of parameters during training, thereby amortizing this cost and enabling faster inference. We observe that, for the NLSE 1D problem, the inference time of the neural network is significantly lower than the run time of a classical numerical solver, making it particularly advantageous in parametric sweeps.
    % \item Comparison of inference time of the proposed neural network with the run time of its classical numerical solver. The inference time of the proposed network is observed to be lower than the run time of its classical numerical solver counterpart, highlighting its advantage
    \item Comparison of the performances of the proposed neural network with some of the existing operator learning networks for PDEs, such as the FNO, DeepONet, and MNN. We observe that the proposed neural network performs competitively with its counterparts at a cheaper parameter count.
\end{enumerate}
The rest of the paper is organized as follows. Section~\ref{sec:directsolver} presents the fast direct solver for linear systems involving HODLR matrix. Section~\ref{sec:direct} presents the neural network architecture that models the direct solver. Section~\ref{sec:PDE_solver} extends the direct solver described in Section~\ref{sec:direct} to solve PDEs. Section~\ref{sec:numerical} presents numerical results on both linear and nonlinear maps. 
% Section~\ref{sec:ablation} shows that our proposed architecture (FDSNet) performs competitively with other PDE solvers across both 1D and 2D problems at a cheaper parameter count, with the gap between other baseline architectures notably narrow on certain problems. 
Finally, Section~\ref{sec:conclusion} provides conclusions and possible future directions.

\section{Fast direct solver for HODLR matrices}\label{sec:directsolver}
The architecture of the neural network proposed in this paper is based on the computational framework of the fast direct solver for HODLR (Hierarchical Off-Diagonal Low-Rank) matrices, introduced in~\cite{ambikasaran2013fast}. In this section, we briefly describe the HODLR matrix structure and the direct solver algorithm. We refer the readers to~\cite{ambikasaran2013fast} for more details.

\subsection{HODLR matrix}
% \begin{definition}[Hierarchical binary tree]
% Let $K$ be the matrix arising out of kernel functions $K(x,y):\mathbb{R}^{d} \times \mathbb{R}^{d}\rightarrow\mathbb{R}$ which is smooth everywhere except on the line $x=y$. Such matrices are often seen in the integral equation formulation of PDEs. Let $N$ be the number of degrees of freedom. Let $\mathcal{I} = \{1, 2, \dots, N\}$ be the index set associated with the set of discretization points in the physical domain $\Omega \subseteq \mathbb{R}^d$. 
Let $K$ be the matrix arising from the discretization of a kernel function $g(\mathbf{x},\mathbf{y}): \mathbb{R}^{d} \times \mathbb{R}^{d} \rightarrow \mathbb{R}$, which is smooth everywhere except along the line $\mathbf{x}=\mathbf{y}$. Such matrices commonly arise in the discretized integral equation formulation of linear elliptic PDEs. Let $N$ denote the number of degrees of freedom, and let $\mathcal{I} = \{1,2,\dots,N\}$ be the index set associated with the discretization points in the physical domain $\Omega \subseteq \mathbb{R}^d$.
A K-D tree is used to recursively bisect \(\mathcal{I}\) and construct the hierarchical binary tree \(\mathcal{T}\). A node $t \in \mathcal{T}$ of the tree represents a spatially localized cluster of indices. These nodes are defined such that the partition of the index set $\mathcal{I}$ reflects the geometric topology of the discretization points, ensuring that every node corresponds to a connected sub-domain in $\Omega$. 

The root node of the tree represents the entire set $\mathcal{I}$. Each non-leaf node $t \in \mathcal{T}$ is bisected into two disjoint children nodes. The recursion continues until the number of indices in a node is below or equal to a specified leaf size $m$. Let the depth of the tree be denoted by $\kappa=\floor{\log_2 \bkt{\frac{N}{m}}}$.
% In this article, we consider $m=2p$ that results in $\kappa = \floor{\log\bkt{\frac{N}{2p}}}$.
% The number of nodes at level $k$ is $2^{k}$.
% \end{definition}

% \begin{definition}[Node]
% A \textbf{node} $t \in \mathcal{T}$ of the tree represents a set of indices. These indices should correspond to a spatially localised cluster of points.
% A \textbf{node} $t \in \mathcal{T}$ is said to be \textbf{geometrically contiguous} if the indices contained within it correspond to a spatially localized cluster of points in $\Omega$. In higher dimensions ($d > 1$), this property is typically enforced via recursive coordinate bisection or space-filling curves, ensuring that sibling nodes represent adjacent but distinct sub-domains.
% \end{definition}

% \begin{definition}[HODLR Matrix]
The HODLR matrix approximates the sub-matrices corresponding to interaction between any two disjoint nodes as a low-rank matrix. At level 0 of the hierarchical tree, the HODLR matrix represented by $K^{(0)}$ is the system matrix $K$ itself. 
\begin{equation}
K^{(0)} = K.
\end{equation}
Let the nodes in level $l$ be indexed by $\cbkt{1,2,\dots,2^l}$.
Let a node $i$ at level $l$ be denoted by $i^{(l)}$.
Let the child nodes of $i^{(l)}$ be indexed by $(2i-1)^{(l+1)}$ and $(2i)^{(l+1)}$. At level 0 of the tree, there is a single node denoted by $1^{(0)}$. 
At level 1 of the tree, there are two nodes, indexed by $1^{(1)}$ and $2^{(1)}$. The HODLR matrix at level 1, represented by $K^{(1)}$, comprises of (i) the self interactions of nodes $1^{(1)}$ and $2^{(1)}$, represented by $K_{1}^{(1)}$ and $K_{2}^{(1)}$ respectively and (ii) the interactions between the nodes $1^{(1)}$ and $2^{(1)}$, represented by $K_{12}^{(1)}$ and $K_{21}^{(1)}$. 
% Let $K_{1}^{(1)}$ and $K_{2}^{(1)}$ represent the self interactions of nodes $1^{(1)}$ and $2^{(1)}$ respectively. Let $K_{12}^{(1)}$ and $K_{21}^{(1)}$ represent the interactions between nodes $1^{(1)}$ and $2^{(1)}$. 
Since nodes $1^{(1)}$ and $2^{(1)}$ are disjoint, the interactions $K_{12}^{(1)}$ and $K_{21}^{(1)}$ are approximated by low-rank matrices and are represented by $K_{12}^{(1)} = U_{1}^{(1)} V_{12}^{(1)^T}$ and $K_{21}^{(1)}=U_{2}^{(1)} V_{21}^{(1)^T}$.

\begin{equation}
K = \begin{bmatrix} 
K_{1}^{(1)} & K_{12}^{(1)} \\ 
K_{21}^{(1)} & K_{2}^{(1)} 
\end{bmatrix}\approx K^{(1)} = \begin{bmatrix} 
K_{1}^{(1)} & U_{1}^{(1)} V_{12}^{(1)^T} \\ 
U_{2}^{(1)} V_{21}^{(1)^T} & K_{2}^{(1)} 
\end{bmatrix}
\end{equation}
At level 2, the sub-matrix $K_{1}^{(1)}$ is further represented using the pairwise interactions between the child nodes of node $1^{(1)}$. The off-diagonal sub-matrices of $K_{1}^{(1)}$ are $K_{12}^{(2)}$ and $K_{21}^{(2)}$ that correspond to interactions between the sibling nodes. The sibling nodes are disjoint, and hence the interactions between them are approximated using low-rank matrices. 
The low-rank representations of the interactions between sibling nodes $i^{(l)}$ and $(i+1)^{(l)}$ are denoted by $K_{i,i+1}^{(l)}\approx U_{i}^{(1)} V_{i,i+1}^{(l)^T}$ and $K_{i+1,i}^{(l)}\approx U_{i+1}^{(1)} V_{i+1,1}^{(l)^T}$.

\begin{equation}
K_{1}^{(1)} = \begin{bmatrix} 
K_{1}^{(2)} & K_{12}^{(2)} \\ 
K_{21}^{(2)} & K_{2}^{(2)} 
\end{bmatrix}\approx
\begin{bmatrix} 
K_{1}^{(2)} & U_{1}^{(2)}V_{12}^{(2)^T} \\ 
U_{2}^{(2)}V_{21}^{(2)^T} & K_{2}^{(2)} 
\end{bmatrix}
\end{equation}
Similarly,
\begin{equation}
K_{2}^{(1)} = \begin{bmatrix} 
K_{3}^{(2)} & K_{34}^{(2)} \\ 
K_{43}^{(2)} & K_{4}^{(2)}
\end{bmatrix} \approx
\begin{bmatrix} 
K_{3}^{(2)} & U_{3}^{(2)}V_{34}^{(2)^T} \\ 
U_{4}^{(2)}V_{43}^{(2)^T} & K_{4}^{(2)}
\end{bmatrix}.
\end{equation}
The HODLR matrix at level 2 is
{\scriptsize
\begin{equation}
K^{(2)} = \begin{bmatrix} 
\begin{bmatrix} 
K_{1}^{(2)} & K_{12}^{(2)} \\ 
K_{21}^{(2)} & K_{2}^{(2)} 
\end{bmatrix} & U_{2}^{(1)} V_{12}^{(1)^T} \\ 
U_{2}^{(1)} V_{21}^{(1)^T} & \begin{bmatrix} 
K_{3}^{(2)} & K_{34}^{(2)} \\ 
K_{43}^{(2)} & K_{4}^{(2)} 
\end{bmatrix} 
\end{bmatrix} \approx 
\begin{bmatrix} 
\begin{bmatrix} 
K_{1}^{(2)} & U_{1}^{(2)} V_{12}^{(2)^T} \\ 
U_{2}^{(2)} V_{21}^{(2)^T} & K_{2}^{(2)} 
\end{bmatrix} & \hspace{-0.1cm}U_{1}^{(1)} V_{12}^{(1)^T} \\ 
U_{2}^{(1)} V_{21}^{(1)^T} & \hspace{-0.1cm}\begin{bmatrix} 
K_{3}^{(2)} & U_{3}^{(2)} V_{34}^{(1)^T} \\ 
U_{4}^{(2)} V_{43}^{(1)^T} & K_{4}^{(2)} 
\end{bmatrix} 
\end{bmatrix}
\end{equation}
}
% \end{definition}

At every level of the hierarchical tree, the off-diagonal blocks representing the interaction between sibling nodes are approximated by low-rank matrices. 
% At level $k$ of the hierarchy, the matrix is subdivided into $2^{k}$ diagonal blocks of size $N/2^{k}$. The child nodes of node $2j$ at level $k$ are nodes $2j-1$ and $2j$ at level $k+1$.
% , and each complementary off-diagonal block admits an approximation of the form $U_{j}^{(i)} {V_{j}^{(i)}}^{T}$. 
The $j^{th}$ diagonal block of the HODLR matrix at level $k$, denoted by $K_j^{(k)}$, is represented at level $k+1$ as
% Let $K_j^{(k)}$ be the $j^{\text{th}}$ diagonal block at level $k$. $K_j^{(k)}$ is approximated at level $k+1$ as:
\begin{equation}
K_{j}^{(k)} \approx
\begin{bmatrix}
    K_{2j-1}^{(k+1)} & U_{2j-1}^{(k+1)} V_{2j-1, 2j}^{(k+1)^{T}} \\
    U_{2j}^{(k+1)} V_{2j,2j-1}^{(k+1)^T} & K_{2j}^{(k+1)}
\end{bmatrix} \label{eq:hodlr}
\end{equation}
where $U_{2j-1}^{(k+1)}, V_{2j-1, 2j}^{(k+1)}, U_{2j}^{(k+1)}, V_{2j, 2j-1}^{(k+1)} \in \mathbb{R}^{\frac{N}{2^{k+1}} \times p}$, $k=0,1,\hdots,\kappa-1$ and rank $p$ follows~\cite{kandappan2023hodlr2d,gujjula2024hodlr3d,khan2024hodlrdd}
\begin{equation}
    p=\begin{cases}
        \mathcal{O}(\log(N)\log(\log(N))) & \text{in 1D}\\
        \mathcal{O}\bkt{N^{\frac{d-1}{d}}} & \text{in dD}\\
        % \mathcal{O}(N^{\frac{d-1}{d}}\log(\frac{d-1}{d\delta})) & \text{in dD}\\
    \end{cases}
\end{equation}
% $p \ll N$ is small and typically constant for 1D problems, or poly-logarithmic in $N$ for higher dimensions.

% where $\operatorname{rank}(U_{2j-1}^{(i+1)} V_{2j-1}^{(i+1)T}) = \operatorname{rank}(U_{2j}^{(i+1)} V_{2j}^{(i+1)T}) = p \ll N/2^{i+1}$, and \(U_{j}^{(i+1)}, V_{j}^{(i+1)} \in \mathbb{R}^{(N/2^{i+1}) \times p}\) are tall, skinny matrices containing compressed representations of long-range interactions between the corresponding subclusters. The diagonal blocks $K_{2j-1}^{(i+1)}, K_{2j}^{(i+1)}$ are full-rank and are stored explicitly. A $2$-level HODLR matrix $K$ can be written as shown in Eq.~\eqref{eq:hodlr2level}.

% \begin{equation}
%     K =
%     \begin{bmatrix}
%         K_{1}^{(1)} & U_{1}^{(1)}{V_{1}^{(1)}}^{T} \\
%         U_{2}^{(1)}{V_{2}^{(1)}}^{T} & K_{2}^{(1)}
%     \end{bmatrix}
%     =
%     \begin{bmatrix}
%         \begin{bmatrix}
%             K_{1}^{(2)} & U_{1}^{(2)}{V_{1}^{(2)}}^{T} \\
%             U_{2}^{(2)}{V_{2}^{(2)}}^{T} & K_{2}^{(2)}
%         \end{bmatrix} & U_{1}^{(1)}{V_{1}^{(1)}}^{T} \\
%         U_{2}^{(1)}{V_{2}^{(1)}}^{T} & \begin{bmatrix}
%             K_{3}^{(2)} & U_{3}^{(2)}{V_{3}^{(2)}}^{T} \\
%             U_{4}^{(2)}{V_{4}^{(2)}}^{T} & K_{4}^{(2)}
%         \end{bmatrix}
%     \end{bmatrix}
%     \label{eq:hodlr2level}
% \end{equation}

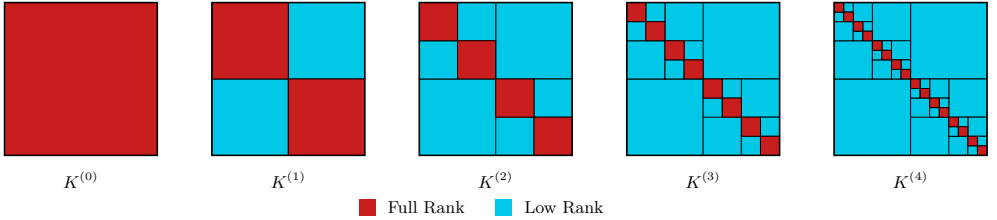
\begin{figure}[htbp]
    \centering
    \begin{tikzpicture}[scale=1]

\def\S{2.8}   % side length of each matrix square
\def\gap{1.0} % gap between matrices

% =================================================================
%  Draw a single filled block with a thin border.
% =================================================================
\newcommand{\blk}[5]{%
  \fill[#5] ({#1*\c},{-#2*\c}) rectangle ++({#3*\c},{-#4*\c});
  \draw[black,very thin] ({#1*\c},{-#2*\c}) rectangle ++({#3*\c},{-#4*\c});
}

% =================================================================
%   κ = 0   (Entire matrix is a single full-rank block)
% =================================================================
\begin{scope}[shift={(0,0)}]
  \pgfmathsetmacro{\c}{\S}
  % Single full rank block representing the root
  \blk{0}{0}{1}{1}{fullrank}
  % Outer border
  \draw[thick] (0,0) rectangle ++({\S},{-\S});
  \node[below=5pt,font=\small] at ({\S/2},{-\S}) {$K^{(0)}$};
\end{scope}

% =================================================================
%   κ = 1   (2 cells)
% =================================================================
\pgfmathsetmacro{\xA}{\S+\gap}
\begin{scope}[shift={(\xA,0)}]
  \pgfmathsetmacro{\c}{\S/2}
  \blk{1}{0}{1}{1}{lowrank}
  \blk{0}{1}{1}{1}{lowrank}
  \blk{0}{0}{1}{1}{fullrank}
  \blk{1}{1}{1}{1}{fullrank}
  \draw[thick] (0,0) rectangle ++({\S},{-\S});
  \node[below=5pt,font=\small] at ({\S/2},{-\S}) {$K^{(1)}$};
\end{scope}

% =================================================================
%   κ = 2   (4 cells)
% =================================================================
\pgfmathsetmacro{\xB}{2*(\S+\gap)}
\begin{scope}[shift={(\xB,0)}]
  \pgfmathsetmacro{\c}{\S/4}
  \blk{2}{0}{2}{2}{lowrank}
  \blk{0}{2}{2}{2}{lowrank}
  \blk{1}{0}{1}{1}{lowrank}
  \blk{0}{1}{1}{1}{lowrank}
  \blk{3}{2}{1}{1}{lowrank}
  \blk{2}{3}{1}{1}{lowrank}
  \foreach \i in {0,...,3}{ \blk{\i}{\i}{1}{1}{fullrank} }
  \draw[thick] (0,0) rectangle ++({\S},{-\S});
  \node[below=5pt,font=\small] at ({\S/2},{-\S}) {$K^{(2)}$};
\end{scope}

% =================================================================
%   κ = 3   (8 cells)
% =================================================================
\pgfmathsetmacro{\xC}{3*(\S+\gap)}
\begin{scope}[shift={(\xC,0)}]
  \pgfmathsetmacro{\c}{\S/8}
  \blk{4}{0}{4}{4}{lowrank}
  \blk{0}{4}{4}{4}{lowrank}
  \blk{2}{0}{2}{2}{lowrank}
  \blk{0}{2}{2}{2}{lowrank}
  \blk{6}{4}{2}{2}{lowrank}
  \blk{4}{6}{2}{2}{lowrank}
  \foreach \g in {0,2,4,6}{
    \pgfmathtruncatemacro{\gp}{\g+1}
    \blk{\gp}{\g}{1}{1}{lowrank}
    \blk{\g}{\gp}{1}{1}{lowrank}
  }
  \foreach \i in {0,...,7}{ \blk{\i}{\i}{1}{1}{fullrank} }
  \draw[thick] (0,0) rectangle ++({\S},{-\S});
  \node[below=5pt,font=\small] at ({\S/2},{-\S}) {$K^{(3)}$};
\end{scope}

% =================================================================
%   κ = 4   (16 cells)
% =================================================================
\pgfmathsetmacro{\xD}{4*(\S+\gap)}
\begin{scope}[shift={(\xD,0)}]
  \pgfmathsetmacro{\c}{\S/16}
  \blk{8}{0}{8}{8}{lowrank}
  \blk{0}{8}{8}{8}{lowrank}
  \foreach \g in {0,8}{
    \pgfmathtruncatemacro{\gp}{\g+4}
    \blk{\gp}{\g}{4}{4}{lowrank}
    \blk{\g}{\gp}{4}{4}{lowrank}
  }
  \foreach \g in {0,4,8,12}{
    \pgfmathtruncatemacro{\gp}{\g+2}
    \blk{\gp}{\g}{2}{2}{lowrank}
    \blk{\g}{\gp}{2}{2}{lowrank}
  }
  \foreach \g in {0,2,...,14}{
    \pgfmathtruncatemacro{\gp}{\g+1}
    \blk{\gp}{\g}{1}{1}{lowrank}
    \blk{\g}{\gp}{1}{1}{lowrank}
  }
  \foreach \i in {0,...,15}{ \blk{\i}{\i}{1}{1}{fullrank} }
  \draw[thick] (0,0) rectangle ++({\S},{-\S});
  \node[below=5pt,font=\small] at ({\S/2},{-\S}) {$K^{(4)}$};
\end{scope}

% =================================================================
%  LEGEND
% =================================================================
\pgfmathsetmacro{\legY}{-\S - 1.1}
\pgfmathsetmacro{\totalW}{\xD + \S}
\pgfmathsetmacro{\legMid}{\totalW/2}

\pgfmathsetmacro{\lA}{\legMid - 2.5}
\fill[fullrank] (\lA,\legY) rectangle ++(0.3,0.3);
\node[right,font=\small] at ({\lA+0.4},{\legY+0.15}) {Full Rank};

\pgfmathsetmacro{\lB}{\lA + 2.5}
\fill[lowrank] (\lB,\legY) rectangle ++(0.3,0.3);
\node[right,font=\small] at ({\lB+0.4},{\legY+0.15}) {Low Rank};

\end{tikzpicture}
    \caption{HODLR matrix at successive levels}
    \label{fig:hodlr}
\end{figure}

The hierarchical matrix is constructed recursively down to the leaf level, where the diagonal blocks are stored explicitly. At intermediate levels, only the low-rank factors are stored, yielding a data-sparse representation with storage complexity \(O(pN \log N)\). Figure \ref{fig:hodlr} illustrates the low-rank structure of the HODLR matrix at successive levels. 

% Unlike nested-basis formats such as HSS or p-HSS, HODLR matrices do not require compatibility constraints across levels; this simplicity makes them attractive for general kernel matrices and operators lacking strict nestedness.

\subsection{Fast direct solver}\label{sec:directsolveralg}
In this section, we discuss the fast direct solver for linear systems involving HODLR matrices proposed in~\cite{ambikasaran2013fast}, which forms the foundation of the neural network proposed in this article. The solver involves two phases i) the factorization phase and ii) the solve phase. 
Before describing the two phases, we describe the Sherman-Morrison-Woodbury (SMW) formula, which is the key tool underlying the fast direct solver.
\subsubsection{Sherman--Morrison--Woodbury (SMW) formula}
SMW~\cite{woodbury1950, hager1989updating} \\formula gives a computationally efficient way of solving a linear system of equations characterized by a low-rank perturbation of the identity matrix. Given the system
\begin{equation}
\bkt{I + UV^T}x = b, \label{eq:smw1}
\end{equation}
where $I \in \mathbb{R}^{N\times N}$, $U,V \in \mathbb{R}^{N\times p}$, and $p \ll N$, the SMW formula provides the solution
\begin{equation}
x = b - U\bkt{I+V^TU}^{-1}V^Tb, \label{eq:smw2}
\end{equation}
where the $p \times p$ matrix $\bkt{I + V^TU}$ only needs to be inverted. The cost of the naive dense solve of Equation~\eqref{eq:smw1} is $O\bkt{N^3}$, whereas the SMW formula (Equation~\eqref{eq:smw2}) reduces it to $O(Np)$. We suggest the reader to refer \cite{woodbury1950, hager1989updating} for the derivation.

\subsubsection{Factorization phase}\label{sssec:factorization}
The objective of the factorization phase is to factorize the $\kappa$-level HODLR matrix $K^{(\kappa)}$ in the form
\begin{equation}
    K^{(\kappa)} = K_{\kappa} K_{\kappa-1} \cdots K_{1} K_{0}. \label{eq:telescopic_factorization}
\end{equation}
The factorization begins at the leaf level. 
% $K^{\kappa}$ consists of the diagonal blocks of the matrix dense blocks which are the self-interactions of the leaf nodes. 
The first factor $K_\kappa$ is a block-diagonal matrix comprising the diagonal blocks of $K^{(\kappa)}$, which is 
\begin{equation}\label{eq:K_kappa}
K_\kappa = \begin{bmatrix} 
K_1^{(\kappa)}\\ & K_2^{(\kappa)}\\ & & \ddots & \\ & & & K_{2^\kappa}^{(\kappa)} \end{bmatrix}
\end{equation}
By factoring out $K_\kappa$ from $K^{(\kappa)}$, we have 

\begin{equation}
\begin{aligned}
% &K_{\kappa-1} \cdots K_{1} K_{0} \\
&K_\kappa^{-1}K^{(\kappa)}\\
&=\resizebox{0.9\textwidth}{!}{%
$ \displaystyle % Keeps the math in display style
\begin{bmatrix} 
K_1^{(\kappa)^{-1}}\\ & K_2^{(\kappa)^{-1}}\\ & & \ddots & \\ & & & K_{2^\kappa}^{(\kappa)^{-1}} 
\end{bmatrix}
\begin{bmatrix} 
\begin{bmatrix}
K_1^{(\kappa)} & U_{1}^{(\kappa)}V_{12}^{(\kappa)^T}\\
U_{2}^{(\kappa)}V_{21}^{(\kappa)T} & K_2^{(\kappa)}
\end{bmatrix} & U_{1}^{(\kappa-1)}V_{12}^{(\kappa-1)^T}
& \dots & \dots \\ 
U_{2}^{(\kappa-1)}V_{21}^{(\kappa-1)^T} & 
\begin{bmatrix}
K_3^{(\kappa)} & U_{3}^{(\kappa)}V_{34}^{(\kappa)^T}\\
U_{4}^{(\kappa)}V_{43}^{(\kappa)T} & K_4^{(\kappa)}
\end{bmatrix} & \dots & \dots \\ 
\vdots & \vdots & \ddots & \vdots \\ 
\dots & \dots & \dots & \begin{bmatrix}
K_{2^{\kappa}-1}^{(\kappa)} & U_{2^{\kappa}-1}^{(\kappa)}V_{2^{\kappa}-1,2^{\kappa}}^{(\kappa)^T}\\
U_{2^{\kappa}}^{(\kappa)}V_{2^{k},2^{k}-1}^{(\kappa)^T} & K_{2^{\kappa}}^{(\kappa)}
\end{bmatrix}
\end{bmatrix}
$
}\\
&= 
\resizebox{0.9\textwidth}{!}{%
$ \displaystyle % Keeps the math in display style
\begin{bmatrix} 
\begin{bmatrix}
I & U_1^{(\kappa,1)}V_{12}^{(\kappa)^T}\\
U_2^{(\kappa,1)}V_{21}^{(\kappa)^T} & I
\end{bmatrix} & U_{1}^{(\kappa-1,1)}V_{12}^{(\kappa-1)^T}
& \dots & \dots \\ 
U_{2}^{(\kappa-1,1)}V_{21}^{(\kappa-1)^T} & 
\begin{bmatrix}
I & U_{3}^{(\kappa,1)}V_{34}^{(\kappa)^T}\\
U_{4}^{(\kappa,1)}V_{43}^{(\kappa)^T} & I
\end{bmatrix} & \dots & \dots \\ 
\vdots & \vdots & \ddots & \vdots \\ 
\dots & \dots & \dots & 
\begin{bmatrix}
I & U_{2^{\kappa}-1}^{(\kappa,1)}V_{2^{\kappa}-1,2^{\kappa}}^{(\kappa)^T}\\
U_{2^{\kappa}}^{(\kappa,1)}V_{2^{\kappa},2^{\kappa}-1}^{(\kappa)^T} & I
\end{bmatrix}
\end{bmatrix}
$
}.
\end{aligned}
\end{equation}
$U_{i}^{(k,1)}$, for all $k\in\{1,2,\dots\kappa\}$, indicates that $U_{i}^{(k)}$ has been updated as a result of applying the inverse of $K_i^{(\kappa)}$, for all $i\in\{1,2,\dots2^{\kappa}\}$.
The size of $K_i^{(\kappa)}$ is $\mathcal{O}(m)\times \mathcal{O}(m)$. Let $m=p$.
% , since $\kappa = \mathcal{O}\left(\log\left(\frac{N}{2p}\right)\right)$.
The cost of applying the inverse of $K_i^{(\kappa)}$ to $r$ columns is $\mathcal{O}\left(p^{3}+p^{2}r\right)$. The inverse of $K_i^{(\kappa)}$ has to be applied to the appropriate rows of $\cbkt{U^{(k)}_{j}}_{j^{(k)}\in\mathcal{A}\bkt{i^{(\kappa)}}}$, where $\mathcal{A}\bkt{i^{(\kappa)}}$ represents the set of ancestors of $i^{(\kappa)}$ including $i^{(\kappa)}$. The cardinality of the set $\mathcal{A}\bkt{i^{(\kappa)}}$ is $\kappa$ and each $U^{(k)}_{j}$ has $p$ columns. Therefore, the number of columns to which the inverse of $K_i^{(\kappa)}$ has to be applied is $p\kappa$. 
% The inverse of $K_i^{(\kappa)}$ has to be applied to the corresponding rows of $U_{i}^{(k,1)}$ for all $k=\{1,2,\dots,\kappa\}$. 
Hence, the cost of applying the inverse of $K_i^{(\kappa)}$ to the $i^{th}$ block row is $\mathcal{O}\bkt{p^{3}+p^{2}p\kappa}=\mathcal{O}(p^{3}\kappa)$. 
% To factor out $K_{\kappa}$, the inverse of $K_{i}^{(\kappa,1)}$ has to be applied to the $i^{th}$ block row of $K^{(\kappa)}$ for all $i\in{0,1,\dots,2^{\kappa}}$. 
To factor out $K_{\kappa}$, the inverse of $K_i^{(\kappa)}$ is to be applied for all $i\in\{1,2,\dots,2^{\kappa}\}$.
Therefore, the total cost of factoring out $K_{\kappa}$ is $\mathcal{O}\bkt{p^{3}\kappa 2^{\kappa}}$. Since $\kappa = \floor{\log_2\bkt{\frac{N}{p}}}$, the total cost of factoring out $K_{\kappa}$ is $\mathcal{O}\bkt{p^{2}N\kappa}$.

% Upon factoring out $K_{\kappa}$ from $K^{(k)}$, $U_{i}^{(\kappa)}$ gets updated to $U_{i}^{(\kappa,1)}$. 
% Upon factoring out $K_{\kappa}$ from $K^{(\kappa)}$, $U_{2j-1}^{(k)}$ gets updated to $U_{2j-1}^{(k,1)}$, $\forall k\in\{1,2,\dots,\kappa\}$. 
Upon factoring out $K_{\kappa}$ from $K^{(k)}$, we obtain the matrix
\begin{equation}
K_\kappa^{-1}K^{(\kappa)}=
    \begin{bmatrix} 
        K_{1}^{(\kappa-1,1)} & U_{1}^{(\kappa-1,1)}V_{12}^{(\kappa-1)^T}
        & \dots & \dots \\ 
        U_{2}^{(\kappa-1,1)}V_{21}^{(\kappa-1)^T} & 
        K_{2}^{(\kappa-1,1)} & \dots & \dots \\ 
        \vdots & \vdots & \ddots & \vdots \\ 
        \dots & \dots & \dots & 
        K_{2^{\kappa-1}}^{(\kappa-1,1)}
    \end{bmatrix} 
\end{equation}
where 
\begin{equation}\label{eq:K_is_I_plus_UV}
    K_{i}^{(\kappa-1,1)}=
    \begin{bmatrix}
        I & U_{2i-1}^{(\kappa,1)}V_{2i-1, 2i}^{(\kappa)^T}\\
        U_{2i}^{(\kappa,1)}V_{2i, 2i-1}^{(\kappa)^T} & I
    \end{bmatrix},\; \forall\; i\in\cbkt{1,2,\dots,2^{\kappa-1}}.
\end{equation}
The second factor $K_{\kappa-1}$ in the factorization of Equation~\eqref{eq:telescopic_factorization},
\begin{equation}\label{eq:K_kappa_minus_1}
K_{\kappa-1} = \begin{bmatrix} 
K_{1}^{(\kappa-1,1)}\\ & K_{2}^{(\kappa-1,1)}\\ & & \ddots & \\ & & & K_{2^{\kappa-1}}^{(\kappa-1,1)} \end{bmatrix},
\end{equation}
is a block diagonal matrix comprising the diagonal blocks of $K_\kappa^{-1}K^{(\kappa)}$.
Upon factoring out $K_{\kappa-1}$ from $K_\kappa^{-1}K^{(\kappa)}$, we obtain the matrix
\begin{equation}
\begin{aligned}
% &K_{\kappa-1} \cdots K_{1} K_{0} \\
&K_{\kappa-1}^{-1}K_\kappa^{-1}K^{(\kappa)}\\
&=
\resizebox{0.9\textwidth}{!}{%
$ \displaystyle % Keeps the math in display style
\begin{bmatrix} 
K_{1}^{(\kappa-1,1)^{-1}}\\ & K_{2}^{(\kappa-1,1)^{-1}}\\ & & \ddots & \\ & & & K_{2^{\kappa-1}}^{(\kappa-1,1)^{-1}} 
\end{bmatrix}
\begin{bmatrix} 
K_{1}^{(\kappa-1,1)} & U_{1}^{(\kappa-1,1)}V_{12}^{(\kappa-1)^T}
& \dots & \dots \\ 
U_{2}^{(\kappa-1,1)}V_{21}^{(\kappa-1)^T} & 
K_{2}^{(\kappa-1,1)} & \dots & \dots \\ 
\vdots & \vdots & \ddots & \vdots \\ 
\dots & \dots & \dots & 
K_{2^{\kappa-1}}^{(\kappa-1,1)}
\end{bmatrix}
$
}\\
\end{aligned}
\end{equation}
From Equation~\eqref{eq:K_is_I_plus_UV}, it can be observed that $K_{i}^{(\kappa-1,1)}$ can be written as $K_{i}^{(\kappa-1,1)} = I+\widetilde{U}_{i}^{(\kappa)}\widetilde{V}_{i}^{(\kappa)^T}$ where
\begin{equation}\label{eq:U_and_V_structure}
    \widetilde{U}_{i}^{(\kappa)} = 
    \begin{bmatrix}
        U_{2i-1}^{(\kappa,1)} & 0\\
        0 & U_{2i}^{(\kappa,1)}
    \end{bmatrix} 
    \text{ and }
    \widetilde{V}_{i}^{(\kappa)^{T}} = 
    \begin{bmatrix}
        0 & V_{2i, 2i-1}^{(\kappa)^{T}}\\
        V_{2i-1, 2i}^{(\kappa)^{T}} & 0
    \end{bmatrix}. 
\end{equation}
$I+\widetilde{U}_{i}^{(\kappa)}\widetilde{V}_{i}^{(\kappa)^T}$ is a rank $2p$ perturbation to an identity matrix of size $\frac{N}{2^{\kappa-1}}\times \frac{N}{2^{\kappa-1}}$.
The inverse of $K_{i}^{(\kappa-1,1)}$ using the SMW formula is
\begin{equation}
    K_{i}^{(\kappa-1,1)^{-1}}=\left(I+\widetilde{U}_{i}^{(\kappa)}\widetilde{V}_{i}^{(\kappa)^T}\right)^{-1} = I - \widetilde U_{j}^{(\kappa)} \Big(I + {\widetilde{V}_{j}^{(\kappa)^T}}\widetilde U_{j}^{(\kappa)}\Big)^{-1} {\widetilde{V}_{j}^{(\kappa)^T}}    .
\end{equation}
Using the SMW formula, the cost to apply the inverse of $I+\widetilde{U}_{i}^{(\kappa)}\widetilde{V}_{i}^{(\kappa)^T}$ to $r$ columns is $\mathcal{O}\left(\frac{N}{2^{\kappa-1}}\times 2p\times r\right)$. The inverse of $K_{i}^{(\kappa-1,1)}$ has to be applied to the appropriate rows of $\cbkt{U^{(k,1)}_{j}}_{j^{(k)}\in\mathcal{A}\bkt{i^{(\kappa-1)}}}$, where $\mathcal{A}\bkt{i^{(\kappa-1)}}$ represents the set of ancestors of $i^{(\kappa-1)}$. The cardinality of the set $\mathcal{A}\bkt{i^{(\kappa-1)}}$ is $\kappa-1$ and each $U^{(k,1)}_{j}$ has $2p$ columns. Therefore, the number of columns to which the inverse of $K_{i}^{(\kappa-1,1)}$ has to be applied is $2p(\kappa-1)$, and the cost is $\mathcal{O}\bkt{\frac{N}{2^{\kappa-1}}\times 2p\times 2p(\kappa-1)}=\mathcal{O}\bkt{\frac{p^{2}N(\kappa-1)}{2^{\kappa-1}}}$. 
To factor out $K_{\kappa-1}$, the inverse of $K_i^{(\kappa-1,1)}$ is to be applied for all $i\in\{1,\dots,2^{\kappa-1}\}$.
Therefore, the total cost of factoring out $K_{\kappa-1}$ is $\mathcal{O}\bkt{p^{2}N(\kappa-1)}$. 

\begin{figure}[htbp]
    \centering
    \begin{tikzpicture}[scale=1]

\def\cs{0.36}  % cell size
\def\N{8}      % 8×8 grid

% =================================================================
%  \blk{Xorigin}{col}{row}{w}{h}{colour}
% =================================================================
\newcommand{\blk}[6]{%
  \fill[#6] ({#1+#2*\cs},{-#3*\cs}) rectangle ++({#4*\cs},{-#5*\cs});
  \draw[black,very thin] ({#1+#2*\cs},{-#3*\cs}) rectangle ++({#4*\cs},{-#5*\cs});
}

% =================================================================
%  \identblk{Xorigin}{col}{row}{size}  — blue diagonal line
% =================================================================
\newcommand{\identblk}[4]{%
  \fill[white] ({#1+#2*\cs},{-#3*\cs}) rectangle ++({#4*\cs},{-#4*\cs});
  \draw[black,very thin] ({#1+#2*\cs},{-#3*\cs}) rectangle ++({#4*\cs},{-#4*\cs});
  \draw[identblue,thin] ({#1+#2*\cs},{-#3*\cs}) -- ++({#4*\cs},{-#4*\cs});
}

% ---------- Horizontal offsets ----------
% K (full red) ≈ K^{(3)} = K_3 × K_2 × K_1 × K_0
\def\matK{0}                                           % K (full red)
\pgfmathsetmacro{\approxX}{\N*\cs+0.35}               % ≈ sign
\pgfmathsetmacro{\matA}{\approxX+0.55}                 % K^{(3)}
\pgfmathsetmacro{\eqX}{\matA+\N*\cs+0.35}             % = sign
\pgfmathsetmacro{\matB}{\eqX+0.55}                     % K_3
\pgfmathsetmacro{\txB}{\matB+\N*\cs+0.2}              % × sign
\pgfmathsetmacro{\matC}{\txB+0.35}                     % K_2
\pgfmathsetmacro{\txC}{\matC+\N*\cs+0.2}
\pgfmathsetmacro{\matD}{\txC+0.35}                     % K_1
\pgfmathsetmacro{\txD}{\matD+\N*\cs+0.2}
\pgfmathsetmacro{\matE}{\txD+0.35}                     % K_0

% ================================================================
%  MATRIX K  (full red — the original dense matrix)
% ================================================================
\blk{\matK}{0}{0}{\N}{\N}{fullrank}
\draw[thick] ({\matK},0) rectangle ++({\N*\cs},{-\N*\cs});
\node[below=5pt,font=\small] at ({\matK+\N*\cs/2},{-\N*\cs}) {$K$};

% ≈ sign
\node[font=\Large] at ({\approxX},{-\N*\cs/2}) {$\approx$};

% ================================================================
%  MATRIX K^{(3)}  (HODLR κ=3, same style as hodlr.tex)
% ================================================================
% Level 0 off-diagonal (4×4 solid blocks)
\blk{\matA}{4}{0}{4}{4}{lowrank}
\blk{\matA}{0}{4}{4}{4}{lowrank}
% Level 1 off-diagonal (2×2 solid blocks)
\blk{\matA}{2}{0}{2}{2}{lowrank}
\blk{\matA}{0}{2}{2}{2}{lowrank}
\blk{\matA}{6}{4}{2}{2}{lowrank}
\blk{\matA}{4}{6}{2}{2}{lowrank}
% Level 2 off-diagonal (1×1 solid blocks)
\foreach \g in {0,2,4,6}{
  \pgfmathtruncatemacro{\gp}{\g+1}
  \blk{\matA}{\gp}{\g}{1}{1}{lowrank}
  \blk{\matA}{\g}{\gp}{1}{1}{lowrank}
}
% Leaf diagonal (8 full-rank blocks)
\foreach \i in {0,...,7}{
  \blk{\matA}{\i}{\i}{1}{1}{fullrank}
}
\draw[thick] ({\matA},0) rectangle ++({\N*\cs},{-\N*\cs});
\node[below=5pt,font=\small] at ({\matA+\N*\cs/2},{-\N*\cs}) {$K^{(3)}$};

% = sign
\node[font=\Large] at ({\eqX},{-\N*\cs/2}) {$=$};

% ================================================================
%  MATRIX K_3  (leaf-level: 8 diagonal full-rank blocks)
% ================================================================
\foreach \i in {0,...,7}{
  \blk{\matB}{\i}{\i}{1}{1}{fullrank}
}
\draw[thick] ({\matB},0) rectangle ++({\N*\cs},{-\N*\cs});
\node[below=5pt,font=\small] at ({\matB+\N*\cs/2},{-\N*\cs}) {$K_3$};

% × sign
\node[font=\Large] at ({\txB},{-\N*\cs/2}) {$\times$};

% ================================================================
%  MATRIX K_2  (4 groups of 2×2)
% ================================================================
\foreach \g in {0,...,3}{
  \pgfmathtruncatemacro{\base}{2*\g}
  \pgfmathtruncatemacro{\basep}{\base+1}
  \identblk{\matC}{\base}{\base}{1}
  \identblk{\matC}{\basep}{\basep}{1}
  \blk{\matC}{\basep}{\base}{1}{1}{lowrank}
  \blk{\matC}{\base}{\basep}{1}{1}{lowrank}
}
\draw[thick] ({\matC},0) rectangle ++({\N*\cs},{-\N*\cs});
\node[below=5pt,font=\small] at ({\matC+\N*\cs/2},{-\N*\cs}) {$K_2$};

% × sign
\node[font=\Large] at ({\txC},{-\N*\cs/2}) {$\times$};

% ================================================================
%  MATRIX K_1  (2 groups of 4×4)
% ================================================================
\identblk{\matD}{0}{0}{2}
\identblk{\matD}{2}{2}{2}
\blk{\matD}{2}{0}{2}{2}{lowrank}
\blk{\matD}{0}{2}{2}{2}{lowrank}
\identblk{\matD}{4}{4}{2}
\identblk{\matD}{6}{6}{2}
\blk{\matD}{6}{4}{2}{2}{lowrank}
\blk{\matD}{4}{6}{2}{2}{lowrank}
\draw[thick] ({\matD},0) rectangle ++({\N*\cs},{-\N*\cs});
\node[below=5pt,font=\small] at ({\matD+\N*\cs/2},{-\N*\cs}) {$K_1$};

% × sign
\node[font=\Large] at ({\txD},{-\N*\cs/2}) {$\times$};

% ================================================================
%  MATRIX K_0  (1 group of 8×8)
% ================================================================
\identblk{\matE}{0}{0}{4}
\identblk{\matE}{4}{4}{4}
\blk{\matE}{4}{0}{4}{4}{lowrank}
\blk{\matE}{0}{4}{4}{4}{lowrank}
\draw[thick] ({\matE},0) rectangle ++({\N*\cs},{-\N*\cs});
\node[below=5pt,font=\small] at ({\matE+\N*\cs/2},{-\N*\cs}) {$K_0$};

% ================================================================
%  LEGEND
% ================================================================
\pgfmathsetmacro{\legY}{-\N*\cs - 1.0}
\pgfmathsetmacro{\legMid}{(\matE+\N*\cs)/2}

\pgfmathsetmacro{\lA}{\legMid - 5.8}
\fill[fullrank] (\lA,\legY) rectangle ++(0.3,0.3);
\node[right,font=\small] at ({\lA+0.4},{\legY+0.15}) {Full Rank};

\pgfmathsetmacro{\lB}{\lA + 2.5}
\fill[lowrank] (\lB,\legY) rectangle ++(0.3,0.3);
\node[right,font=\small] at ({\lB+0.4},{\legY+0.15}) {Low Rank};

\pgfmathsetmacro{\lC}{\lB + 2.5}
\draw[thick] (\lC,\legY) rectangle ++(0.3,0.3);
\node[right,font=\small] at ({\lC+0.4},{\legY+0.15}) {Zero Matrix};

\pgfmathsetmacro{\lD}{\lC + 2.8}
\draw[thick] (\lD,\legY) rectangle ++(0.3,0.3);
\draw[identblue,thin] (\lD,{\legY+0.3}) -- ++(0.3,-0.3);
\node[right,font=\small] at ({\lD+0.4},{\legY+0.15}) {Identity Matrix};

\end{tikzpicture}
    \caption{A $3$-level HODLR matrix factorization}
    \label{fig:directfactorization}
\end{figure}
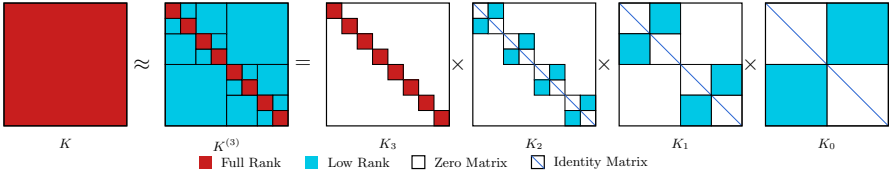

% $U_{i}^{(k,1)}$ for all $k=\{1,2,\dots,\kappa-1\}$. Since $U_{i}^{(\kappa,1)}$ is of size $\frac{N}{2^{\kappa}}\times 2p$, the number of columns $r$ to which the inverse of $I+\widetilde{U}_{i}^{(\kappa)}\widetilde{V}_{i}^{(\kappa)^T}$ is to be applied is $2p$. Therefore the total cost of applying the inverse of $I+\widetilde{U}_{i}^{(\kappa)}\widetilde{V}_{i}^{(\kappa)^T}$ to the $i^{th}$ block row is $\mathcal{O}\left(\frac{N}{2^{\kappa-1}}\times 2p\times 2p \times (\kappa-1)\right)=\mathcal{O}(\frac{p^{2}N(\kappa-1)}{2^{\kappa-1}})$. To factor out $K_{\kappa-1}$, the inverse of $K_{i}^{(\kappa-1,1)} = I+\widetilde{U}_{i}^{(\kappa)}\widetilde{V}_{i}^{(\kappa)^T}$ has to be applied to the $i^{th}$ block row of $K_\kappa^{-1}K^{(\kappa)}$ for all $i\in{0,1,\dots,2^{\kappa-1}}$. Therefore the total cost of factoring out $K_{\kappa-1}$ is $\mathcal{O}(\frac{p^{2}N(\kappa-1)}{2^{\kappa-1}}(2^{\kappa-1})) = \mathcal{O}(p^{2}N(\kappa-1))$. 
This is repeated until level 0 is reached, to obtain the subsequent factors \\$\cbkt{K_{\kappa-2},\hdots,K_0}$. 
% \begin{equation}
%     K_{0}^{-1} K_{1}^{-1}\dots K_{\kappa-1}^{-1}K_\kappa^{-1}K^{(\kappa)} = I
% \end{equation}
As a result, this process yields the factorization
\begin{equation}\label{eq:HODLRfactorization}
    K^{(\kappa)} = K_\kappa K_{\kappa-1} \dots K_{1} K_{0}.
\end{equation}
The structure of the factorization is illustrated in Figure~\ref{fig:directfactorization}.

On similar lines of deriving the cost of factoring out $K_{\kappa-1}$, the cost of factoring out $K_{k}$ for a $k\in\{1,2\dots \kappa-2\}$ is $\mathcal{O}(p^{2}Nk)$. Hence, the total cost of factorization is equal to the sum of costs of factoring out $K_{k}$ for all $k\in\{1,2\dots \kappa\}$, which is $\mathcal{O}\bkt{\sum_{k=1}^{\kappa}p^{2}Nk} = \mathcal{O}\bkt{p^{2}N\log^2(N)}$.

\subsubsection{Solve phase} \label{sec:solve_phase}
% Once the factorization of the HODLR matrix $K$ is complete, solving for $K^{(\kappa)}x=b$ becomes straightforward. 
% gives
% \begin{equation*}
%     K^{(\kappa)} = K_\kappa K_{\kappa-1} \dots K_{1} K_{0}.
% \end{equation*}
Using the factorization of $K^{(\kappa)}$ obtained in Equation~\eqref{eq:HODLRfactorization}, the solution to $K^{(\kappa)}x=b$ is $x = K^{(\kappa)^{-1}}b = K_{0}^{-1} K_{1}^{-1}\dots K_{\kappa-1}^{-1}K_\kappa^{-1}b$.
To solve for $x$, each of the $K_{k}^{-1}$ has to be applied to vector $b$ starting from $k=\kappa$ to $k=0$. Applying $K_k^{-1}$ to a vector is similar to applying $K_k^{-1}$ in the factorization phase (while factoring them out).

\textbf{Structure of $K_{\kappa}^{-1}$.} From Equation~\eqref{eq:K_kappa}, it follows that $K_{\kappa}$ is a block diagonal matrix with diagonal blocks $K_{i}^{(\kappa)}$, for all $i\in\{1,2,\dots 2^\kappa\}$. Thus $K_{\kappa}^{-1}$ is also a block-diagonal matrix with diagonal blocks $K_{i}^{(\kappa)^{-1}}$, for all $i\in\{1,2,\dots 2^\kappa\}$. 
% The diagonal blocks of $K_{\kappa}^{-1}$ are denoted by $K_i^{\kappa^{-1}}$, $i\in\{1,2,\dots,2^{\kappa}\}$, and each of them are dense full-rank matrices.

\textbf{Structure of $K_{k}^{-1}$ for all $k\in\{\kappa-1,\hdots,0\}$.} From Section~\ref{sssec:factorization}, it follows that $K_{k}$ is a block diagonal matrix with diagonal blocks 
\begin{equation}\label{eq:K_block_diag_at_level_k}
    K_{i}^{(k,\kappa-k)} = I+\widetilde{U}_{i}^{(k+1)}\widetilde{V}_{i}^{(k+1)^T},\;\; \forall\; i\in\{1,2,\dots 2^k\}.
\end{equation}
% $K_{i}^{(k,\kappa-k)}$, for all $i\in\{1,2,\dots 2^k\}$. 
The notation $K_{i}^{(k,\kappa-k)}$ indicates that it is the $i^{th}$ diagonal block of matrix $K_{k}$, 
\begin{equation}\label{eq:U_and_V_structure_at_level_k}
    \widetilde{U}_{i}^{(k+1)} = 
    \begin{bmatrix}
        U_{2i-1}^{(k+1,\kappa-k))} & 0\\
        0 & U_{2i}^{(k+1,\kappa-k))}
    \end{bmatrix} 
    \text{ and }
    \widetilde{V}_{i}^{(k+1)^{T}} = 
    \begin{bmatrix}
        0 & V_{2i, 2i-1}^{(k+1)^{T}}\\
        V_{2i-1, 2i}^{(k+1)^{T}} & 0
    \end{bmatrix}. 
\end{equation}
and $U_{2i-1}^{(k+1)}$ and $U_{2i}^{(k+1)}$ had undergone updation $\kappa-k$ times as a result of factoring out $K_{\kappa},K_{\kappa},\hdots,\text{ and } K_{k+1}$. It follows that $K_{k}^{-1}$ is a block-diagonal matrix with diagonal blocks $K_{i}^{(k,\kappa-k)^{-1}}$. Using the SMW formula, it is
\begin{equation}
\underbrace{K_{i}^{(k,\kappa-k)^{-1}}}_{\frac{N}{2^{k}}\times \frac{N}{2^{k}}}
=
\underbrace{I}_{\frac{N}{2^{k}}\times \frac{N}{2^{k}}}
-
\underbrace{\widetilde{U}_{i}^{(k+1)}}_{\frac{N}{2^{k}}\times 2p}
\;
\underbrace{S_{i}^{(k+1)}}_{2p\times 2p}
\;
\underbrace{\widetilde{V}_{i}^{(k+1)T}}_{2p\times \frac{N}{2^{k}}},\;\; \forall i\in\{1,2,\dots 2^k\},
\label{eq:K_inverse_sizes_k}
\end{equation}
where 
$S^{(k+1)}_i = \Big(I + {\widetilde{V}_{i}^{(k+1)^T}}\widetilde U_{i}^{(k+1)}\Big)^{-1}$ and the dimensions of the matrices are indicated by the underbraces beneath each matrix. It can be observed that $K_{i}^{(k,\kappa-k)^{-1}}$ is a rank-$2p$ perturbation of the identity matrix.
Further, the matrices $\widetilde{U}_{i}^{(k+1)}$ and $\widetilde{V}_{i}^{(k+1)^T}$, for all $k\in\{\kappa-1,\hdots,0\}$, are block diagonal and block anti-diagonal matrices respectively.

\begin{figure}[ht]
    \centering
    \begin{tikzpicture}[scale=1]

\def\cs{0.36}
\def\N{8}

% =================================================================
%  \blk{Xorigin}{col}{row}{w}{h}{colour}
% =================================================================
\newcommand{\blk}[6]{%
  \fill[#6] ({#1+#2*\cs},{-#3*\cs}) rectangle ++({#4*\cs},{-#5*\cs});
  \draw[black,very thin] ({#1+#2*\cs},{-#3*\cs}) rectangle ++({#4*\cs},{-#5*\cs});
}

% =================================================================
%  \identblk{Xorigin}{col}{row}{size}
% =================================================================
\newcommand{\identblk}[4]{%
  \fill[white] ({#1+#2*\cs},{-#3*\cs}) rectangle ++({#4*\cs},{-#4*\cs});
  \draw[black,very thin] ({#1+#2*\cs},{-#3*\cs}) rectangle ++({#4*\cs},{-#4*\cs});
  \draw[identblue,thin] ({#1+#2*\cs},{-#3*\cs}) -- ++({#4*\cs},{-#4*\cs});
}

% ================================================================
%  HORIZONTAL OFFSETS (SWAPPED ORDER)
% ================================================================
\def\matK{0}

\pgfmathsetmacro{\approxX}{\N*\cs+0.35}
\pgfmathsetmacro{\matA}{\approxX+0.55}

\pgfmathsetmacro{\eqX}{\matA+\N*\cs+0.35}

% --- swapped order ---
\pgfmathsetmacro{\matE}{\eqX+0.55}        % K_0
\pgfmathsetmacro{\txB}{\matE+\N*\cs+0.2}

\pgfmathsetmacro{\matD}{\txB+0.35}        % K_1
\pgfmathsetmacro{\txC}{\matD+\N*\cs+0.2}

\pgfmathsetmacro{\matC}{\txC+0.35}        % K_2
\pgfmathsetmacro{\txD}{\matC+\N*\cs+0.2}

\pgfmathsetmacro{\matB}{\txD+0.35}        % K_3

% ================================================================
%  MATRIX K (full red)
% ================================================================
\blk{\matK}{0}{0}{\N}{\N}{fullrank}
\draw[thick] ({\matK},0) rectangle ++({\N*\cs},{-\N*\cs});
\node[below=5pt,font=\small] at ({\matK+\N*\cs/2},{-\N*\cs}) {$K^{-1}$};

\node[font=\Large] at ({\approxX},{-\N*\cs/2}) {$\approx$};

% ================================================================
%  MATRIX K^{(3)}
% ================================================================
\blk{\matA}{4}{0}{4}{4}{lowrank}
\blk{\matA}{0}{4}{4}{4}{lowrank}

\blk{\matA}{2}{0}{2}{2}{lowrank}
\blk{\matA}{0}{2}{2}{2}{lowrank}
\blk{\matA}{6}{4}{2}{2}{lowrank}
\blk{\matA}{4}{6}{2}{2}{lowrank}

\foreach \g in {0,2,4,6}{
  \pgfmathtruncatemacro{\gp}{\g+1}
  \blk{\matA}{\gp}{\g}{1}{1}{lowrank}
  \blk{\matA}{\g}{\gp}{1}{1}{lowrank}
}

\foreach \i in {0,...,7}{
  \blk{\matA}{\i}{\i}{1}{1}{fullrank}
}

\draw[thick] ({\matA},0) rectangle ++({\N*\cs},{-\N*\cs});
\node[below=5pt,font=\small] at ({\matA+\N*\cs/2},{-\N*\cs}) {$K^{(3)^{-1}}$};

\node[font=\Large] at ({\eqX},{-\N*\cs/2}) {$=$};

% ================================================================
%  MATRIX K_0
% ================================================================
\identblk{\matE}{0}{0}{4}
\identblk{\matE}{4}{4}{4}
\blk{\matE}{4}{0}{4}{4}{lowrank}
\blk{\matE}{0}{4}{4}{4}{lowrank}
\draw[thick] ({\matE},0) rectangle ++({\N*\cs},{-\N*\cs});
\node[below=5pt,font=\small] at ({\matE+\N*\cs/2},{-\N*\cs}) {$K_0^{-1}$};

\node[font=\Large] at ({\txB},{-\N*\cs/2}) {$\times$};

% ================================================================
%  MATRIX K_1
% ================================================================
\identblk{\matD}{0}{0}{2}
\identblk{\matD}{2}{2}{2}
\blk{\matD}{2}{0}{2}{2}{lowrank}
\blk{\matD}{0}{2}{2}{2}{lowrank}
\identblk{\matD}{4}{4}{2}
\identblk{\matD}{6}{6}{2}
\blk{\matD}{6}{4}{2}{2}{lowrank}
\blk{\matD}{4}{6}{2}{2}{lowrank}
\draw[thick] ({\matD},0) rectangle ++({\N*\cs},{-\N*\cs});
\node[below=5pt,font=\small] at ({\matD+\N*\cs/2},{-\N*\cs}) {$K_1^{-1}$};

\node[font=\Large] at ({\txC},{-\N*\cs/2}) {$\times$};

% ================================================================
%  MATRIX K_2
% ================================================================
\foreach \g in {0,...,3}{
  \pgfmathtruncatemacro{\base}{2*\g}
  \pgfmathtruncatemacro{\basep}{\base+1}
  \identblk{\matC}{\base}{\base}{1}
  \identblk{\matC}{\basep}{\basep}{1}
  \blk{\matC}{\basep}{\base}{1}{1}{lowrank}
  \blk{\matC}{\base}{\basep}{1}{1}{lowrank}
}
\draw[thick] ({\matC},0) rectangle ++({\N*\cs},{-\N*\cs});
\node[below=5pt,font=\small] at ({\matC+\N*\cs/2},{-\N*\cs}) {$K_2^{-1}$};

\node[font=\Large] at ({\txD},{-\N*\cs/2}) {$\times$};

% ================================================================
%  MATRIX K_3
% ================================================================
\foreach \i in {0,...,7}{
  \blk{\matB}{\i}{\i}{1}{1}{fullrank}
}
\draw[thick] ({\matB},0) rectangle ++({\N*\cs},{-\N*\cs});
\node[below=5pt,font=\small] at ({\matB+\N*\cs/2},{-\N*\cs}) {$K_3^{-1}$};

% ================================================================
%  LEGEND
% ================================================================
\pgfmathsetmacro{\legY}{-\N*\cs - 1.25}
\pgfmathsetmacro{\legMid}{(\matB+\N*\cs)/2}

\pgfmathsetmacro{\lA}{\legMid - 5.8}
\fill[fullrank] (\lA,\legY) rectangle ++(0.3,0.3);
\node[right,font=\small] at ({\lA+0.4},{\legY+0.15}) {Full Rank};

\pgfmathsetmacro{\lB}{\lA + 2.5}
\fill[lowrank] (\lB,\legY) rectangle ++(0.3,0.3);
\node[right,font=\small] at ({\lB+0.4},{\legY+0.15}) {Low Rank};

\pgfmathsetmacro{\lC}{\lB + 2.5}
\draw[thick] (\lC,\legY) rectangle ++(0.3,0.3);
\node[right,font=\small] at ({\lC+0.4},{\legY+0.15}) {Zero Matrix};

\pgfmathsetmacro{\lD}{\lC + 2.8}
\draw[thick] (\lD,\legY) rectangle ++(0.3,0.3);
\draw[identblue,thin] (\lD,{\legY+0.3}) -- ++(0.3,-0.3);
\node[right,font=\small] at ({\lD+0.4},{\legY+0.15}) {Identity Matrix};

\end{tikzpicture}
    \caption{Factorization of inverse of HODLR matrix}
    \label{fig:inversefactorization}
\end{figure}
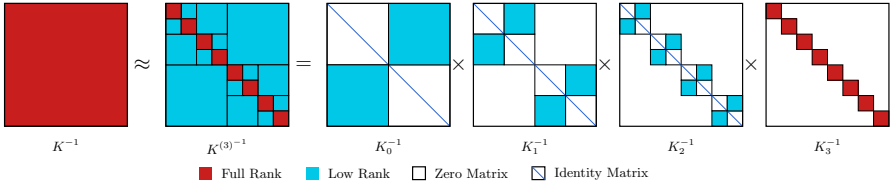

Figure~\ref{fig:inversefactorization} illustrates the structure of the factorization of \(K^{(\kappa)^{-1}}\). The matrices appearing in the factorization are block diagonal, where \(K_k^{-1}\) contains \(2^k\) diagonal blocks, each of size \(N/2^k \times N/2^k\), for all $k\in\{0,1,\hdots,\kappa\}$. The diagonal blocks of \(K_{\kappa}^{-1}\) are dense matrices, whereas the diagonal blocks of the remaining factors are rank-\(2p\) perturbations of the identity matrix.

Applying $K_{k}^{-1}$ to a vector costs $\mathcal{O}(pN)$, where $k\in\{0,1,\hdots,\kappa\}$ (the derivation follows the same lines as the derivation we have seen in the factorization phase). Hence, the total cost of the solve phase is $\mathcal{O}(pN\log(N))$.

\section{Fast direct solver as a neural network} \label{sec:direct}
% In this section, we discuss the design of a neural network that solves a linear system involving a HODLR matrix. 
In this section, we present the design of a neural network for solving linear systems involving HODLR matrices, based on the solver described in Section~\ref{sec:directsolver}. First, we recall the standard feedforward network formulation.
Consider a feed-forward neural network, denoted by $\mathcal{N}$, with input $v$, output $u$, and $D+1$ layers, denoted by
$$
u = \mathcal{N}(v), \qquad v \in \mathbb{R}^{n_D}, u \in \mathbb{R}^{n_0},
$$
which can be expressed recursively as
\begin{align}
\xi^{(D)} &= v, \\
\xi^{(k)} &= \phi\!\left(W^{(k)} \xi^{(k+1)} + b^{(k)}\right), \qquad k \in\cbkt{ D-1, D-2,\dots, 0} \\
u &= \xi^{(0)}.
\end{align}
The layers are indexed by $\cbkt{D, D-1,\hdots,0}$, where $D$ denotes the input layer and $0$ denotes the output layer. The number of neurons in layer $l$ is denoted by $n_{l}$.
In the terminology of neural networks, the vector \(\xi^{(k)} \in \mathbb{R}^{n_k}\) denotes the \textit{activations}, the matrix \(W^{(k)} \in \mathbb{R}^{n_k \times n_{k+1}}\) denotes the \textit{weights}, and the vector \(b^{(k)} \in \mathbb{R}^{n_k}\) denotes the \textit{biases}, for all $k\in\{ 1,2,\dots D\}$. The function \(\phi\) is the activation function, applied element-wise, and typical choices include the identity (linear), the rectified linear unit (ReLU), and the sigmoid.

\begin{remark}\label{remark:ffn}
If the activation function is linear (i.e., $\phi(x)=x$) and the bias is set to zero, then each feedforward layer reduces to a matrix-vector multiplication $\xi^{(k)} = W^{(k)} \xi^{(k+1)}$. Hence, a multi-layer neural network built with linear activation and zero biases is a sequence of matrix-vector multiplications.   
\end{remark}
% , which is the structure of the factorized inverse $K^{-1} = \prod_{i} K_i^{-1}$.

\subsection{Locally connected networks}\label{sec:LCNN}
% We now describe the representation of the solve phase of the direct solver through the framework of locally connected networks.
% \begin{itemize}[label={}, leftmargin=1em, itemsep=0pt, parsep=0pt, topsep=0pt]
    % \item 
It follows from Section~\ref{sec:solve_phase} that to solve for $x$, we need to compute 
\begin{equation}\label{eq:solve}
    x = K^{(\kappa)^{-1}}b = K_{0}^{-1} K_{1}^{-1}\dots K_{\kappa-1}^{-1}K_\kappa^{-1}b,
\end{equation}
where $K_\kappa^{-1}$ is first applied to $b$, followed successively by $K_{\kappa-1}^{-1}$ and continuing up to  $K_0^{-1}$. Thus, $x$ is obtained by performing a sequence of matrix-multiplications on $b$, where all the matrices are block diagonal. Following Remark~\ref{remark:ffn}, the mapping from \(b\) to \(x\) can be represented by a multilayer neural network that takes \(b\) as input and produces \(x\) as output, with \(D=\kappa+1\) layers using linear activations and zero biases.

Another important observation is that the sequence of matrices \\$K_{0}^{-1}, K_{1}^{-1},\dots K_{\kappa-1}^{-1},K_\kappa^{-1}$ are block-diagonal. At level $k$, $K_{k}^{-1}$ is applied to a vector, say, $h$, and let $y=K_{k}^{-1}h$. Let $h$ and $y$ be partitioned into $2^{k}$ segments as 
\begin{equation}
    h=\begin{bmatrix} 
h_1 \\ 
h_2 \\
\vdots\\
h_{2^{k}}
\end{bmatrix}, 
y=\begin{bmatrix} 
y_1 \\ 
y_2 \\
\vdots\\
y_{2^{k}}
\end{bmatrix}, \text{ where}
\end{equation}

\begin{equation}\label{eq:LCNN_segmenting}
   y =
   \begin{bmatrix} 
y_1 \\ 
y_2 \\
\vdots\\
y_{2^{k}}
\end{bmatrix}
= K_{k}^{-1}h = 
\begin{bmatrix} 
K_{1}^{(k,\kappa-k)^{-1}}\\ & K_{2}^{(k,\kappa-k)^{-1}}\\ & & \ddots & \\ & & & K_{2^k}^{(k,\kappa-k)^{-1}} \end{bmatrix}
\begin{bmatrix} 
h_1 \\ 
h_2 \\
\vdots\\
h_{2^{k}}
\end{bmatrix}
\end{equation}
Since \(K_k^{-1}\) is block diagonal, the \(i^{\text{th}}\) segment of \(y\),
\[
y_i = K_i^{(k,\kappa-k)^{-1}} h_i,\;\;\; \forall \;i\in\cbkt{1,2,\dots,2^{k}},\;\; \forall\; k\;\in\cbkt{0,1,\dots,\kappa}
\]
depends only on the corresponding input segment \(h_i\). This localised dependence can be represented using \textit{locally connected networks (LCNNs)} instead of standard feedforward networks. In the terminology of neural networks, LCNNs are those networks where neurons in each layer are connected only to a local subset of neurons in the previous layer. Hence, the computation in Equation~\eqref{eq:solve} can be represented using an LCNN framework.

\subsection{Neural network architecture} \label{sec:3.2}
% We now describe the neural network architecture whose forward pass has the same computational framework as that of the solve phase of the direct solver described in Section~\ref{sec:solve_phase}.
We now describe a locally connected neural network (LCNN) architecture whose forward pass follows the same computational framework as that of the solve phase of the direct solver presented in Section~\ref{sec:solve_phase}.
% The solve phase computes $x$ as described in Equation~\eqref{eq:solve}, restated here for reference.
% \begin{equation}\label{eq:solve2}
%     x = K_{0}^{-1} K_{1}^{-1}\dots K_{\kappa-1}^{-1}K_\kappa^{-1}b,
% \end{equation}
% To apply $K^{(\kappa)^{-1}}$ to $b$, first we apply $K_{\kappa}^{-1}$ to 
% $b$ followed by $K_{\kappa-1}^{-1}$ all the way upto $K_{0}^{-1}$. Each of the $K_{k}^{-1}$ for all $k\in\{0,1\dots,\kappa\}$ has block diagonal structure. The diagonal blocks of $K_{\kappa}^{-1}$ are full-rank matrices. 
% The diagonal blocks of $K_{k}^{-1}$ for all $k\in\{0,1\dots,\kappa-1\}$ have the structure as described in Equation~\eqref{eq:K_inverse_sizes_k}, rewritten here for easy reference.
% We consider 
The following describes the sequence of operations, or equivalently the various layers present in the network.% To construct the neural network, we model the sequence of operations as follows:
\begin{enumerate}[leftmargin=2.5em]
\item \textbf{Initialization:} The input layer takes vector \(b\) as input to the network, which is represented by \(\xi^{(\kappa+1)}\).
\begin{align}
\text{\textbf{Input layer:\;\;\;}}  \xi^{(\kappa+1)} = b\label{eq:input_layer}
\end{align}
% \xi^{(\kappa)} &= K_{\kappa}^{-1}\,\xi^{(\kappa+1)}.
\item \textbf{$K^{(\kappa)}$ layer:}
% The initial hidden layer applies $K_{\kappa,i}^{-1}$ to its $i^{th}$ input. 
Let $\xi^{(\kappa+1)}$ be partitioned as $\xi^{(\kappa+1)} = \sbkt{\xi^{(\kappa+1)^T}_1,\dots,\xi^{(\kappa+1)^T}_{2^{\kappa}}}^T$ such that each partition corresponds to a leaf node of the HODLR tree. 
Each input segment \(\xi_i^{(\kappa+1)}\) is multiplied by the corresponding block inverse matrix \(K^{(\kappa)^{-1}}_i\). Since each segment is processed independently, this layer follows a locally connected structure. 
\begin{equation}
% &\mathllap{\text{(i)}\quad} &\text{\textbf{$K^{(\kappa)}$-layer:\;\;}} 
\xi^{(\kappa)}_i = \textstyle\underbrace{K^{(\kappa)^{-1}}_i}_{\frac{N}{2^{\kappa}}\times\frac{N}{2^{\kappa}}} \xi^{(\kappa+1)}_i,
\quad \forall i\in\{1,\dots,2^\kappa\}, \label{eq:K_kappa_layer}
\end{equation}
% Each of the partitioned $\xi^{(\kappa+1)}$ is given as input to the network. 
\item \textbf{Flatten layer:}
The resulting collection of segment-wise outputs from the $K^{(\kappa)}$-layer is then flattened to form the vector \(\xi^{(\kappa)}\), facilitating the tensor operations in the subsequent layer.
\begin{equation}
\xi^{(\kappa)} := \mathrm{Flatten}\bkt{\sbkt{\xi^{(\kappa)}_{1}, \hdots,\xi^{(\kappa)}_{2^\kappa}}}\label{eq:flatten_layer1}
\end{equation}

% \begin{alignat}{2}
% &\mathllap{\text{(i)}\quad} &\text{\textbf{$K^{(\kappa)}$-layer:\;\;}} 
% & \xi^{(\kappa)}_i = \textstyle\underbrace{K^{(\kappa)^{-1}}_i}_{\frac{N}{2^{\kappa}}\times\frac{N}{2^{\kappa}}} \xi^{(\kappa+1)}_i,
% \quad \forall i\in\{1,\dots,2^\kappa\}, \label{eq:K_kappa_layer} \\
% &\mathllap{\text{(ii)}\quad} &\text{\textbf{Flatten layer:\;\;}} 
% & \xi^{(\kappa)} := \mathrm{flatten}\bkt{\sbkt{\xi^{(\kappa)}_{2^\kappa}, \hdots,\xi^{(\kappa)}_{2^\kappa}}}\label{eq:flatten_layer1}
% \end{alignat}

% $\xi^{(\kappa)} = \sbkt{\xi^{(\kappa)^T}_1, \dots, \xi^{(\kappa)^T}_{2^k}}^T_{\frac{N}{2^{\kappa}}$
\item \textbf{Iterative layers:}
The subsequent hidden layers apply $K_k^{-1}$ for \\$k \in \{0,1,\dots,\kappa-1\}$. \(K_k^{-1}\) is applied to a vector, say \(h\), in a segment-wise manner. 
To avoid notational clutter, from here on, we refer to $K_{i}^{(k,\kappa-k)^{-1}}$, the $i^{th}$ diagonal block of $K_{k}^{-1}$,  as $K^{(k)^{-1}}_i$. 
\begin{equation}
    K^{(k)^{-1}}_i = I - \widetilde{U}^{(k+1)}_i S^{(k+1)}_i \widetilde{V}^{(k+1)^T}_i, \;\;\forall\; k\in\cbkt{\kappa-1,\hdots,0}, \;\;\forall\; i\in\{1,\dots,2^k\}. 
\end{equation}
For each segment $h_i$ of the vector $h$, the operations are performed sequentially by first applying \(\widetilde{V}_{i}^{(k+1)^T}\), followed by \(S_i^{(k+1)}\), and then \(\widetilde{U}_i^{(k+1)}\). The resulting vector is subsequently subtracted from \(h_i\). For all \(k \in \cbkt{\kappa-1,\hdots,0}\), these operations are performed iteratively and are implemented through the following layers.
% For each level $k \in\cbkt{ \kappa-1, \ldots, 0}$, five sequential operations are performed on the segmented input.
% which are to apply $\widetilde{V}_{i}^{(k)\mathsf T}$, followed by $W^{(k)}_i$, followed by $\Tilde{U}^{(k)}_i$, and to subtract the result from its input. 
% $V^{(k)}$-layer, $W^{(k)}$-layer, $U^{(k)}$-layer and the residual layer.
\begin{enumerate}[leftmargin=*]
\item \textbf{Swap layer:}
    % Let $\xi^{(k+1)}$ be partitioned into $2^{k+1}$ segments as $\xi^{(k+1)} = \sbkt{\xi^{(k+1)^T}_1,\dots,\xi^{(k+1)^T}_{2^{k+1}}}^T$. 
Let the $2^{k+1}$ segments of $\xi^{(k+1)}$ be denoted by $\xi^{(k+1)} = \sbkt{\xi^{(k+1)^T}_1,\dots,\xi^{(k+1)^T}_{2^{k+1}}}^T$. 
    % where each segment corresponds to a node at level $k$ of the HODLR tree structure.
    % Since $\widetilde{V}_{i}^{(k)\mathsf T}$ is block anti-diagonal with its two blocks of size $p\times \frac{N}{2^{k}}$, applying it to a vector $[\xi^{(\kappa)^{T}}_{2i-1},\; \xi^{(\kappa)^T}_{2i}]^T$ results in 
    $\widetilde{V}_{i}^{(k+1)^ T}$ is block anti-diagonal with its two anti-diagonal blocks of size $p\times \frac{N}{2^{k+1}}$. Applying
    \[
    \widetilde{V}_{i}^{(k+1)^T} =
    \begin{bmatrix}
    0 & V_{2i, 2i-1}^{(k+1)^T} \\
    V_{2i-1, 2i}^{(k+1)^T} & 0
    \end{bmatrix}
    \]
    to a vector consisting of two segments is equivalent to applying a transformed matrix to a vector with swapped segments. The transformed matrix is obtained by rewriting the anti-diagonal block matrix as a diagonal block matrix, as illustrated in Equation~\eqref{eq:swap_V}. 
    % Applying \begin{equation*}
    %     \widetilde{V}_{i}^{(k)^ T} = \begin{bmatrix}
    %     0 & V_{2i, 2i-1}^{(k)^{T}}\\
    %     V_{2i-1, 2i}^{(k)^{T}} & 0
    % \end{bmatrix}
    % \end{equation*}
    % to a vector of two segments is equivalent to applying a transformed matrix, that is obtained by writing the anti-diagonal block matrix as a diagonal block matrix, to a vector whose segments are swapped, as shown in Equation~\eqref{eq:swap_V}. 
    % {\small
    \begin{align}\label{eq:swap_V}
    {
    \scriptsize
    \begin{bmatrix}
        0 & V_{2i, 2i-1}^{(k+1)^{T}}\\
        V_{2i-1, 2i}^{(k+1)^{T}} & 0
    \end{bmatrix}
    \begin{bmatrix}
        \xi^{(k+1)}_{2i-1}\\
        \xi^{(k+1)}_{2i}
    \end{bmatrix} =         \begin{bmatrix}
        V_{2i, 2i-1}^{(k+1)^{T}} & 0\\
        0 & V_{2i-1, 2i}^{(k+1)^{T}}
    \end{bmatrix}
    \begin{bmatrix}
        \xi^{(k+1)}_{2i}\\
        \xi^{(k+1)}_{2i-1}
    \end{bmatrix}=
    \begin{bmatrix}
        V_{2i, 2i-1}^{(k+1)^{T}}\xi^{(k+1)}_{2i}\\
        V_{2i-1, 2i}^{(k+1)^{T}}\xi^{(k+1)}_{2i-1}
    \end{bmatrix}
    }
\end{align}
% }
Rewriting the matrix as a diagonal matrix makes it amenable to model it using an LCNN layer.
Therefore, the following swap operation is performed, before applying the transformed matrix.
    \begin{equation}
        \bkt{\xi^{(k+1)}_{2i-1}, \xi^{(k+1)}_{2i}} := \bkt{\xi^{(k+1)}_{2i}, \xi^{(k+1)}_{2i-1}},\; \forall\;i\in\{1,\hdots,2^{k}\} 
    \end{equation}
    
\item \text{\textbf{$V^{(k)}$ layer:}} Apply $\widetilde{V}_{i}^{(k+1)^ T}$ as described in Equation~\eqref{eq:swap_V}, to compute $\eta^{(k)}_{i}$ \\
    \begin{equation}
            \eta^{(k)}_{i} = 
        \begin{bmatrix}
            V_{2i, 2i-1}^{(k+1)^{T}}\xi^{(k+1)}_{2i}\\
            V_{2i-1, 2i}^{(k+1)^{T}}\xi^{(k+1)}_{2i-1}
        \end{bmatrix}
            , \quad \forall \;i\in\{1,\dots,2^{k}\}
            % ,\quad \eta^{(k)} = \sbkt{\eta^{(k)}_1, \dots, \eta^{(k)}_{2^{k}}}_{2p \times 2^{k}} \label{eq:V_layer}
        % \eta^{(k)}_i = \textstyle\underbrace{V_{i}^{(k)\mathsf T}}_{p\times\frac{N}{2^{k+1}}}  \xi^{(k+1)}_i, \quad \forall i\in\{1,\dots,2^{k+1}\},\quad \eta^{(k)} = \sbkt{\eta^{(k)}_1, \dots, \eta^{(k)}_{2^{k+1}}}_{p \times 2^{k+1}}
        \label{eq:V_layer}
    \end{equation}

\item \text{\textbf{$S^{(k)}$ layer:}} Compute\\
    \begin{equation}
        \zeta^{(k)}_i = \textstyle\underbrace{S^{(k+1)}_i}_{2p\times 2p} \eta^{(k)}_i, \quad \forall\; i\in\{1,\dots,2^k\}
        % , \quad \zeta^{(k)} = \sbkt{\zeta^{(k)}_1, \dots, \zeta^{(k)}_{2^k}}_{2p \times 2^k}
        \label{eq:W_layer}
    \end{equation}

\item \text{\textbf{Reshape layer:}} 
The output of the \(S^{(k)}\)-layer consists of \(2^k\) segments, each of dimension \(2p\). Since the subsequent \(U^{(k)}\)-layer operates on \(2^{k+1}\) segments of dimension \(p\), the tensor organization of \(\zeta^{(k)}\) must be modified accordingly. Each segment \(\zeta_i^{(k)} \in \mathbb{R}^{2p}\) is partitioned into two \(p\)-dimensional segments, thereby reshaping the tensor dimensions from \((2^k,\,2p)\) to\\ \((2^{k+1},\,p)\), i.e.,
\[
\zeta^{(k)} := \mathrm{reshape}\bigl(\zeta^{(k)},\,(2^{k+1},\,p)\bigr).
\]
% The output of the \(S^{(k)}\)-layer consists of \(2^k\) segments, each of dimension \(2p\). Since the subsequent \(U^{(k)}\)-layer operates on \(2^{k+1}\) segments of dimension \(p\), each segment \(\zeta_i^{(k)} \in \mathbb{R}^{2p}\) is partitioned into two \(p\)-dimensional segments. Thus, the collection of vectors is reshaped from dimensions \((2^k,\,2p)\) to \((2^{k+1},\,p)\), i.e.,
% \[
% \zeta^{(k)} = \mathrm{reshape}\bigl(\zeta^{(k)},\,(2^{k+1},\,p)\bigr).
% \]

\item \text{\textbf{$U^{(k)}$ layer:}} 
$\widetilde{U}^{(k+1)}_i$ is to be applied to $\zeta^{(k)}_i$. Let $\zeta^{(k)}_i$ be partitioned into two blocks as $\zeta^{(k)}_i = \begin{bmatrix}
    \zeta^{(k)}_{i,1}\\
    \zeta^{(k)}_{i,2}
\end{bmatrix}$.
Since 
\begin{equation}
    \widetilde{U}_{i}^{(k+1)} = 
    \begin{bmatrix}
        U_{2i-1}^{(k+1,\kappa-k)} & 0\\
        0 & U_{2i}^{(k+1,\kappa-k)}
    \end{bmatrix},
\end{equation}
\begin{flalign}
& \text{compute } \chi^{(k)}_i =  \widetilde{U}_{i}^{(k+1)}\zeta^{(k)}_i = 
\begin{bmatrix}
\widetilde{U}_{2i-1}^{(k+1,\kappa-k)}\zeta^{(k)}_{i,1}\\
\widetilde{U}_{2i}^{(k+1,\kappa-k)}\zeta^{(k)}_{i,2}
\end{bmatrix}
, \quad \forall\; i\in\{1,\dots,2^k\} &
\label{eq:U_layer}
\end{flalign}
    % , \quad \chi^{(k)} = \sbkt{\chi^{(k)}_1, \dots, \chi^{(k)}_{2^k}}_{\frac{N}{2^{k}} \times 2^k}\label{eq:U_layer} 
    % \\[6pt]
\item \text{\textbf{Residual layer:}} The resulting collection of segment-wise outputs from the $U^{(k)}$-layer is then flattened and then subtracted from \(\xi^{(k+1)}\).
    \begin{equation}
        \xi^{(k)} = \xi^{(k+1)} - \text{Flatten}\bkt{\chi^{(k)}}. \label{eq:res_layer}
    \end{equation}
\end{enumerate}
We refer to the collection of all the layers described in steps (a) to (f) above, for iteration $k$, by \textit{SMW block} $k$.

\item \textbf{Output layer:}
The final output \(x\) is obtained at \(k=0\), and is given by
\begin{align}
&x = \xi^{(0)}
\label{eq:output_layer}
\end{align}
\end{enumerate}

Each of the \(K^{(\kappa)}\), \(V^{(k)}\), \(S^{(k)}\), and \(U^{(k)}\) layers described above is implemented using locally connected layers using tensor operations.
The LCNN representation of these layers is described below.
\begin{itemize}[leftmargin=*]
    \item 
    The $K^{(\kappa)}$ layer, described in Equation~\eqref{eq:K_kappa_layer}, is modeled using a locally connected layer with stride length $s=\frac{N}{2^\kappa}$, kernel window size $w=\frac{N}{2^\kappa}$, and number of filters $f=\frac{N}{2^\kappa}$. 
    
    \item 
    The $V^{(k)}$ layer, described in Equation~\eqref{eq:V_layer}, is modeled using a locally connected layer with stride length $s$, kernel window size $w$ as $s=w=\frac{N}{2^{k+1}}$, and number of filters $f=p$. 
    \item 
    The $S^{(k)}$ layer, described in Equation~\eqref{eq:W_layer}, is modeled using a locally connected layer with stride length and kernel window size $s=w=2$, and number of filters $f=2p$. 
    \item 
    The $U^{(k)}$ layer, described in Equation~\eqref{eq:U_layer}, is modeled using a locally connected layer with stride length and kernel window size $s=w=1$, and number of filters $f=\frac{N}{2^{k+1}}$.
\end{itemize}
% Combining Equations~\eqref{eq:input_layer}, ~\eqref{eq:K_kappa_layer}, ~\eqref{eq:flatten_layer1}, ~\eqref{eq:V_layer}, ~\eqref{eq:W_layer}, ~\eqref{eq:U_layer}, ~\eqref{eq:res_layer}, and~\eqref{eq:output_layer}, we describe the neural network architecture that represents the fast direct solver algorithm in Algorithm~\ref{alg:linear}.
% Algorithm~\ref{alg:linear} describes the forward pass of the proposed neural network architecture that emulates the computational steps of the fast direct solver for computing \(K^{-1}b\).
\begin{algorithm}[htbp]
\caption{Forward pass of the neural network that computes $K^{(\kappa)^{-1}}b$. Here, LC denotes locally connected layer.}
\label{alg:linear}
\begin{algorithmic}[1]
\REQUIRE The input vector $b\in \mathbb{R}^{N}$
\REQUIRE The number of levels $\kappa$ in HODLR tree structure
\STATE $\xi^{(\kappa+1)}= b;$
\STATE $\xi^{(\kappa)} = \text{LC}\sbkt{\text{activation=linear}; s=\frac{N}{2^\kappa}; w=\frac{N}{2^\kappa}; f=\frac{N}{2^\kappa}}\bkt{\xi^{(\kappa+1)}};$ \mycomment{$K^{(\kappa)}$ layer}
\STATE $\xi^{(\kappa)} = \text{Flatten}\sbkt{\xi^{(\kappa)}};$
\mycomment{Flatten layer}
% \STATE $\xi^{(\kappa)}  = \textsf{Flatten}(\xi^{(\kappa)});$
\FOR{ $k = \kappa-1$ to $0$ }
    \STATE Swap sibling blocks of $\xi^{(k+1)} = \sbkt{\xi^{(k+1)}_1, \xi^{(k+1)}_2, \dots, \xi^{(k+1)}_{2^{k+1}}}$, i.e. swap $\xi^{(k+1)}_{2i-1}$ with $\xi^{(k+1)}_{2i}$, $\forall\;i\in\{1,\hdots2^{k}\}$
    
    \mycomment{Swap layer}
    \STATE $\eta^{(k)} = \text{LC}\sbkt{\text{activation=linear}; s=\frac{N}{2^{k+1}}; w=\frac{N}{2^{k+1}}; f=p}\bkt{\xi^{(k+1)}};$ \mycomment{$V^{(k)}$ layer}
    \STATE $\zeta^{(k)} = \text{LC}\sbkt{\text{activation=linear}; s=2; w=2; f=2p}\bkt{\eta^{(k)}};$ \mycomment{$S^{(k)}$ layer}
    \STATE $\zeta^{(k)} = \text{Reshape}\bkt{\zeta^{(k)},\,\bkt{2^{k+1},\,p}}$
    \mycomment{Reshape layer}
    \STATE $\chi^{(k)} = \text{LC}\sbkt{\text{activation=linear}; s=1; w=1; f=\frac{N}{2^{k+1}}}\bkt{\zeta^{(k)}};$ \mycomment{$U^{(k)}$ layer}
    % \STATE $\zeta^{(i)} = \textsf{Flatten}(\chi^{(i)});$
    \STATE $\xi^{(k)} = \xi^{(k+1)}  - \text{Flatten}\bkt{\chi^{(k)}};$
    \mycomment{Residual layer}
\ENDFOR
\STATE $x_{N} = \xi^{(0)}$
\RETURN $x_{N}$
\end{algorithmic}
\end{algorithm}

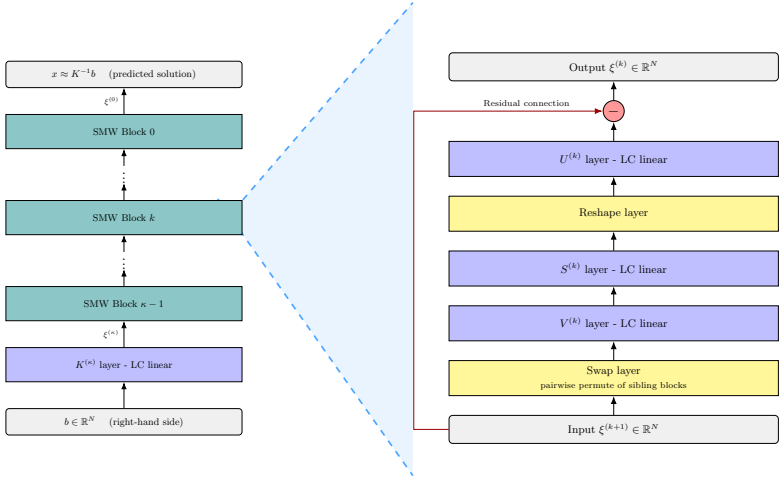
\begin{figure}[htbp]
% \hspace*{1.1cm} 
% \centering
% Set a common visual height for both subfigure bodies.
% Subfigures use [b]: captions align at the bottom baseline.
% Inside each subfigure, the tikz is centered in a fixed-height \vbox.
\newlength{\panelht}\setlength{\panelht}{8.2cm}
\centering
\resizebox{0.8\linewidth}{!}{
\begin{tikzpicture}

% --- Subfigure (a): high-level FDSNet architecture ---
\node[inner sep=0pt, anchor=south] (subA) at (1.6cm,0) {%
\begin{subfigure}[b]{0.42\textwidth}
\centering
\begin{minipage}[c][\panelht][c]{\linewidth}
\centering
\resizebox{\linewidth}{!}{%
\begin{tikzpicture}[
    >=Latex,
    node distance=0.85cm,
    every node/.style={font=\small},
    io/.style    ={rectangle, rounded corners=3pt, draw, fill=gray!12,
                   minimum width=7.6cm, minimum height=0.85cm, align=center},
    leaf/.style  ={rectangle, draw, fill=blue!25,
                   minimum width=7.6cm, minimum height=1.1cm, align=center},
    smw/.style   ={rectangle, draw, fill=teal!45,
                   minimum width=7.6cm, minimum height=1.1cm, align=center},
    dimlbl/.style={font=\scriptsize, anchor=east},
    sidelbl/.style={font=\scriptsize\itshape}
]
\node (input)  [io] {$b\in\mathbb{R}^{N}$ \quad (right-hand side)};
\node (kleaf)  [leaf,  above=of input]
       {$K^{(\kappa)}$ layer - LC linear};
\node (smwk1)  [smw,   above=of kleaf] {SMW Block $\kappa-1$};
\node (sdots)  [above=0.55cm of smwk1, font=\Large] {$\vdots$};
\node (smwj)   [smw,   above=0.55cm of sdots] {SMW Block $k$};
\node (sdots2) [above=0.55cm of smwj,  font=\Large] {$\vdots$};
\node (smw0)   [smw,   above=0.55cm of sdots2] {SMW Block 0};
\node (output) [io,    above=of smw0]
       {$x\approx K^{-1}b$ \quad (predicted solution)};

\draw[->, thick] (input)  -- node[dimlbl, xshift=-0.15cm] {}    (kleaf);
\draw[->, thick] (kleaf)  -- node[dimlbl, xshift=-0.15cm] {$\xi^{(\kappa)}$} (smwk1);
\draw[->, thick] (smwk1)  -- (sdots);
\draw[->, thick] (sdots)  -- (smwj);
\draw[->, thick] (smwj)   -- (sdots2);
\draw[->, thick] (sdots2) -- (smw0);
\draw[->, thick] (smw0)   -- node[dimlbl, xshift=-0.15cm] {$\xi^{(0)}$}   (output);

\draw [opacity=0, decorate, decoration={brace, amplitude=6pt, raise=2pt}, very thick]
  ($(smw0.north east)+(0.35,0)$) -- ($(smwk1.south east)+(0.35,0)$)
  node[midway, sidelbl, right=10pt, align=left, text width=2.4cm]
       {$\kappa$ residual SMW levels ($k=\kappa{-}1,\dots,0$) \\encoding\\$K_0^{-1}K_1^{-1}\!\cdots K_{\kappa-1}^{-1}$};
       % \\encoding $K_0^{-1}K_1^{-1}\!\cdots K_{\kappa-1}^{-1}$

% \node[sidelbl, right=0.55cm of kleaf, align=left, text width=3.6cm]
  % {leaf factor $K_\kappa$\\(single LC layer)};
\end{tikzpicture}}
\end{minipage}
% \caption*{\hspace*{-2.5cm}(a) \hspace{0.05cm}Linear FDSNet architecture}
\label{fig:FDSNet-lin-highlevel}
\end{subfigure}%
};

% --- Subfigure (b): detailed SMW_Block_k ---
\node[inner sep=0pt, anchor=south] (subB) at ($(subA.south east)+(4.2cm,0)$) {%
\begin{subfigure}[b]{0.46\textwidth}
\centering
\begin{minipage}[c][\panelht][c]{\linewidth}
\centering
\resizebox{\linewidth}{!}{%
\begin{tikzpicture}[
    >=Latex,
    node distance=0.55cm,
    every node/.style={font=\small},
    io/.style       ={rectangle, rounded corners=3pt, draw, fill=gray!12,
                      minimum width=9.4cm, minimum height=0.8cm, align=center},
    layer/.style    ={rectangle, draw, fill=blue!25,
                      minimum width=9.4cm, minimum height=1.0cm, align=center},
    reshapelayer/.style={rectangle, draw, fill=yellow!50,
                      minimum width=9.4cm, minimum height=1.0cm, align=center},
    swaplayer/.style={rectangle, draw, fill=yellow!50,
                      minimum width=9.4cm, minimum height=1.0cm, align=center},
    op/.style       ={circle, draw, fill=red!40, font=\large,
                      inner sep=1pt, minimum size=0.6cm},
    dimlbl/.style   ={font=\scriptsize, anchor=west}
]
\node (input)  [io] {Input $\xi^{(k+1)} \in \mathbb{R}^{N}$};
\node (swap)   [swaplayer, above=of input]
       {Swap layer \\ \scriptsize pairwise permute of sibling blocks};
\node (vlayer) [layer, above=of swap]
       {$V^{(k)}$ layer - LC linear};
\node (slayer) [layer, above=of vlayer]
       {$S^{(k)}$ layer - LC linear};
\node (reshapelayer) [reshapelayer, above=of slayer]
       {Reshape layer};
\node (ulayer) [layer, above=of reshapelayer]
       {$U^{(k)}$ layer - LC linear};
\node (sub)    [op, above=of ulayer] {$-$};
\node (output) [io, above=of sub] {Output $\xi^{(k)} \in \mathbb{R}^{N}$};

\draw[->, thick] (input)  -- node[dimlbl, xshift=0.15cm] {}                       (swap);
\draw[->, thick] (swap)   -- node[dimlbl, xshift=0.15cm] {}              (vlayer);
\draw[->, thick] (vlayer) -- node[dimlbl, xshift=0.15cm] {}                  (slayer);
\draw[->, thick] (slayer) -- node[dimlbl, xshift=0.15cm] {}                  (reshapelayer);
\draw[->, thick] (reshapelayer) -- node[dimlbl, xshift=0.15cm] {}                  (ulayer);
\draw[->, thick] (ulayer) -- node[dimlbl, xshift=0.15cm] {}        (sub);
\draw[->, thick] (sub)    -- (output);

\coordinate (skipX)  at ($(input.west) + (-1.0,0)$);
\coordinate (skipBot) at (skipX |- input.west);
\coordinate (skipTop) at (skipX |- sub.west);
\draw[->, thick, red!60!black]
     (input.west) -- (skipBot)
                  -- (skipTop)
                  -- (sub.west)
     node[pos=1.0, font=\scriptsize, black, xshift=-62pt, yshift=4pt, rotate=0, anchor=south]
        {Residual connection};

% \draw [decorate, decoration={brace, amplitude=5pt, raise=2pt}, thick]
%   ($(ulayer.north east)+(0.35,0)$) -- ($(swap.south east)+(0.35,0)$)
%   node[midway, font=\scriptsize\itshape, right=8pt, align=left, text width=2.4cm]
%        {SMW residual update:\\$I-\widetilde{U}^{(k)}_i W^{(k)}_i \widetilde{V}^{(k)\top}_i$};
\end{tikzpicture}}
\end{minipage}
% \caption*{\hspace*{-1cm}(b) \hspace{0.05cm} SMW Block $k$}
\label{fig:FDSNet-lin-smwblock}
\end{subfigure}%
};

% --- Wedge: from (a)'s SMW_Block_k area to (b)'s frame ---
\begin{scope}[on background layer]
  \coordinate (narrowTop) at ($(subA.south east)!0.54!(subA.north east) + (-2.0cm,0)$);
  \coordinate (narrowBot) at ($(subA.south east)!0.60!(subA.north east) + (-2.0cm,0)$);
  \coordinate (wideTop)   at ($(subB.south west)!0.99!(subB.north west) + (0.0cm,0)$);
  \coordinate (wideBot)   at ($(subB.south west)!0.05!(subB.north west) + (0.0cm,0)$);
  \fill[blue!50!cyan, opacity=0.09]
      (narrowTop) -- (wideTop) -- (wideBot) -- (narrowBot) -- cycle;
  \draw[blue!50!cyan!70, dashed, thick] (narrowTop) -- (wideTop);
  \draw[blue!50!cyan!70, dashed, thick] (narrowBot) -- (wideBot);
\end{scope}

\end{tikzpicture}
}
\caption{Linear FDSNet architecture: A neural network to compute the inverse operation of a HODLR matrix.}
\label{fig:FDSNet-linear}
\end{figure}

Algorithm~\ref{alg:linear} outlines the architecture of the neural network described above.     
We refer to this neural network as linear FDSNet (Fast Direct Solver based Neural Network). 
The learnable parameters of this network include the parameters of all the LCNN layers described in the above steps, which includes entries of the matrices 
{\scriptsize$\cbkt{V_{2i, 2i-1}^{(k+1)^{T}}, V_{2i-1, 2i}^{(k+1)^{T}}, S_i^{(k+1)}, \widetilde{U}_{2i-1}^{(k+1,\kappa-k)}, \widetilde{U}_{2i}^{(k+1,\kappa-k)}, \forall\; i\in\cbkt{1,\dots,2^k}, \forall\; k\in\{\kappa-1,\dots,0\}}$} and\\
$\cbkt{K_{i}^{(\kappa)^{-1}}, \forall\; i\in\cbkt{1,\dots,2^\kappa}}$.
% , and scales as $\mathcal{O}\bkt{pN\log\bkt{N}}$. 
Let $x_{N}$ denote the output of Algorithm~\ref{alg:linear} corresponding to the input $b$, and let $K_{N}^{-1}$ represent the map learned by the network. Then, $x_{N} = K_{N}^{-1}(b)$ approximates the inverse operation associated with the HODLR matrix $K^{(\kappa)}$. Figure~\ref{fig:FDSNet-linear} illustrates the architecture described in Algorithm~\ref{alg:linear}.

% =================================================================
% FIGURE 1 -- LINEAR FDSNet ARCHITECTURE (Algorithm 1)
% 2 subfigures. Subfigure bodies centered to a common height;
% captions bottom-aligned at the same baseline.
% =================================================================

% \begin{remark}
% Since
% \begin{equation}
%    S_i^{(k)} = \bkt{I + {\widetilde{V}_{i}^{(k)^T}}\widetilde U_{i}^{(k)}}^{-1},
%    \label{eq:S_equation} 
% \end{equation}
% it can be treated as a non-learnable parameter by explicitly defining it in terms of \(\widetilde{V}_{i}^{(k)^T}\) and \(\widetilde U_{i}^{(k)}\). However, we choose to keep \(S_i^{(k)}\) learnable, since it was observed that enforcing the above relation results in a higher training time compared to the learnable formulation. After training, \(S_i^{(k)}\) is expected to approximate the relation in Equation~\eqref{eq:S_equation}.
% \end{remark}

\begin{remark}
In Algorithm~\ref{alg:linear}, 
% Algorithm~\ref{alg:linear} establishes the same computational framework as that of the fast direct solver for HODLR matrices, with the exception that 
$S_i^{(k+1)}$ is treated as a learnable parameter rather than being explicitly enforced via the relation
\begin{equation}
   S_i^{(k+1)} = \left(I + \widetilde{V}_{i}^{(k+1)\top}\widetilde{U}_{i}^{(k+1)}\right)^{-1}.
   \label{eq:S_equation}
\end{equation}
We choose the learnable formulation since enforcing the above relation was observed to result in higher training times. After training, $S_i^{(k+1)}$ is expected to approximate the relation in Equation~\eqref{eq:S_equation}. 
\end{remark}

\section{Neural network as PDE solver}\label{sec:PDE_solver}
% Consider an elliptic PDE $\mathcal{N}[u;\,\mu] = 0$, where $\mathcal{N}$ is a an elliptic differential operator, $u$ is the solution and $\mu$ encodes the problem data such as the problem parameters, boundary condition, forcing function, etc. 
Consider a PDE $\mathcal{D}[u;\,\mu] = f$ on a domain $\Omega$, subject to 
appropriate boundary and initial conditions, where $\mathcal{D}$ is a 
differential operator, $u$ is the solution, $f$ is the forcing function, and $\mu$ 
encodes the problem parameters such as material coefficients.
% Let $\mathcal{S}:\bkt{\mu,f} \mapsto u$ denote the associated solution operator.

\subsection{Solution operators arising in linear elliptic PDEs}\label{ssec:linear}
% Consider a solution operator $\mathcal{S}$ arising from a linear elliptic PDE. We consider its integral equation formulation using the associated Green's function $g(\mathbf{x},\mathbf{y})$. 
% Let $\mathcal{S}:f \mapsto u$ be a solution operator arising from a linear elliptic PDE. 
Consider a linear elliptic PDE with fixed $\mu$ and solution operator $\mathcal{S}:f \mapsto u$.
We consider its integral-equation formulation based on the associated Green's function $g(\mathbf{x},\mathbf{y})$.
The discretization of the integral equation using an $N$ point grid over $\Omega$, results in a linear system $Ku = f$. The matrix \(K\) can be efficiently approximated by a HODLR matrix $K^{(\kappa)}$, since the Green's function of an elliptic operator with smooth coefficients is asymptotically smooth away from the singularity \cite{bebendorf2008hierarchical,hackbusch1999sparse}. Therefore, the solution operator $f \mapsto u \approx K^{(\kappa)^{-1}}f$ can be learnt using the architecture described in Algorithm~\ref{alg:linear}.

\subsection{Non-linear solution operators}\label{ssec:nonlinear}
\begin{figure}[htbp]
\newlength{\panelhtB}\setlength{\panelhtB}{8.2cm}
\centering
\resizebox{1.2\linewidth}{!}{
\begin{tikzpicture}

% --- Subfigure (a): K_kappa^{-1} leaf block internals ---
\node[inner sep=0pt, anchor=south] (subNew) at (0,0) {%
\begin{subfigure}[b]{0.30\textwidth}
\centering
\begin{minipage}[c][\panelhtB][c]{\linewidth}
\centering
\resizebox{\linewidth}{!}{%
\begin{tikzpicture}[
    >=Latex,
    node distance=0.55cm,
    every node/.style={font=\small},
    io/.style    ={rectangle, rounded corners=3pt, draw, fill=gray!12,
                   minimum width=5.4cm, minimum height=0.85cm, align=center},
    lclin/.style ={rectangle, draw, fill=blue!25,
                   minimum width=5.4cm, minimum height=1.0cm, align=center},
    lcrelu/.style={rectangle, draw, fill=orange!25,
                   minimum width=5.4cm, minimum height=1.0cm, align=center},
    nlmark/.style={dashed, draw=gray!60!black, very thick,
                   rounded corners=3pt, inner sep=2.5pt},
    sidelbl/.style={font=\scriptsize\itshape}
]
\node (input)   [io] {$f \in \mathbb{R}^{N}$ \\ (forcing function)};
\node (lc1)     [lcrelu, above=of input]
        {LC layer 1 - non-linear};
\node (lcdots1) [above=0.40cm of lc1, font=\Large] {$\vdots$};
\node (lcmid)   [lcrelu, above=0.40cm of lcdots1]
        {LC layer $j$ - non-linear};
\node (lcdots2) [above=0.40cm of lcmid, font=\Large] {$\vdots$};
\node (lcK)     [lclin, above=0.40cm of lcdots2]
        {LC layer $D$ - linear};
\node (output)  [io, above=of lcK] {$\xi^{(\kappa)} \in \mathbb{R}^{N}$ \\ (output)};

\draw[->, thick] (input)   -- (lc1);
\draw[->, thick] (lc1)     -- (lcdots1);
\draw[->, thick] (lcdots1) -- (lcmid);
\draw[->, thick] (lcmid)   -- (lcdots2);
\draw[->, thick] (lcdots2) -- (lcK);
\draw[->, thick] (lcK)     -- (output);

% \begin{scope}[on background layer]
%   \node[nlmark, fit=(lc1)(lcdots1)(lcmid)(lcdots2)(lcK)] (lcbox) {};
% \end{scope}

\draw [opacity=0, decorate, decoration={brace, amplitude=5pt, raise=2pt}, thick]
  ($(lcK.north east)+(0.35,0)$) -- ($(lc1.south east)+(0.35,0)$)
  node[midway, sidelbl, right=8pt, align=left, text width=2cm, font=\small]
       {Stack of $D$ LC layers with non-linear activations except the top layer};
\end{tikzpicture}}
\end{minipage}
% \caption{$K_\kappa^{-1}$ leaf block internals.}
\label{fig:FDSNet-nl-leaf}
\end{subfigure}%
};

% --- Subfigure (b): high-level FDSNet ---
\node[inner sep=0pt, anchor=south, xshift=-1cm] (subA) at ($(subNew.south east)+(3cm,0)$) {%
\begin{subfigure}[b]{0.34\textwidth}
\centering
\begin{minipage}[c][\panelhtB][c]{\linewidth}
\centering
\resizebox{\linewidth}{!}{%
\begin{tikzpicture}[
    >=Latex,
    node distance=0.85cm,
    every node/.style={font=\small},
    io/.style    ={rectangle, rounded corners=3pt, draw, fill=gray!12,
                   minimum width=7.6cm, minimum height=0.85cm, align=center},
    leaf/.style  ={rectangle, draw, fill=purple!20,
                   minimum width=7.6cm, minimum height=1.1cm, align=center},
    smw/.style   ={rectangle, draw, fill=teal!45,
                   minimum width=7.6cm, minimum height=1.1cm, align=center},
    nlmark/.style={dashed, draw=gray!60!black, very thick,
                   rounded corners=3pt, inner sep=2.5pt},
    dimlbl/.style={font=\scriptsize, anchor=east},
    sidelbl/.style={font=\scriptsize\itshape}
]
\node (input)  [io] {$f\in\mathbb{R}^{N}$ \quad (forcing function)};
\node (kleaf)  [leaf,  above=of input]
       {\scriptsize $K_\kappa^{-1}$ Block\\  Stack of $D$ LC layers with non-linear activations \\ except the top layer};
\node (smwk1)  [smw,   above=of kleaf] {SMW Block $\kappa-1$};
\node (sdots)  [above=0.55cm of smwk1, font=\Large] {$\vdots$};
\node (smwj)   [smw,   above=0.55cm of sdots] {SMW Block $k$};
\node (sdots2) [above=0.55cm of smwj,  font=\Large] {$\vdots$};
\node (smw0)   [smw,   above=0.55cm of sdots2] {SMW Block 0};
\node (output) [io,    above=of smw0]  {$u$ \quad (predicted solution)};

\draw[->, thick] (input)  -- node[dimlbl, xshift=-0.15cm] {}    (kleaf);
\draw[->, thick] (kleaf)  -- node[dimlbl, xshift=-0.15cm] {$\xi^{(\kappa)}$} (smwk1);
\draw[->, thick] (smwk1)  -- (sdots);
\draw[->, thick] (sdots)  -- (smwj);
\draw[->, thick] (smwj)   -- (sdots2);
\draw[->, thick] (sdots2) -- (smw0);
\draw[->, thick] (smw0)   -- node[dimlbl, xshift=-0.15cm] {$\xi^{(0)}$}   (output);

% \begin{scope}[on background layer]
%   \node[nlmark, fit=(kleaf)]                                  (nl1) {};
%   \node[nlmark, fit=(smwk1)(sdots)(smwj)(sdots2)(smw0)]       (nl2) {};
% \end{scope}

\draw [opacity=0, decorate, decoration={calligraphic brace, amplitude=6pt, raise=2pt}, very thick]
  ($(smw0.north east)+(0.35,0)$) -- ($(smwk1.south east)+(0.35,0)$)
  node[midway, sidelbl, right=10pt, align=left, text width=2cm, font=\small]
       {$\kappa$-residual SMW levels ($k=\kappa{-}1,\dots,0$)
       \\encoding $K_0^{-1}K_1^{-1}$ \\ $\!\cdots K_{\kappa-1}^{-1}$};

% \node[sidelbl, right=0.55cm of kleaf, align=left, text width=2cm, font=\small]
%   {leaf matrix factor\\encoding $K_\kappa^{-1}$};

% \node[font=\normalsize, black!75, align=center,
%       below=0.55cm of input.south, text width=9cm] (legend)
%   {Note: dashed boxes contain an internal $\times K$ stack of LC layers with ReLU on middle layers and linear first/last.};
\end{tikzpicture}}
\end{minipage}
% \caption{FDSNet architecture.}
\label{fig:FDSNet-nl-highlevel}
\end{subfigure}%
};

% --- Subfigure (c): SMW_Block_k ---
\node[inner sep=0pt, anchor=south, xshift=-1cm] (subB) at ($(subA.south east)+(3cm,0)$) {%
\begin{subfigure}[b]{0.32\textwidth}
\centering
\begin{minipage}[c][\panelhtB][c]{\linewidth}
\centering
\resizebox{\linewidth}{!}{%
\begin{tikzpicture}[
    >=Latex,
    node distance=0.55cm,
    every node/.style={font=\small},
    io/.style       ={rectangle, rounded corners=3pt, draw, fill=gray!12,
                      minimum width=9.4cm, minimum height=0.8cm, align=center},
    layer/.style    ={rectangle, draw, fill=blue!25,
                      minimum width=9.4cm, minimum height=1.0cm, align=center},
    swaplayer/.style={rectangle, draw, fill=yellow!50,
                      minimum width=9.4cm, minimum height=1.0cm, align=center},
    nlstack/.style  ={rectangle, draw, fill=purple!20,
                      minimum width=9.4cm, minimum height=1.25cm, align=center},
    op/.style       ={circle, draw, fill=red!40, font=\large,
                      inner sep=1pt, minimum size=0.6cm},
    dimlbl/.style   ={font=\scriptsize, anchor=west}
]
\node (input)  [io] {Input $\xi^{(k+1)} \in \mathbb{R}^{N}$};
\node (swap)   [swaplayer, above=of input]
       {\textbf{Swap layer} \\ \scriptsize pairwise permute of sibling blocks};
\node (vlayer) [layer, above=of swap]
       {\textbf{$V^{(k)}$ layer} - LC linear};
\node (sstack) [nlstack, above=of vlayer, align=center]
       {\textbf{$S^{(k)}$ layer - LC stack }\\[-1pt]
        Stack of $D$ LC layers with non-linear activations \\ except the top layer};
\node (ulayer) [layer, above=of sstack]
       {\textbf{$U^{(k)}$ layer} - LC linear};
\node (sub)    [op, above=of ulayer] {$-$};
\node (output) [io, above=of sub] {Output $\xi^{(k)} \in \mathbb{R}^{N}$};

\draw[->, thick] (input)  -- node[dimlbl, xshift=0.15cm] {}                       (swap);
\draw[->, thick] (swap)   -- node[dimlbl, xshift=0.15cm] {}              (vlayer);
\draw[->, thick] (vlayer) -- node[dimlbl, xshift=0.15cm] {}                  (sstack);
\draw[->, thick] (sstack) -- node[dimlbl, xshift=0.15cm] {}                  (ulayer);
\draw[->, thick] (ulayer) -- node[dimlbl, xshift=0.15cm] {}        (sub);
\draw[->, thick] (sub)    -- (output);

\coordinate (skipX)  at ($(input.west) + (-1.0,0)$);
\coordinate (skipBot) at (skipX |- input.west);
\coordinate (skipTop) at (skipX |- sub.west);
\draw[->, thick, red!60!black]
     (input.west) -- (skipBot)
                  -- (skipTop)
                  -- (sub.west)
     node[pos=1.0, font=\scriptsize, black, xshift=-65pt, yshift=4pt, rotate=0, anchor=south, font=\small]
        {Residual connection};

% \draw [decorate, decoration={brace, amplitude=5pt, raise=2pt}, thick]
%   ($(ulayer.north east)+(0.35,0)$) -- ($(swap.south east)+(0.35,0)$)
%   node[midway, font=\small\itshape, right=8pt, align=left, text width=1cm]
%        {SMW residual update};
       % \\$I-\widetilde{U}^{(k)}_i W^{(k)}_i \widetilde{V}^{(k)\top}_i$
\end{tikzpicture}}
\end{minipage}
% \caption{Detailed SMW\_Block$_k$.}
\label{fig:FDSNet-nl-smwblock}
\end{subfigure}%
};

% --- Wedges ---
\begin{scope}[on background layer]
  \coordinate (w1WideTop)   at ($(subNew.south east)!0.78!(subNew.north east) + (-1.3cm,0)$);
  \coordinate (w1WideBot)   at ($(subNew.south east)!0.23!(subNew.north east) + (-1.3cm,0)$);
  \coordinate (w1NarrowTop) at ($(subA.south west)!0.31!(subA.north west) + (0.10cm,0)$);
  \coordinate (w1NarrowBot) at ($(subA.south west)!0.29!(subA.north west) + (0.19cm,0)$);
  \fill[blue!50!cyan, opacity=0.09]
      (w1WideTop) -- (w1NarrowTop) -- (w1NarrowBot) -- (w1WideBot) -- cycle;
  \draw[blue!50!cyan!70, dashed, thick] (w1WideTop) -- (w1NarrowTop);
  \draw[blue!50!cyan!70, dashed, thick] (w1WideBot) -- (w1NarrowBot);
\end{scope}

\begin{scope}[on background layer]
  \coordinate (narrowTop) at ($(subA.south east)!0.55!(subA.north east) + (-1.4cm,0)$);
  \coordinate (narrowBot) at ($(subA.south east)!0.55!(subA.north east) + (-1.30cm,0)$);
  \coordinate (wideTop)   at ($(subB.south west)!0.80!(subB.north west) + (0.4cm,0)$);
  \coordinate (wideBot)   at ($(subB.south west)!0.20!(subB.north west) + (0.4cm,0)$);
  \fill[blue!50!cyan, opacity=0.09]
      (narrowTop) -- (wideTop) -- (wideBot) -- (narrowBot) -- cycle;
  \draw[blue!50!cyan!70, dashed, thick] (narrowTop) -- (wideTop);
  \draw[blue!50!cyan!70, dashed, thick] (narrowBot) -- (wideBot);
\end{scope}
\path ([xshift=3cm]current bounding box.east);
\end{tikzpicture}
}
\vspace{-1.5cm}
\caption{Non-linear FDSNet: A neural network to learn non-linear solution operators associated with PDEs}
% \caption{Nonlinear FDSNet architecture (Algorithm~2), drawn bottom-to-top. (a) Leaf $K_\kappa^{-1}$ block internals: a stack of $D$ LC layers with $k\!=\!s\!=\!f\!=\!m$, $m\!=\!N/2^{\kappa}$; first and last layers linear (blue), middle layers ReLU (orange). (b) Sequential composition: a leaf $K_\kappa^{-1}$ block followed by $\kappa$ residual SMW blocks realizing $K^{-1}=K_0^{-1}\!\cdots K_\kappa^{-1}$. (c) One SMW level implementing the residual update $\xi^{(k)}=\xi^{(k+1)}-\widetilde{U}^{(k)}_iW^{(k)}_i\widetilde{V}^{(k)\top}_i\xi^{(k+1)}$; the dashed box marks the $\times D$ $S^{(k)}$-layer stack with ReLU on middle layers and linear first/last.}
\label{fig:FDSNet-nonlinear}
\end{figure}
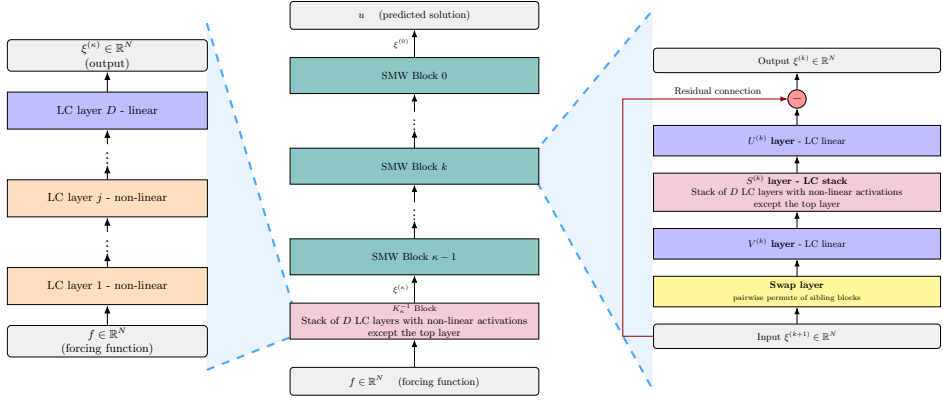

Consider a PDE with non-linear solution operator $\mathcal{S}:\bkt{\mu,f} \mapsto u$. Its solution in general cannot be represented through the application of a fixed matrix inverse.
% For a PDE $\mathcal{D}[u;\,\mu] = f$, with non-linear solution operator $\mathcal{S}:\bkt{\mu,f} \mapsto u$, the solution in general cannot be represented through the application of a fixed matrix inverse.
Therefore we extend Algorithm~\ref{alg:linear}
to Algorithm~\ref{alg:non_linear} to learn the non-linear solution operator. The $V^{(k)}$ and $U^{(k)}$ layers are retained as single linear maps.
These layers are considered to encode the geometry of the hierarchical partition, which is operator-independent. The operator dependence is placed entirely on the \(K^{(\kappa)}\) and \(S^{(k)}\) layers, where the single LC (locally connected) layer with linear activation in Algorithm~\ref{alg:linear} is replaced by a stack of \(D\) LC layers with nonlinear activations (e.g., ReLU) except the top layer.
% Each LC layer in the stack is an MLP  from $\mathbb{R}^{2p}$ to $\mathbb{R}^{2p}$, giving the architecture the capacity to learn the non-linear operator. 
Each LC layer in the stack is implemented as a fully connected layer.
% , mapping \(\mathbb{R}^{2p}\) to \(\mathbb{R}^{2p}\). 
We refer to the neural network described in Algorithm~\ref{alg:non_linear} as non-linear FDSNet. 
% Following a similar analysis as that for the linear FDSNet, the number of parameters in the nonlinear FDSNet scales as $\mathcal{O}\bkt{pN\log\bkt{N} + p^2DN}$.
Figure~\ref{fig:FDSNet-nonlinear} illustrates the architecture  described in Algorithm~\ref{alg:non_linear}.
The effectiveness of this network is demonstrated on a range of PDE benchmark problems in Section~\ref{sec:numerical}.

\begin{algorithm}[htbp]
\caption{Neural network to learn non-linear solution operators associated with PDEs.}
% \caption{Application of Direct Solver NN of inverse of an $\mathcal{H}$-matrix to a vector $v \in \mathbb{R}^N$.}
\label{alg:non_linear}
\begin{algorithmic}[1]
\REQUIRE The input vector $b\in \mathbb{R}^{N}$
\REQUIRE The number of levels $\kappa$ in HODLR tree structure
\STATE $\xi^{(\kappa+1)}_0= b;$
% \FOR{ $i = 1$ to $2^{\kappa}$ }
    % \STATE $\xi^{(\kappa+1)}_i= b_i;$
% \ENDFOR
\FOR{ $j = 0$ to $D-1$ }
    \IF{$j < D-1$}
        \STATE $\xi^{(\kappa)}_{j+1} = \text{LC}\sbkt{\text{activation=non-linear}; s=\frac{N}{2^\kappa}; w=\frac{N}{2^\kappa}; f=\frac{N}{2^\kappa}}\bkt{\xi^{(\kappa+1)}_j};$ 
    \ELSE
        \STATE $\xi^{(\kappa)}_{j+1} = \text{LC}\sbkt{\text{activation=linear}; s=\frac{N}{2^\kappa}; w=\frac{N}{2^\kappa}; f=\frac{N}{2^\kappa}}\bkt{\xi^{(\kappa+1)}_j};$ 
    \ENDIF
\ENDFOR
        \mycomment{A stack of $D$ $K^{(\kappa)}$ layers}
\STATE $\xi^{(\kappa)} = \text{Flatten}\bkt{\xi^{(\kappa)}_D};$
\mycomment{Flatten layer}
% \STATE $\xi^{(\kappa)}  = \textsf{Flatten}(\xi^{(\kappa)});$
\FOR{ $k = \kappa-1$ to $0$ }
    % \STATE Swap sibling blocks of $\xi^{(\kappa)} = \sbkt{\xi^{(\kappa)}_1, \xi^{(\kappa)}_2, \dots, \xi^{(\kappa)}_{2^{\kappa}}}$, i.e. swap $\xi^{(\kappa)}_{2i-1}$ with $\xi^{(\kappa)}_{2i}$ $\forall\;i\in\{1,\hdots2^{\kappa-1}\}$
    \STATE Swap sibling blocks of $\xi^{(k+1)} = \sbkt{\xi^{(k+1)}_1, \xi^{(k+1)}_2, \dots, \xi^{(k+1)}_{2^{k+1}}}$, i.e. swap $\xi^{(k+1)}_{2i-1}$ with $\xi^{(k+1)}_{2i}$, $\forall\;i\in\{1,\hdots2^{k}\}$

    \mycomment{Swap layer}
    \STATE $\eta^{(k)}_0 = \text{LC}\sbkt{\text{activation=linear}; s=\frac{N}{2^{k+1}}; w=\frac{N}{2^{k+1}}; f=p}\bkt{\xi^{(k+1)}};$ 
    \mycomment{$V^{(k)}$ layer}
    % \STATE $\eta^{(k)}_0 = $
    % \STATE \mycomment{Apply a stack of $D$ $W^{(k)}$ layers with non-linear activations, i.e.,}
    \FOR{ $j = 0$ to $D-1$ } 
        \IF{$j < D-1$}
            \STATE $\eta^{(k)}_{j+1} = \text{LC}\sbkt{\text{activation=non-linear}; s=2; w=2; f=2p}\bkt{\eta^{(k)}_j};$ 
        \ELSE
            \STATE $\eta^{(k)}_{j+1} = \text{LC}\sbkt{\text{activation=linear}; s=2; w=2; f=2p}\bkt{\eta^{(k)}_j};$ 
        \ENDIF
    \ENDFOR
    \mycomment{A stack of $D$ $S^{(k)}$ layers}
    \STATE $\zeta^{(k)} = \text{Reshape}\bkt{\eta^{(k)}_D,\,\bkt{2^{k+1},\,p}}$
    \mycomment{Reshape layer}
    % \STATE $\zeta^{(k)} = \eta^{(k)}_D$
    \STATE $\chi^{(k)} = \text{LC}\sbkt{\text{activation=linear}; s=1; w=1; f=\frac{N}{2^{k+1}}}\bkt{\zeta^{(k)}};$ 
    \mycomment{$U^{(k)}$ layer}
    % \STATE $\zeta^{(i)} = \textsf{Flatten}(\chi^{(i)});$
    \STATE $\xi^{(k)} = \xi^{(k+1)}  - \textsf{Flatten}\bkt{\chi^{(k)}};$
    \mycomment{Residual layer}
\ENDFOR
\STATE $u = \xi^{(0)}$
\RETURN $u$
\end{algorithmic}
\end{algorithm}

% =================================================================
% FIGURE 2 -- NONLINEAR FDSNet ARCHITECTURE (Algorithm 2)
% =================================================================

\subsection{Translation-invariant case}
\label{sec:3.2.3}
% In the one dimensional case, if the kernel is translation invariant, i.e., $g(x, y) = g(x-y)$, then the resulting matrix is a Toeplitz matrix. For Toeplitz matrices, $A[i, j] = A[i+1, j+1]$ for all $i,j \in \Omega$, where $\Omega$ is the domain. This implies that all elements along each diagonal are identical. 
\textbf{Linear Case.} For a translation-invariant kernel $g(\mathbf{x},\mathbf{y}) = g(\mathbf{x}-\mathbf{y})$ on a uniform grid, the low-rank factors $U_i^{(k)}$ and $V_i^{(k)^T}$ of the off-diagonal blocks at level $k$ in the HODLR matrix are exactly identical across all $i \in \{1, \dots, 2^k\}$. This is because, on an uniform grid, all off-diagonal blocks at the same level correspond to interactions over exactly the same set of pairwise differences between the coordinates of the discretization points, and hence share the same low-rank structure. 
% Consequently, the $\textsf{LC}$ layers at each level can be replaced by $\textsf{CNN}$ layers with a shared kernel, removing the dependence on block index $j$ while retaining a distinct kernel per hierarchical level $i$.
Further, in the factorization step of the direct solver algorithm, where the matrix $K^{(\kappa)}$ is factorized into $K_i$'s for all $i=\{0,\dots,\kappa\}$, with each factor having $2^i$ diagonal blocks in $\mathbb{R}^{N/2^i\times N/2^i}$, the translation invariance implies that all block matrices within each factor share the same entries, i.e., the blocks $K_i^{(k)}$ for all $i \in \cbkt{1, \dots, 2^k}$ are identical. 
Furthermore, \(\widetilde{V}_{i}^{(k+1)^{T}}\) are identical for all \(i \in \cbkt{1,\dots,2^k}\) and the same holds true fro \(\widetilde{U}_{i}^{(k+1)}\) and \(S_{i}^{(k+1)}\).
Consequently, within each of \(K^{(\kappa)}\), \(V^{(k)}\), \(S^{(k)}\), and \(U^{(k)}\) layers, for all \(k \in \cbkt{0,\dots,\kappa-1}\), all the local connections share common weights. This leads to a key architectural simplification: the LC layers can be replaced with CNN (convolutional neural network) layers, which use identical weights for all spatial positions.

\textbf{Non-linear case.}
In the nonlinear case, suppose the solution operator $\mathcal{S}$ is translation equivariant, i.e., for a translation $T$, $T \mathcal{S}(v) = \mathcal{S}(Tv)$.
Under a uniform grid discretization, the LC layers in Algorithm~\ref{alg:non_linear} can be replaced by CNN layers, analogous to the linear case.

The scaling of the number of parameters in the networks discussed above is summarized in Table~\ref{tab:scaling_of_parameters}.

\begin{table}[htbp]
\centering
\small
\begin{tabular}{ll}
% \begin{tabularx}{\linewidth}{XX}
\midrule
Network & Number of parameters \\[0.2em]
% \noalign{\hrule height 1pt}
\toprule 
% \hline
Linear FDSNet with LC layers & $\mathcal{O}\bkt{2pN\log\bkt{N} + p^2 N}$ \\
Linear FDSNet with CNN layers & $\mathcal{O}\bkt{2pN + p^2 \log\bkt{N}}$ \\
Non-linear FDSNet with LC layers & $\mathcal{O}\bkt{2pN\log\bkt{N} + p^2 DN}$ \\
Non-linear FDSNet with CNN layers & $\mathcal{O}\bkt{2pN + p^2 D\log\bkt{N}}$ \\
\midrule
\end{tabular}
% \end{tabularx}
\caption{Scaling of number of parameters in FDSNet}
\label{tab:scaling_of_parameters}
\end{table}

\section{Numerical results}\label{sec:numerical}
In this section, we describe the implementation details of FDSNet and demonstrate its performance through a comprehensive set of experiments.
% \begin{enumerate}
    % \item 
    % We describe the implementation details of FDSNet and demonstrate its accuracy on a comprehensive set of problems governed by partial differential and integral equations.
    % \item We compare the proposed architecture with existing neural PDE solvers, including the hierarchical multiscale neural network (MNN), the Fourier neural operator (FNO), and DeepONet.
    % \item Parameter sweep
% \end{enumerate}
The following problems are considered for the various experiments presented in this section.
% The experiments are conducted on the following four problems:
\begin{enumerate}
    \item 
    Fredholm integral equation of the second kind in 2D, 
    \item 
    Nonlinear Schr\"odinger equation (NLSE) in both 1D and 2D, 
    \item 
    Burgers' equation in 1D, and 
    \item 
    Steady-state Darcy's flow equation in 2D. 
\end{enumerate}
% Overview of implementation is provided in Section~\ref{sec:implementation}, followed by results for each test case. 

% \subsection{Implementation Overview}\label{sec:implementation}
FDSNet is implemented in Python using PyTorch~\cite{paszke2019pytorch}. For non-linear problems, rectified linear unit (ReLU)
activation~\cite{nair2010rectified} is used in the non-linear layers of $K^{(\kappa)}$ and $S^{(k)}$ stack, while a linear activation is used elsewhere.
The network is trained end-to-end by minimizing the mean squared relative
\(\ell^{2}\) error over a minibatch of size \(B\). Let
\(u_{NN}^{(i)}\) and \(u_{G}^{(i)}\) denote the network prediction and
ground-truth corresponding to the \(i\)-th sample in a minibatch of size $B$,
respectively. The training loss \(\mathcal{L}\) is defined as
\begin{equation}
\mathcal{L} = \frac{1}{B}\sum_{i=1}^{B} \epsilon_i^2,
\qquad
\text{where}
\qquad
\epsilon_i
=
\frac{
{\norm{ u_{NN}^{(i)} - u_{G}^{(i)} }}_{2}
}{
{\norm{ u_{G}^{(i)} }}_{2}
}.
\label{eq:loss}
\end{equation}
% and $\lVert \cdot \rVert_{2}$ denotes the $\ell^{2}$-norm.
We use NAdam as the optimization algorithm~\cite{dozat2016incorporating}.
The evaluation metric is defined as $\frac{1}{B}\sum_{i=1}^{B}\epsilon_i$ and is reported as the error in all subsequent experiments.
% For a single sample, the relative $\ell^{2}$ error is defined as
% \begin{equation}
%     \epsilon_i \;=\; \frac{\lVert u_{i}^{(\mathrm{NN})} - u_{i}^{(G)} \rVert_{\ell^{2}}}
%                         {\lVert u^{(G)}_i \rVert_{\ell^{2}}},
%     \label{eq:rel_l2}
% \end{equation}
%

All experiments are conducted on GPUs in 32-bit floating-point precision. The hyperparameters $p$ 
% (off-diagonal block rank in the HODLR structure),
and $D$ 
% (the number of non-linear layers in the $K^{(\kappa)}$ and $S^{(k)}$ layers~\ref{alg:non_linear}) 
are tuned per test case. 
% The reported results are the outcomes of multiple independent training runs on identical datasets across varying $(r, K)$ configurations, with seeds varied between runs. 
The source code and training scripts for all experiments reported in this paper have been made publicly available at {\url{https://github.com/JashwanthKadaru/directsolver-neural-network}}.
% \\\url{https://github.com/JashwanthKadaru/directsolver-neural-network}.

% \subsection{Results on PDE Benchmarks}
\subsection{Fredholm integral equation} \label{sec:fredholm}
The Fredholm integral equation of the second kind arises in many fields such as potential theory, radiative transfer, and scattering, and is widely used as a benchmark for solvers involving hierarchical matrices. 
We consider its two-dimensional form,
\begin{equation}
u(\mathbf{x}) + \int_{\Omega} g(\mathbf{x}, \mathbf{y}) u(\mathbf{y}) \, d\mathbf{y} = f(\mathbf{x}), \quad \mathbf{x} \in \Omega \subset \mathbb{R}^2
\label{eq:fredholm}
\end{equation}

\noindent where $g(\mathbf{x}, \mathbf{y})$ is the kernel, $f(\mathbf{x})$ is the right-hand side, and $u(\mathbf{x})$ is the unknown solution. We consider $g(\mathbf{x}, \mathbf{y})=\log\bkt{\|\mathbf{x}-\mathbf{y}\|_2}$ and $\Omega = [-1,1]^2$. Discretizing on an $n \times n$ uniform grid reduces Equation~\eqref{eq:fredholm} to a dense linear system $(I + A)\mathbf{u} = \mathbf{f}$ with $N = n^2$ degrees of freedom, where $\mathbf{u}$ and $\mathbf{f}$ are the discretized solution and right-hand side vectors respectively. The matrix \(I + A\) admits HODLR structure, and hence the proposed linear FDSNet can be applied to the problem.

\begin{table}[htbp]
\centering
\footnotesize
\begin{tabular}{ccccccc}
\midrule
$N$ & $\kappa$ & $p$ & $D$ & $m$ & Parameter count & Test error \\
\toprule
 % \multirow{4}{*}{1,600}  & 6 & 10 & 1 & 25 & 36,305  & $4.50 \times 10^{-3}$ \\
 \multirow{3}{*}{1,600}  & 6 & 12 & 1 & 25 & 43,697  & $4.63 \times 10^{-3}$ \\
                         & 6 & 14 & 1 & 25 & 51,281  & $4.25 \times 10^{-3}$ \\
                         & 6 & 16 & 1 & 25 & 59,057  & $3.73 \times 10^{-3}$ \\
\midrule
 \multirow{3}{*}{6,400}  & 8 & 10 & 1 & 25 & 137,965 & $2.96 \times 10^{-3}$ \\
                         & 8 & 12 & 1 & 25 & 164,921 & $3.49 \times 10^{-3}$ \\
                         & 8 & 14 & 1 & 25 & 192,133 & $3.41 \times 10^{-3}$ \\
                         % & 8 & 16 & 1 & 25 & 219,601 & $3.05 \times 10^{-3}$ \\
\midrule
 % \multirow{4}{*}{14,400} & 6 & 10 & 1 & 225 & 351,105 & $4.27 \times 10^{-3}$ \\
 \multirow{3}{*}{14,400} & 6 & 12 & 1 & 225 & 408,897 & $4.35 \times 10^{-3}$ \\
                         & 6 & 14 & 1 & 225 & 466,881 & $5.03 \times 10^{-3}$ \\
                         & 6 & 16 & 1 & 225 & 525,057 & $3.54 \times 10^{-3}$ \\
\midrule
\end{tabular}
\caption{Results for Fredholm integral equation using FDSNet}
\label{tab:fredholm_results}
\end{table}

\paragraph{Experimental Setup \& Results}
We learn the mapping $f\rightarrow u$ using the linear FDSNet. 
Since the kernel is translation invariant and the discretization grid is uniform, the LC layers are replaced with CNN layers.
% In each case, the domain \(\Omega\) is discretized  with \(N\) nodes. 
The dataset consists of 10000 \((f,u)\) pairs, with 5000 samples used for training and 5000 for testing. Training is performed using NAdam for a maximum of 2000 epochs, learning rate $1.0\times 10^{-3}$, batch size $128$ together with early stopping using a patience of 150 epochs.
The ground-truth solutions are generated by discretizing Equation~\eqref{eq:fredholm} using the Nyström method, followed by directly solving the resulting discrete linear systems.
FDSNet is evaluated on three grid resolutions. 
Table~\ref{tab:fredholm_results} lists (i) the hyperparameters for each problem size: $\kappa$ (number of levels in the HODLR structure), $p$ (off-diagonal block rank), $D$ (depth of the $S^{(k)}$ and $K^{(\kappa)}$ stacks), $m$ (leaf size), and the model parameter count and (ii) the test error (relative \(\ell^{2}\) error). It can be observed that FDSNet learns the $f\rightarrow u$ map well across all three resolutions, achieving test errors at or below \(0.50\%\), with a best-case error of nearly \(0.30\%\) at \(N = 6400\).

\subsection{Nonlinear Schr\"odinger Equation (NLSE 1D \& 2D)}\label{sec:nlse}

The nonlinear Schr\"odinger equation (NLSE) governs ground states of trapped Bose--Einstein condensates and a range of nonlinear-optical phenomena~\cite{sulem1999}. 
We consider the stationary defocusing nonlinear Schr\"odinger equation with 
an inhomogeneous potential $V(x)$,
\begin{equation}
    -\Delta u(x) + V(x)\,u(x) + \beta\,|u(x)|^{2}u(x) \;=\; E\,u(x),
    \qquad x \in \Omega,
    \label{eq:nlse}
\end{equation}
on $\Omega = [0,1)^{d}$, $d \in \{1,2\}$, with periodic boundary conditions, 
where $E \in \mathbb{R}$ is the eigen value (chemical potential), subject to 
the normalization $\int_{\Omega} u^{2}\,\mathrm{d}x = 1$. We fix $\beta = 10$, 
corresponding to a strongly defocusing regime, in which the ground state 
solution $u_{G}$ is real and positive. 
% Generalization across $\beta$ is studied separately in Section~\ref{sec:beta_generalization}.

The objective is to learn the solution operator $\mathcal{S}\colon V(\cdot) \mapsto u_{G}$. 
Owing to the periodic boundary conditions and the translation-invariant structure of the differential operators in Equation~\eqref{eq:nlse}, the solution operator $\mathcal{S}$ is translation-equivariant, and its discretization on a uniform grid exhibits a convolutional structure in both 1D and 2D. We exploit this property by replacing the locally-connected blocks of the non-linear FDSNet with shared convolutional kernels.
% When the domain is discretized on a uniform grid with periodic boundary conditions, $\mathcal{M}$ is translation-equivariant. 
% The discretization of \(V\) on a uniform grid, combined with periodic boundary conditions, ensures that \(\mathcal{M}\) is translation-equivariant. We exploit this by replacing the locally-connected blocks of FDSNet with shared convolutional kernels.
% The learning task is the operator $\mathcal{M}\colon V \mapsto u^{(G)}$. 
% With periodic boundary conditions and a translation-invariant data distribution, $\mathcal{M}$ is translation-equivariant and the discretised system carries a Toeplitz (block-Toeplitz in 2D) structure. 

% \paragraph{Dataset.} The potential is a periodic sum of Gaussian
% wells
% \[
% V(x) = -\sum_{i=1}^{n_{g}} \sum_{j\in\mathbb{Z}}
% \frac{\rho^{(i)}}{\sqrt{2\pi T}}
% \exp\!\Big(-\tfrac{|x - j - c^{(i)}|^{2}}{2T}\Big).
% \]
% Sampling parameters are $n_{g}\in\{1,\dots,4\}$, $\rho\sim U(1,4)$,
% $c\sim U(0,1)$, $T\sim U(2,4)\!\times\!10^{-3}$.
% The ground states are computed using normalized gradient flow (imaginary-time propagation)~\cite{baodu2004}.We use 20{,}000 training and 20{,}000 test pairs
% $(V, u^{(G)})$ in both dimensions.

\paragraph{Dataset} The potential $V(x)$ is a periodic sum of Gaussian wells,
\[
V(x) = -\sum_{i=1}^{n_{g}} \sum_{j\in\mathbb{Z}}
\frac{\rho^{(i)}}{\sqrt{2\pi\sigma^2}}
\exp\!\Big(-\frac{|x - j - c^{(i)}|^{2}}{2\sigma^2}\Big),
\]
where $n_{g} \in \{1,\dots,4\}$ is the number of Gaussian wells, 
$\rho^{(i)} \sim U(1,4)$ are the amplitudes, $c^{(i)} \sim U(0,1)$ are 
the centers, and $\sigma^2 \sim U(2,4)\times 10^{-3}$ is the variance. 
Since the centers $c^{(i)}$ are drawn uniformly on $[0,1)^d$, the 
distribution of $V$ is translation-invariant on the torus, consistent 
with the periodic boundary conditions. 
We use $20{,}000$ training and 
$20{,}000$ test pairs $(V, u_{G})$ in both dimensions and train for 2000 maximum epochs with a learning rate of $1\times10^{-3}$. The ground-truth solutions are computed using normalized gradient flow~\cite{baodu2004}.

\subsubsection{NLSE 1D}\label{sec:nlse_1d} 
\begin{figure}[H]
  \centering
  \begin{subfigure}[t]{0.45\linewidth}
    \centering
    \includegraphics[width=0.8\linewidth]{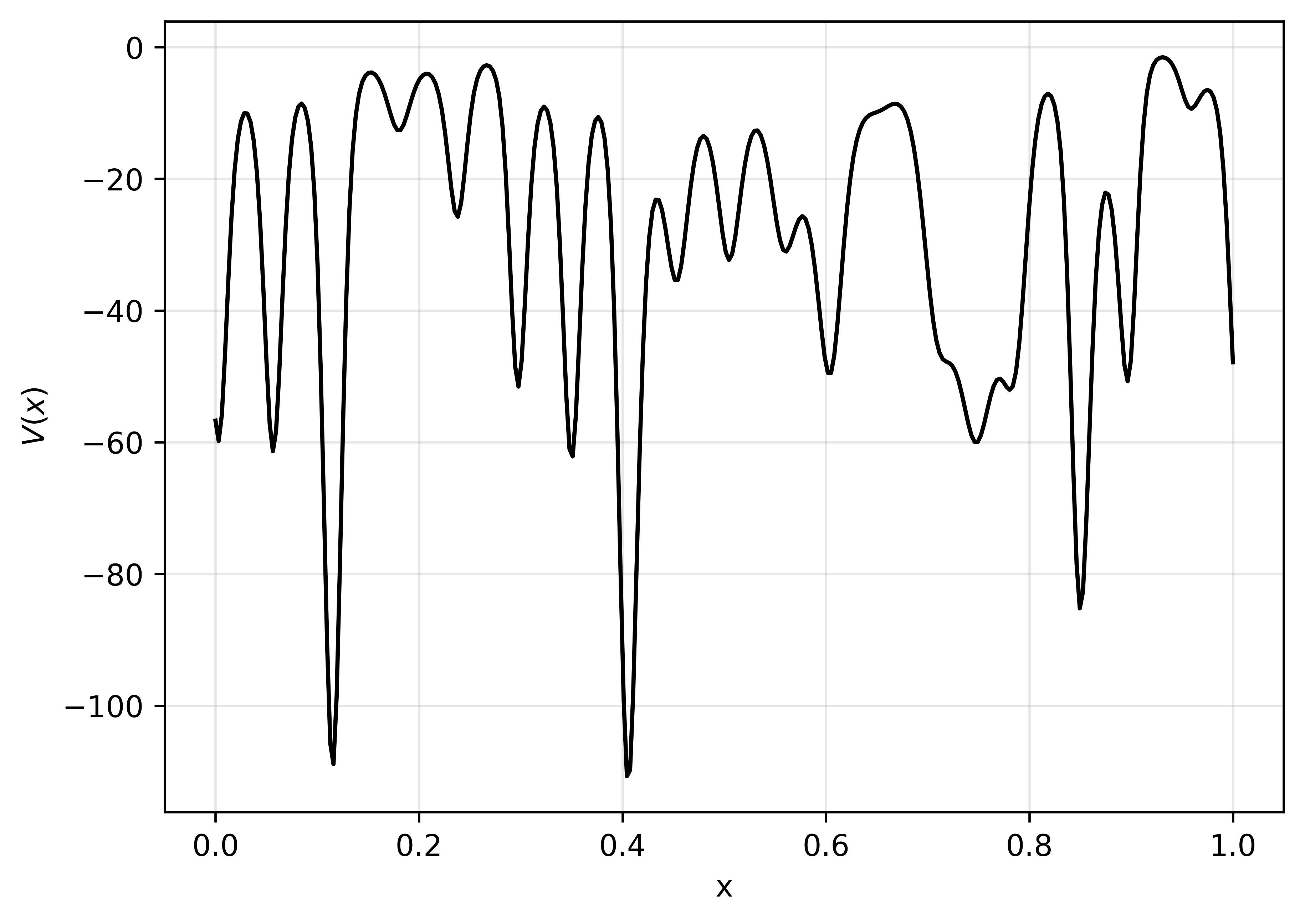}
    \caption{$V$}
  \end{subfigure}\hfill
  \begin{subfigure}[t]{0.45\linewidth}
    \centering
    \includegraphics[width=0.8\linewidth]{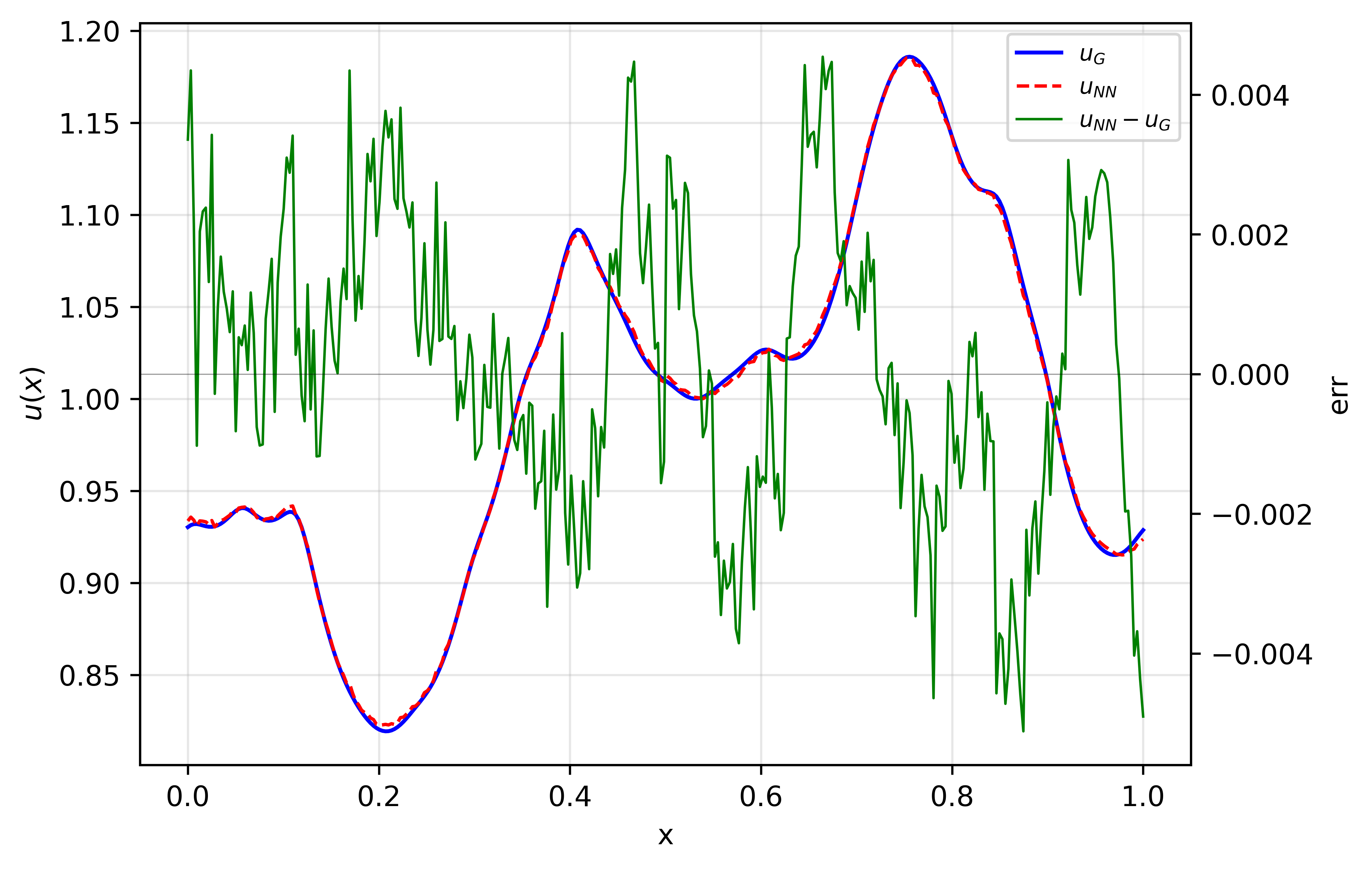}
    \caption{$u_G$, $u_{NN}$, and $u_{NN}-u_G$}
  \end{subfigure}
  \caption{NLSE 1D: FDSNet prediction $u_{NN}$, ground truth $u_G$, and the error function $u_{NN} - u_G$ for a test sample $V$. The test error is $2.10\times 10^{-3}$ with $N=320$, $\kappa=6$, $p=10$, $D=5$.}
  \label{fig:nlse1d_fdsnet_prediction}
\end{figure}

The domain is discretised using a uniform grid with $N = 320$ points. A hierarchical depth of $\kappa = 6$ is used in the FDSNet configuration, resulting in a leaf size of $m = 5$. We experimented with varying depths $D$ of $K^{(\kappa)}$ and $S^{(k)}$ stacks and ranks $p$. The network is trained using the NAdam optimiser with batch size $128$ and early stopping with patience $200$. 
An illustration of a sample of $V$ from the test data, the ground truth $u_G$, the prediction $u_{NN}$ by the FDSNet, and the error with respect to $u_G$ is given in Figure~\ref{fig:nlse1d_fdsnet_prediction}.
The performance of FDSNet is tabulated in Table~\ref{tab:FDSNet_sweep_all}. 
Increasing the rank from $2$ to $10$ reduces the error by $\approx 75\%$ at $\approx 8\times$ more parameters, while the depth sweep is non-monotone, with $D = 5$ outperforming both $D = 3$ and $D = 7$. 

% for up to $2000$ epochs
% The best epoch in training is determined by minimizing the validation loss and the epoch's checkpoint is used for inference and metrics.   

% ---- NLSE 1D: best test sample (lowest rel-l2 error) ----

\begin{figure}[htbp]
  \centering
  \begin{subfigure}[t]{0.45\linewidth}
    \centering
    \includegraphics[width=0.8\linewidth]{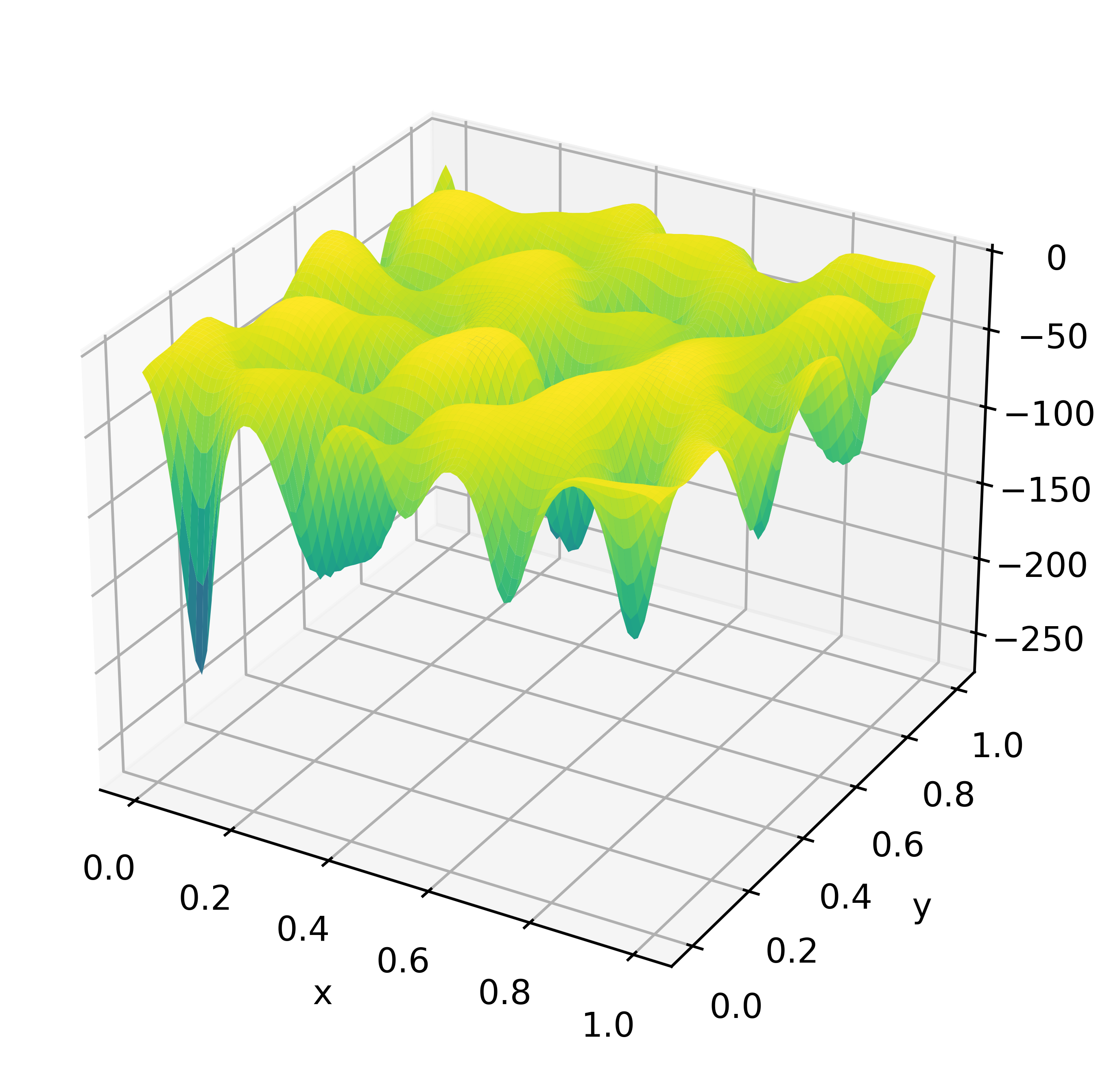}
    \caption{$V$}
  \end{subfigure}\hfill
  \begin{subfigure}[t]{0.45\linewidth}
    \centering
    \includegraphics[width=0.8\linewidth]{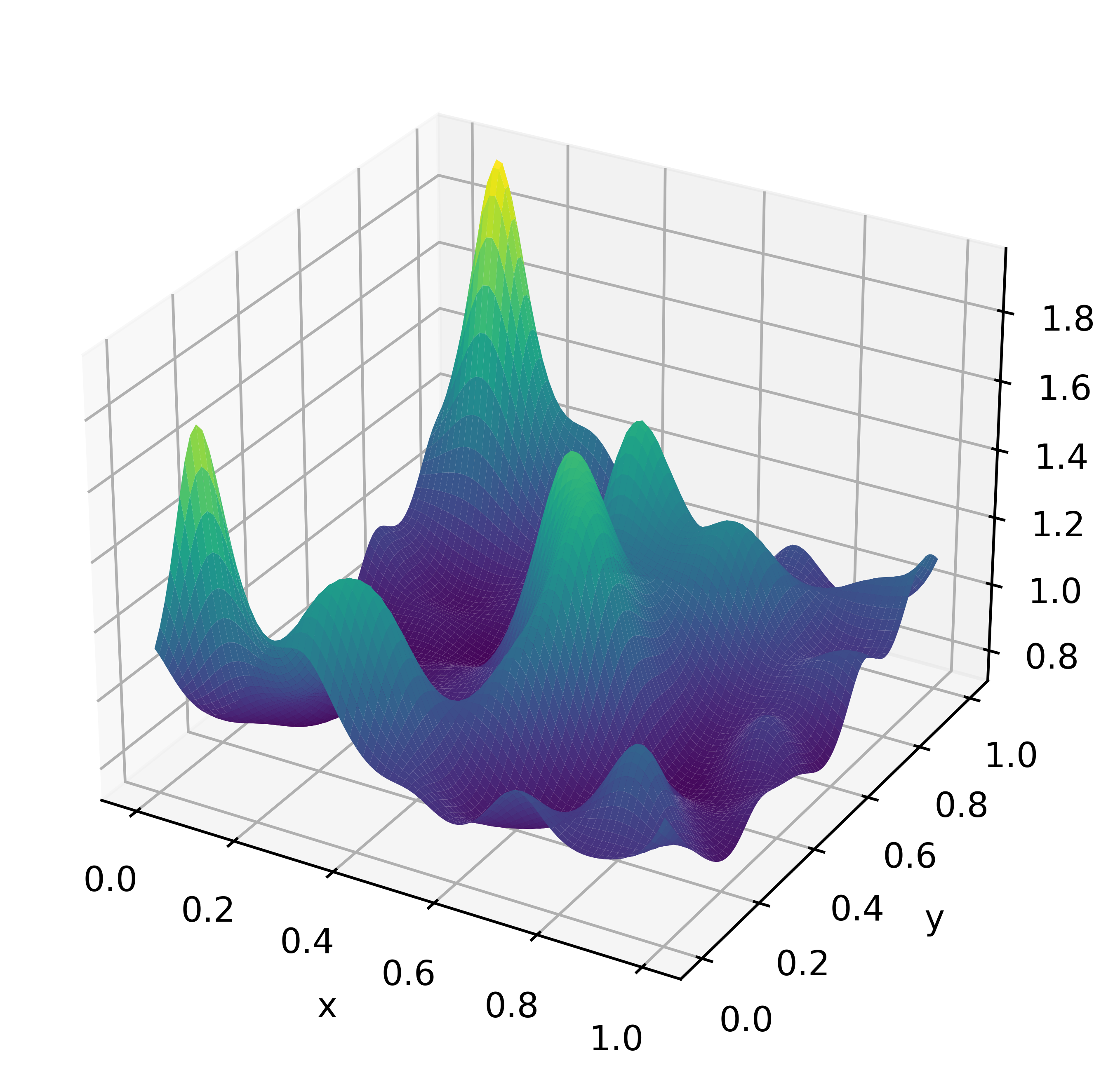}
    \caption{$u_G$}
  \end{subfigure}

  \begin{subfigure}[t]{0.45\linewidth}
    \centering
    \includegraphics[width=0.8\linewidth]{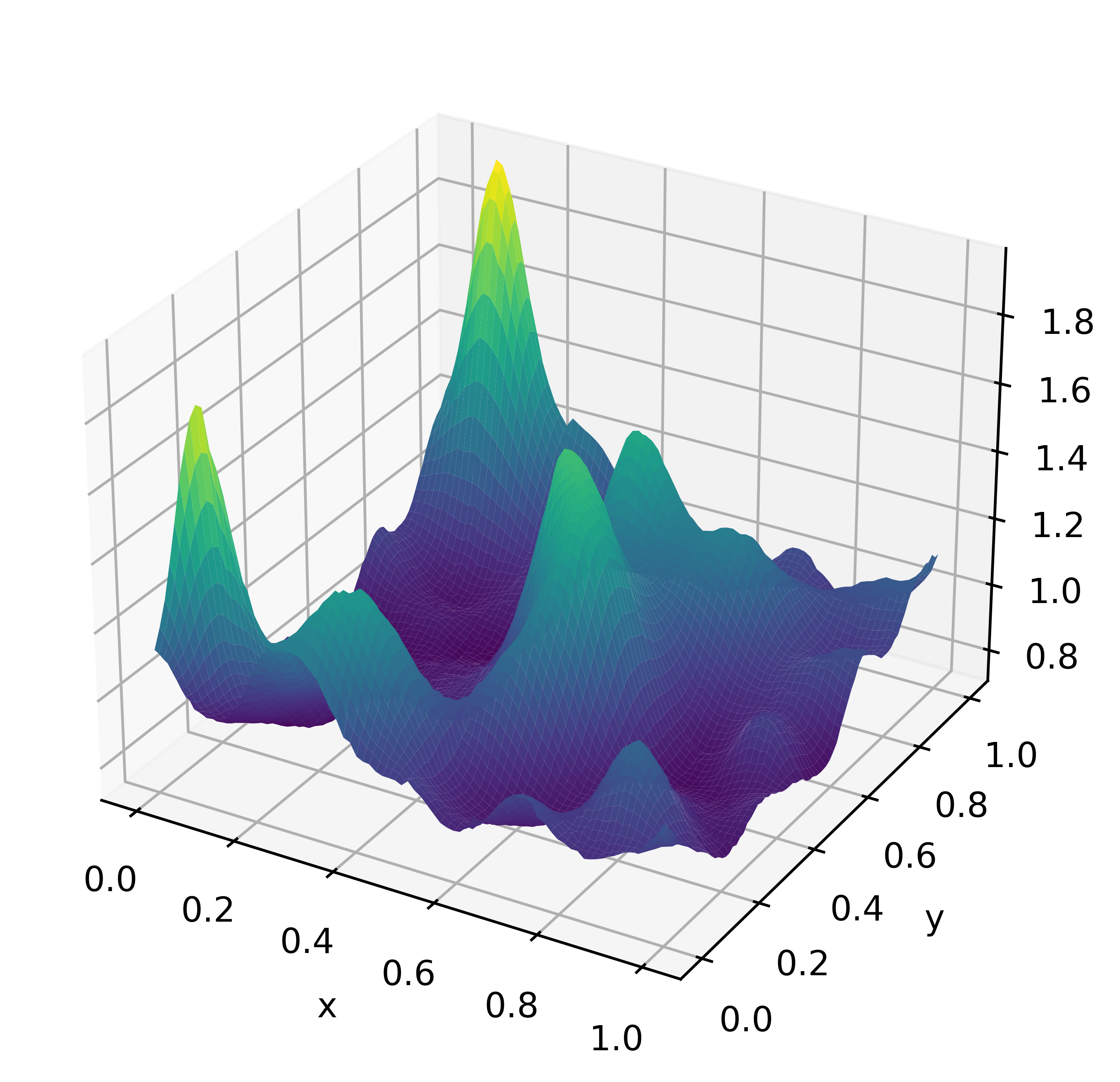}
    \caption{$u_{NN}$}
  \end{subfigure}\hfill
  \begin{subfigure}[t]{0.45\linewidth}
    \centering
    \includegraphics[width=0.8\linewidth]{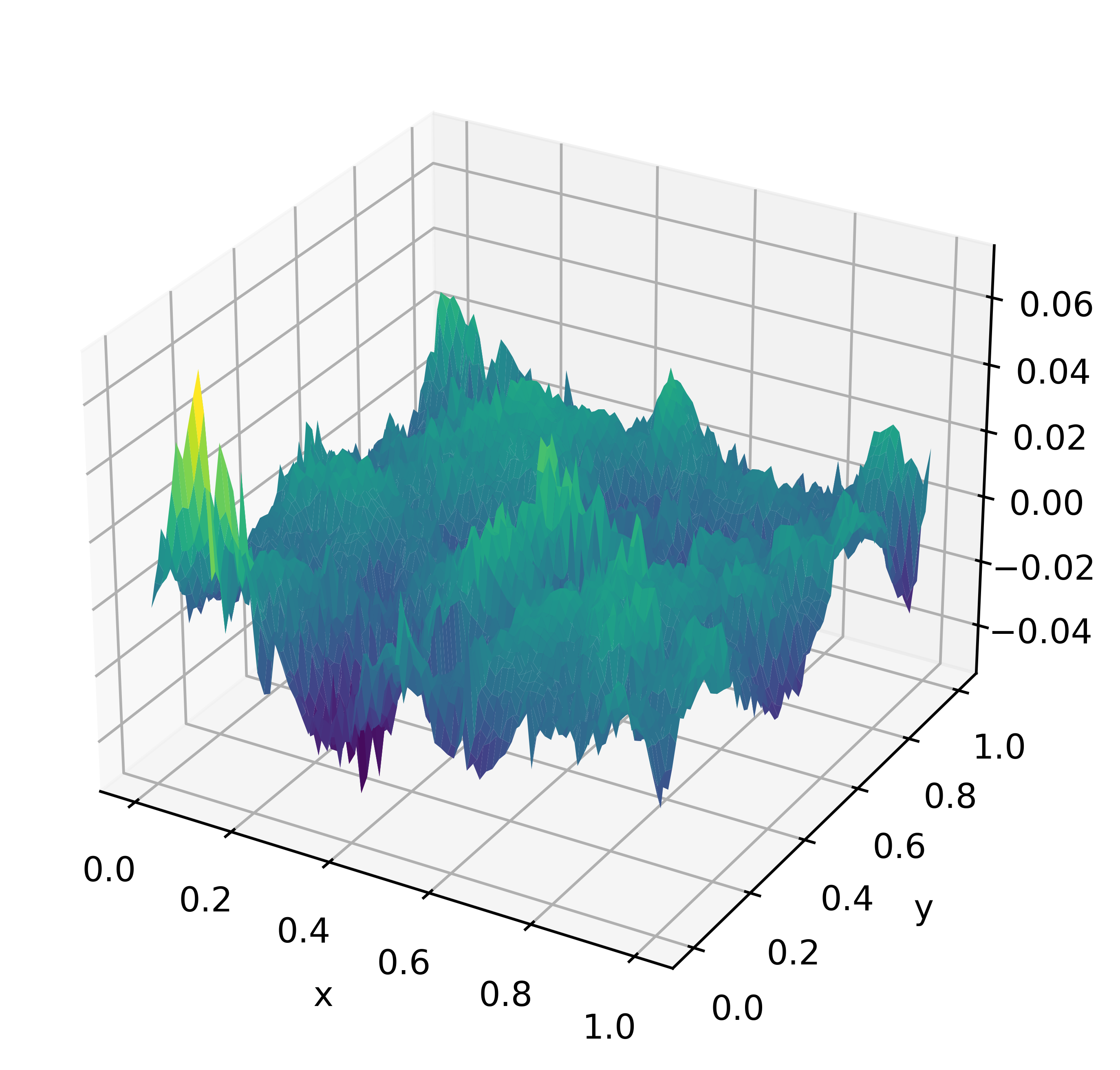}
    \caption{$u_{NN}-u_G$}
  \end{subfigure}

  \caption{NLSE 2D: FDSNet prediction $u_{NN}$, ground truth $u_G$, and the error function $u_{NN} - u_G$ for a test sample $V$. The test error is $1.03\times 10^{-2}$ with $N = 6400$, $\kappa=8$, $p=10$, $D=5$.}
  \label{fig:nlse2d_fdsnet_prediction}
\end{figure}

\subsubsection{NLSE 2D}\label{sec:nlse_2d} 
The domain is discretised on an $N_{x} \times N_{y}$ 2D uniform grid with $N_{x} = N_y = 80$, resulting in $N = 6400$ grid points. 
% $\beta$ is set to 10. 
FDSNet is configured to $\kappa = 8$, which results in a leaf size of $m = 25$. 
% K-D Tree based space-partitioning and ordering technique is used to convert $80\times 80$ 2D-grid input-output pairs into $6400$ 1D grid input-output pairs. The K-D tree has $\kappa$ levels of hierarchy and each of the leaf nodes in last level contain $m$ values. 
We trained FDSNet using NAdam optimiser with batch size 128, and early stopping with patience 200 epochs. 
An illustration of a sample $V$ from the test data, the ground truth $u_G$, the prediction $u_{NN}$ by FDSNet, and the error with respect to $u_G$ is given in Figure~\ref{fig:nlse2d_fdsnet_prediction}.
We experimented with multiple depth $D$ and rank $p$ configurations and the results are tabulated in Table \ref{tab:FDSNet_sweep_all}. 
Increasing the rank from $2$ to $10$ reduces the error by $\approx 62\%$ at $\approx 4\times$ more parameters, while the depth sweep is non-monotone.

\subsubsection{Generalization across $\beta$ values (NLSE 1D)}\label{sec:beta_generalization}
In Equation~\eqref{eq:nlse}, $\beta$ scales the term $|u(x)|^{2}u(x)$ and controls how far the ground state $u_{G}$ departs from the linear Schr\"odinger regime. 
% As $\beta$ increases, $u_{G}$ deviates more strongly from a smooth Gaussian profile and instead tracks the potential wells of $V$ more sharply. 
% A natural question is whether a single trained FDSNet can predict $u_{G}$ accurately across a range of $\beta$ values, including values it has never seen during training.
In Section~\ref{sec:nlse_1d}, $\beta$ is set to 10. 
In this section, we study whether a single FDSNet can predict $u_{G}$ across a range of $\beta$, including values not seen during training.
To test this, we sample 101 values of $\beta$ on a uniform grid over $[5.0, 10.0]$ and 
generate $2{,}000$ potential and ground-state pairs $(V, u_{G})$ for each value of $\beta$. 
Four FDSNets, each with $N=320$, $\kappa=6$, $p=12$, $D=5$, and $26{,}127$ 
parameters, are trained and evaluated under the following 
train/test split protocol of the sampled $\beta$ values:
\begin{itemize}[leftmargin=1em]
    \item \textsf{P1} (dense interpolation): trains on every other $\beta$ value, 
          tests on the remaining $50$.
    \item \textsf{P2} (sparse interpolation): trains on every fifth $\beta$ value, 
          tests on the remaining $80$.
    \item \textsf{P3} (one-sided extrapolation): trains on $\beta \in [5.0, 8.0]$, 
          tests on $\beta \in [8.05, 10.0]$.
    \item \textsf{P4} (two-sided extrapolation): trains on $\beta \in [6.5, 8.5]$, 
          tests on both tails outside this interval.
\end{itemize}%
The models use NAdam optimizer with batch size $128$ and early stopping.
% , with at most $40{,}000$ training samples drawn from the training $\beta$ values. 
Validation error is monitored on a fixed $8{,}000$-sample 
subset drawn at the start of training from the held-out $\beta$ pool, and the weights achieving lowest validation loss are retained. 
% Final test metrics are computed over all samples at every held-out $\beta$ value.

The results are summarized in Table~\ref{tab:nlse_betas_results}. The following observations are noteworthy: (i) \textsf{P1} and \textsf{P2} achieve nearly identical test 
errors ($1.99\%$ and $1.98\%$, respectively), even though \textsf{P2} trains on roughly 
one-fifth as many distinct $\beta$ values as that of \textsf{P1}. This suggests the network has learned a representation that varies smoothly in $\beta$ rather than interpolating between memorized samples. (ii) Both extrapolation protocols 
(\textsf{P3} and \textsf{P4}) achieve $2.50$--$2.52\%$ test errors that are only slightly higher than that of the interpolation runs, despite testing entirely outside the training range. 
(iii) The fixed-$\beta$ result at $\beta = 10$ from Section~\ref{sec:nlse_1d} with $p=10$, $D=5$, (refer Table~\ref{tab:FDSNet_sweep_all}) gives roughly $0.98\%$ test error at the same configuration.
% , so conditioning on an unknown $\beta$ raises the error by a factor of roughly $2$--$2.5$.

% , a mild and predictable cost. 
\begin{table}[htbp]
\centering
% \renewcommand{\arraystretch}{1.1}
% \vspace{4pt}
\renewcommand{\arraystretch}{1.1}
\footnotesize
\begin{tabular}{llcc}
\midrule
\textbf{Protocol} 
    & \textbf{Train $\beta$ range} 
    & \textbf{No. of test} 
    & \textbf{Test error (mean$\pm$sd)} \\
    & 
    & \textbf{$\beta$ samples} 
    &  \\
\toprule 
% \hline
\textsf{P1} (dense interpolation)   
    & $[5.00, 10.00]$ & $50$ 
    & $1.99\!\times\!10^{-2}\!\pm\!7.9\!\times\!10^{-3}$ \\
\textsf{P2} (sparse interpolation)  
    & $[5.00, 10.00]$ & $80$ 
    & $1.98\!\times\!10^{-2}\!\pm\!8.0\!\times\!10^{-3}$ \\
\textsf{P3} (1-sided extrapolation) 
    & $[5.00, 8.00]$  & $40$ 
    & $2.52\!\times\!10^{-2}\!\pm\!5.7\!\times\!10^{-3}$ \\
\textsf{P4} (2-sided extrapolation) 
    & $[6.50, 8.50]$  & $60$ 
    & $2.50\!\times\!10^{-2}\!\pm\!7.4\!\times\!10^{-3}$ \\
\midrule
\end{tabular}
\caption{Generalization of FDSNet across $\beta$ in the NLSE 1D problem.
% All protocols use the same FDSNet with $(N{=}320,\,\kappa{=}6,\,p{=}12,\,D{=}5)$ and $26{,}127$ parameters.
% Test MSE and Rel.\ $\ell^{2}$ aggregate over all samples at all test $\beta$ values; the last column gives the mean $\pm$ standard deviation of the per-$\beta$ relative $\ell^{2}$ error across held-out $\beta$ values.
The mean and standard deviation of the test errors are reported in the last column.
}
\label{tab:nlse_betas_results}
\end{table}

\subsubsection{Comparison of FDSNet with a classical numerical solver}
We solved the ground state of the NLSE problem using a normalized gradient flow solver\\~\cite{baodu2004} and compared its run time to that of the inference time of FDSNet, at a matched accuracy. For this experiment, we trained FDSNet using the \textsf{P1} protocol described in Section~\ref{sec:beta_generalization}. 
When averaged over test samples spanning several $\beta$ values, the normalized gradient flow solver took $0.2373$~ms per sample, while FDSNet took $0.034$~ms, on the same GPU. FDSNet is roughly $~7\times$ faster, justifying the offline training cost.

\begin{table}[htbp]
\centering
\renewcommand{\arraystretch}{1.1}
\setlength{\tabcolsep}{7pt}
\footnotesize

\begin{tabular}{llccrc}
\midrule
\textbf{PDE} & \textbf{Sweep} & $p$ & $D$ 
& \textbf{Parameter} & \textbf{Test error} \\
 & &  &  
& \textbf{count} &  \\
\toprule
\multirow{8}{*}{\shortstack{NLSE 1D\\{\footnotesize $(N=320,\ \kappa=6,\ m=5)$}}}
& \multirow{5}{*}{$p$ with $D=5$}
& 2  & 5 & 2,367  & $3.95 \times 10^{-2}$ \\
& & 4  & 5 & 5,199  & $1.72 \times 10^{-2}$ \\
& & 6  & 5 & 8,991  & $1.39 \times 10^{-2}$ \\
& & 8  & 5 & 13,743 & $1.16 \times 10^{-2}$ \\
& & 10 & 5 & 19,455 & $9.80 \times 10^{-3}$ \\

\cmidrule(lr){2-6}

& \multirow{3}{*}{$D$ with $p=10$}
& 10 & 3 & 14,355 & $1.16 \times 10^{-2}$ \\
& & 10 & 5 & 19,455 & $9.80 \times 10^{-3}$ \\
& & 10 & 7 & 24,555 & $1.01 \times 10^{-2}$ \\

\toprule

\multirow{8}{*}{\shortstack{NLSE 2D\\{\footnotesize $(N=6400,\ \kappa=8,\ m=25)$}}}
& \multirow{5}{*}{$p$ with $D=7$}
& 2  & 7 & 38,211  & $7.17 \times 10^{-2}$ \\
& & 4  & 7 & 66,639  & $4.85 \times 10^{-2}$ \\
& & 6  & 7 & 96,859  & $3.61 \times 10^{-2}$ \\
& & 8  & 7 & 128,871 & $2.85 \times 10^{-2}$ \\
& & 10 & 7 & 162,675 & $2.72 \times 10^{-2}$ \\

\cmidrule(lr){2-6}

& \multirow{3}{*}{$D$ with $p=10$}
& 10 & 3 & 146,635 & $2.72 \times 10^{-2}$ \\
& & 10 & 5 & 154,655 & $2.78 \times 10^{-2}$ \\
& & 10 & 7 & 162,675 & $2.72 \times 10^{-2}$ \\

\toprule

\multirow{8}{*}{\shortstack{Burgers' 1D\\{\footnotesize $(N=1024,\ \kappa=7,\ m=8)$}}}
& \multirow{5}{*}{$p$ with $D=5$}
& 2  & 5 & 6,226  & $1.33 \times 10^{-2}$ \\
& & 4  & 5 & 12,124 & $1.86 \times 10^{-3}$ \\
& & 6  & 5 & 19,142 & $1.41 \times 10^{-3}$ \\
& & 8  & 5 & 27,280 & $1.28 \times 10^{-3}$ \\
& & 10 & 5 & 36,538 & $1.11 \times 10^{-3}$ \\

\cmidrule(lr){2-6}

& \multirow{3}{*}{$D$ with $p=6$}
& 6 & 3 & 16,814 & $1.73 \times 10^{-3}$ \\
& & 6 & 5 & 19,142 & $1.41 \times 10^{-3}$ \\
& & 6 & 7 & 21,470 & $2.72 \times 10^{-3}$ \\

\toprule

\multirow{6}{*}{\shortstack{Darcy flow 2D\\{\footnotesize $(N=9216,\ \kappa=10,\ m=9)$}}}
& \multirow{3}{*}{$p$ with $D=7$}
& 6  & 7 & 131,391 & $5.04 \times 10^{-2}$ \\
& & 9  & 7 & 199,683 & $4.68 \times 10^{-2}$ \\
& & 12 & 7 & 273,015 & $7.77 \times 10^{-2}$ \\

\cmidrule(lr){2-6}

& \multirow{3}{*}{$D$ with $p=6$}
& 6 & 3 & 124,791 & $5.91 \times 10^{-2}$ \\
& & 6 & 5 & 128,091 & $5.61 \times 10^{-2}$ \\
& & 6 & 7 & 131,391 & $5.04 \times 10^{-2}$ \\

\bottomrule
\end{tabular}

\caption{FDSNet hyperparameter sweeps across four PDE benchmarks. For each problem, the upper sub-block varies $p$ at fixed $D$, while the lower sub-block varies $D$ at fixed $p$.}

\label{tab:FDSNet_sweep_all}
\end{table}

\subsection{Burgers' equation in 1D}\label{sec:burgers}
% The viscous Burgers' equation is a canonical nonlinear advection--diffusion
% PDE exhibiting shock formation and viscous dissipation. We consider the 1D
% form
% \begin{equation}
% \frac{\partial u}{\partial t} + u\,\frac{\partial u}{\partial x}
% = \nu\,\frac{\partial^{2} u}{\partial x^{2}},
% \qquad x \in [0, 2\pi),\quad t > 0,
% \label{eq:burgers}
% \end{equation}
% with periodic boundary conditions. We learn the solution operator
% $\mathcal{M}\colon u_{0}(x)\mapsto u(x,T)$ mapping initial conditions to the
% solution at a fixed terminal time $T$. As in NLSE, the periodic translation
% invariance and resulting (block-)Toeplitz structure justify replacing the
% locally-connected blocks of FDSNet with shared convolutional kernels.

The viscous Burgers' equation is a standard nonlinear advection-diffusion PDE exhibiting steep gradient formation and viscous dissipation. 
We consider the one-dimensional form
\begin{equation}
    \frac{\partial u}{\partial t} + u\,\frac{\partial u}{\partial x}
    = \nu\,\frac{\partial^{2} u}{\partial x^{2}},
    \qquad x \in [0, 2\pi),\quad t > 0,
    \label{eq:burgers}
\end{equation}
with periodic boundary conditions, where $\nu > 0$ is the viscosity 
coefficient. We learn the solution operator $\mathcal{S}\colon u_{0}(x) 
\mapsto u(x, T)$ mapping the initial condition $u_0(x) = u(x,0)$ to the solution at a fixed terminal time $T$. 
% When the distribution of $u_0$ is translation-invariant on $[0, 2\pi)$, the solution operator $\mathcal{M}$ is translation-equivariant, and the discretized operator carries a Toeplitz structure on a uniform grid. 
The solution operator $\mathcal{S}$ is translation-equivariant on $[0, 2\pi)$ due to the periodic and spatially homogeneous structure of Equation~\eqref{eq:burgers}, and the discretized operator exhibits a Toeplitz structure on a uniform grid. We exploit this structure by replacing the locally-connected blocks of FDSNet with shared convolutional kernels. 
The ground-truth solutions are obtained by numerically solving Equation~\eqref{eq:burgers} using an explicit forward-Euler with second-order central finite difference scheme.

\paragraph{Setup} 
% Analogous to the NLSE case, we exploit this by replacing the locally-connected blocks of FDSNet with shared convolutional kernels.
We set viscosity $\nu = 0.03$ and terminal time $T = 1.0$.
FDSNet is configured with $N=1024$, and $\kappa=7$, resulting in leaf size $m=8$. 
% We experimented with multiple $D$ and $p$ values for the network. 
We use NAdam optimiser for training with minibatch size $B=128$, maximum epochs 2000, learning rate $1\times 10^{-3}$ and early stopping with patience of 200 epochs. Both train and test datasets contain $20$K samples.

%The best epoch is determined by minimising validation loss and the model weights of best epoch are used for inference and evaluation metrics. 
\paragraph{Results} 
Table~\ref{tab:FDSNet_sweep_all} reports the rank sweep 
($p \in \{2, 4, 6, 8, 10\}$ at fixed $D = 5$) and depth sweep 
($D \in \{3, 5, 7\}$ at fixed $p = 6$). 
Increasing the rank from $4$ to $10$ reduces the error by $\approx 40\%$ at $\approx 3\times$ more parameters, while the depth sweep is non-monotone, with $D = 5$ outperforming both $D = 3$ and $D = 7$. 
The best configuration $(p = 10,\, D = 5)$ achieves a test error of $\approx 0.11\%$, despite the underlying steep-gradient dynamics.

%The rank sweep shows that increasing the rank from $4$ to $10$ reduces the error by roughly $40\%$ at the cost of roughly $3\times$ more parameters. The depth sweep exhibits a non-monotone trend, with $D = 5$ outperforming both $D = 3$ and $D = 7$, suggesting that deeper nonlinear stacks at this rank begin to overfit. The best configuration $(p = 10,\, D = 5)$ achieves a test error of $ \approx 0.11\%$, demonstrating that FDSNet effectively captures the Burgers solution operator.
% despite the underlying steep-gradient dynamics.

\subsection{Darcy's flow in 2D}\label{sec:darcy2dflow}
Darcy's flow describes the steady-state pressure $u(x,y)$ driven by a 
heterogeneous permeability field $a(x,y)$ in a porous medium. We consider 
the elliptic boundary value problem
\begin{equation}
    -\nabla \cdot \bigl(a(x,y)\,\nabla u(x,y)\bigr) = f(x,y),
    \qquad (x,y) \in (0,1)^{2}, \quad u\big|_{\partial(0,1)^{2}} = 0,
    \label{eq:darcy}
\end{equation}
with constant forcing $f = 1$. The permeability field is sampled from a 
Gaussian random field (GRF) with Whittle--Mat\'{e}rn covariance
$C = \sigma^{2}(-\Delta + \tau^{2}I)^{-\alpha}$, with $\sigma^2 = 1$, 
$\alpha = 2$, $\tau = 3$, and thresholded to yield a piecewise-constant 
high-contrast field
\[
a(x,y) = \begin{cases} 12 & \text{if } \mathrm{GRF}(x,y) \geq 0, \\
                        4  & \text{otherwise.} \end{cases}
\]
The task is to learn the solution operator $\mathcal{S}\colon a \mapsto u$. 
Although Equation~\eqref{eq:darcy} is linear in $u$, the operator 
$\mathcal{S}$ is nonlinear in $a$.
Therefore we model it using the non-linear FDSNet.
% The permeability field $a$ is input-dependent and its distribution is not translation-invariant. Moreover, the Dirichlet boundary conditions break periodicity. We therefore retain the locally-connected layers of FDSNet without replacement by convolutional kernels. 
Although the homogeneous Dirichlet boundary conditions break translation equivariance, we impose translation equivariance as an inductive bias by replacing the locally-connected layers of FDSNet with shared convolutional kernels. 

% Darcy's flow models steady-state pressure $u(x,y)$ driven by a heterogeneous
% permeability field $a(x,y)$ in a porous medium. We consider the elliptic
% boundary value problem
% \begin{equation}
% -\nabla \cdot \bigl(a(x,y)\,\nabla u(x,y)\bigr) = f(x,y),
% \qquad (x,y)\in(0,1)^{2},\quad u\big|_{\partial(0,1)^{2}}=0,
% \label{eq:darcy}
% \end{equation}
% with constant forcing $f= 1$. The permeability field is sampled from a
% Gaussian random field with Mat\'{e}rn-type covariance
% $C = \sigma^{2}(-\Delta + \tau^{2} I)^{-\alpha}$, with $\alpha=2$, $\tau=3$,
% % matching the FNO benchmark~\cite{li2021fourier})
% and thresholded to yield a
% piecewise-constant high-contrast field
% \[
% a(x,y) = \begin{cases} 12 & \text{if } \mathrm{GRF}(x,y) \geq 0,\\
%                        4  & \text{otherwise.}\end{cases}
% \]
% The task is to learn the operator $\mathcal{M}\colon a \mapsto u$. Unlike the NLSE and Burgers' problems, Darcy's permeability field $a$ is sample-specific and is not translation-invariant, so we use the locally-connected layers in the FDSNet.

\paragraph{Setup} The domain $(0,1)^{2}$ is discretised on a uniform
$N_{x}\times N_{y}$ interior grid with $N_{x}=96$, resulting in $N=9{,}216$ degrees of freedom. 
Ground-truth solutions are generated by discretising Equation~\eqref{eq:darcy} with a five-point cell-centred finite-difference scheme using arithmetic-mean interface coefficients on an interior grid, and the resulting sparse symmetric positive-definite linear system is solved directly via a sparse LU factorisation.
% Ground-truth pressures are obtained by finite-difference discretisation of~\eqref{eq:darcy}.
% ; data-generation code is provided in the accompanying repository. 
The dataset comprises of $20{,}000$ training and $20{,}000$ test pairs $(a, u)$.
% , where $a$ and $u$ denote the discretisations of $a$ and $u$, respectively. 
% The 2D input is linearised to a length-$N$ vector via the $k$-d tree ordering ($k=2$) of Section~\ref{ssec:linear_2d} before being passed to Algorithm~\ref{alg:linear}. 
Training is performed using NAdam with \(\text{learning rate } 2\times10^{-4}\), batch size \(128\), and early stopping with a patience of \(200\), for a maximum of \(1000\) epochs. 
The hierarchical depth is set to \(\kappa=10\), resulting in a leaf size of \(m = 9\).

\paragraph{Results} Table~\ref{tab:FDSNet_sweep_all} reports the rank sweep ($p\in\{6,9,12\}$ at fixed $D=7$) and depth sweep ($D\in\{3,5,7\}$ at fixed $p=6$). The depth sweep at $p = 6$ is monotone, with deeper nonlinear stacks consistently improving performance, and $D = 7$ improves over $D = 3$ by roughly $15\%$ in test error. The rank sweep is non-monotone.
% since $p = 9$ improves over $p = 6$, but $p = 12$ diverged early in training. 
The best configuration $(p = 9,\, D = 7)$ achieves a test error of $4.68\%$ at just $199\text{K}$ parameters. 

\begin{table}[htbp]
\begin{minipage}{\textwidth}
% \centering
\renewcommand{\arraystretch}{1.2}
\resizebox{\textwidth}{!}{%
\begin{tabular}{lllrrrr}
\midrule
\textbf{PDE} & \textbf{Model} & \textbf{Configuration} & \textbf{Parameter} & \textbf{Train error} & \textbf{Test error} & \textbf{Inference} \\
 &  &  & \textbf{count} &  &  & \textbf{time (ms)} \\
\toprule
% --- nlse_1d 20k ---
NLSE 1D & FNO      & $M{=}16,\,W{=}64,\,D{=}4$ & $287{,}425$ & $1.07\!\times\!10^{-3}$ & $6.45\!\times\!10^{-4}$ & $0.031$ \\
 & MNN & $L{=}6,\,\alpha{=}10,\,K{=}5$ & $20{,}455$ & $4.21\!\times\!10^{-3}$ & $3.27\!\times\!10^{-3}$ & $0.032$ \\
 & FDSNet                            & $\kappa{=}6,\,p{=}10,\,D{=}5$ & $19{,}455$ & $1.11\!\times\!10^{-2}$ & $1.06\!\times\!10^{-2}$ & $0.035$ \\
 & MLP                                    & $W{=}35,\,D{=}3$ & $25{,}275$ & $7.44\!\times\!10^{-2}$ & $7.48\!\times\!10^{-2}$ & $0.005$ \\
 & DeepONet         & $p{=}128,\,D{=}3$ & $124{,}033$ & $1.02\!\times\!10^{-1}$ & $1.04\!\times\!10^{-1}$ & $0.005$ \\
\toprule
% --- burgers_1d 20k ---
Burgers 1D & FNO      & $M{=}16,\,W{=}64,\,D{=}4$ & $287{,}425$ & $3.77\!\times\!10^{-3}$ & $8.23\!\times\!10^{-4}$ & $0.042$ \\
 & FDSNet                & $\kappa{=}7,\,p{=}10,\,D{=}5$ & $36{,}538$ & $3.07\!\times\!10^{-3}$ & $9.57\!\times\!10^{-4}$ & $0.043$ \\
 & MNN         & $L{=}7,\,\alpha{=}10,\,K{=}5$ & $31{,}944$ & $3.96\!\times\!10^{-3}$ & $9.94\!\times\!10^{-4}$ & $0.039$ \\
 & MLP                                    & $W{=}18,\,D{=}3$ & $38{,}590$ & $2.39\!\times\!10^{-2}$ & $1.85\!\times\!10^{-2}$ & $0.003$ \\
 & DeepONet         & $p{=}128,\,D{=}3$ & $214{,}145$ & $1.08\!\times\!10^{-1}$ & $8.56\!\times\!10^{-2}$ & $0.005$ \\
\toprule
% % --- darcy_2d 5k ---
% Darcy 2D & 5k  & FNO        & $M{=}12,\,W{=}32,\,D{=}4$ & $1{,}188{,}000$ & $1.80\!\times\!10^{-2}$ & $2.04\!\times\!10^{-2}$ & $0.606$ \\
%  &   & FDSNet                  & $\kappa{=}10,\,p{=}9,\,D{=}7$ & $200{,}000$ & $5.96\!\times\!10^{-2}$ & $6.10\!\times\!10^{-2}$ & $0.137$ \\
%  &   & DeepONet           & $p{=}128,\,D{=}4$ (CNN br.) & $7{,}498{,}000$ & $4.93\!\times\!10^{-2}$ & $7.82\!\times\!10^{-2}$ & $0.035$ \\
%  &   & MLP                                      & $W{=}11,\,D{=}3$ & $212{,}000$ & $1.50\!\times\!10^{-1}$ & $1.50\!\times\!10^{-1}$ & $0.004$ \\
%  &   & MNN           & $L{=}5,\,\alpha{=}15,\,K{=}5$ & $221{,}000$ & $5.37\!\times\!10^{-1}$ & $5.17\!\times\!10^{-1}$ & $0.088$ \\
% \hline
% --- darcy_2d 20k ---
Darcy 2D & FNO        & $M{=}12,\,W{=}32,\,D{=}4$ & $1{,}188{,}353$ & $1.34\!\times\!10^{-2}$ & $6.65\!\times\!10^{-3}$ & $0.275$ \\
 & FDSNet                  & $\kappa{=}10,\,p{=}9,\,D{=}7$ & $199{,}683$ & $5.69\!\times\!10^{-2}$ & $4.89\!\times\!10^{-2}$ & $0.087$ \\
 & DeepONet           & $p{=}128,\,D{=}4$ (CNN br.) & $7{,}498{,}497$ & $6.13\!\times\!10^{-2}$ & $6.29\!\times\!10^{-2}$ & $0.042$ \\
 & MLP                                      & $W{=}11,\,D{=}3$ & $212{,}243$ & $1.50\!\times\!10^{-1}$ & $1.50\!\times\!10^{-1}$ & $0.004$ \\
 & MNN           & $L{=}5,\,\alpha{=}15,\,K{=}5$ & $221{,}265$ & $2.99\!\times\!10^{-1}$ & $4.54\!\times\!10^{-1}$ & $0.053$ \\
\toprule
% --- nlse_2d 20k ---
NLSE 2D & FNO         & $M{=}12,\,W{=}32,\,D{=}4$ & $1{,}188{,}353$ & $2.63\!\times\!10^{-3}$ & $1.19\!\times\!10^{-3}$ & $0.362$ \\
 & FDSNet                   & $\kappa{=}8,\,p{=}10,\,D{=}7$ & $162{,}675$ & $2.71\!\times\!10^{-2}$ & $2.71\!\times\!10^{-2}$ & $0.067$ \\
 & MNN            & $L{=}4,\,\alpha{=}12,\,K{=}7$ & $176{,}947$ & $4.21\!\times\!10^{-2}$ & $3.52\!\times\!10^{-2}$ & $0.037$ \\
 & MLP                                       & $W{=}13,\,D{=}3$ & $173{,}177$ & $3.26\!\times\!10^{-1}$ & $3.25\!\times\!10^{-1}$ & $0.004$ \\
 & DeepONet            & $p{=}128,\,D{=}4$ (CNN br.) & $5{,}008{,}385$ & $3.44\!\times\!10^{-1}$ & $3.43\!\times\!10^{-1}$ & $0.007$ \\
\bottomrule
\end{tabular}
}
\caption{\footnotesize Comparison of FDSNet (the proposed architecture) with existing neural operator learning architectures across various PDEs in 1D and 2D. For each PDE, networks are listed in ascending order of test error.\\
    \textit{Symbols.} For FNO, $M$ is the number of retained Fourier modes per spatial axis, $W$ is the channel width, and $D$ is the number of Fourier layers~\cite{li2021fourier}. For FDSNet, $\kappa$ is the number of HODLR levels, $p$ is the off-diagonal block rank, and $D$ is the depth of $K^{(\kappa)}$ and $S^{(k)}$ stacks. For MNN, $L$ is the number of hierarchical levels, $\alpha$ is the channel/rank parameter, and $K$ is the depth of the non-linear stack~\cite{fan2019multiscale}. 
    % For 2D problems, MNN operates on the native 2D grid, so its $L$ is approximately half of FDSNet's flattened-1D $L$ (e.g.\ FDSNet $L{=}10$ on a $96\!\times\!96$ grid versus MNN $L{=}5$). 
    For DeepONet, $p$ is the output basis dimension and $D$ is the trunk depth~\cite{lu2021deeponet}. For MLP, $W$ is the hidden width and $D$ is the depth.
    % }
    % Inference times are wall-clock per-sample timings on the test set (ms) under identical hardware.
    % \item 
    % {\footnotesize
    \newline
    \textit{Configuration.} 
    % The architectures differ in how 1D and 2D problems are handled. 
    FNO uses distinct 1D and 2D variants, with the spectral convolution applied via \texttt{rfft} (1D) or \texttt{rfft2} (2D)~\cite{li2021fourier}. The spectral-weight tensor scales as $O(W^{2}M)$ in 1D and $O(W^{2}M^{2})$ in 2D, which is why the reported 2D configurations use smaller $M$ and $W$. MNN similarly uses separate 1D and 2D variants and operates on a 2D grid for 2D problems via a quadtree.
    % (one MNN level corresponds to a $2\!\times\!2$ partition, so $4^{L}$ leaves at depth $L$). 
    DeepONet uses an MLP branch for 1D function inputs and a CNN branch for 2D function inputs (denoted ``CNN br.''\ in the table), as recommended by Lu et al.~\cite{lu2021deeponet}. FDSNet, by contrast, uses the same architecture irrespective of dimensions. 
    It uses a K-D tree with depth $\kappa$. 
    % regimes: 2D inputs are linearized via a $k$-d tree ordering and processed as length-$N$ vectors. 
    As a consequence, on 2D problems FDSNet's $\kappa$ counts levels of a binary partition of the index set and is roughly twice the corresponding 2D MNN's count of levels $L$ at a matched leaf size (e.g.\ FDSNet with $\kappa{=}10$ and MNN with $L{=}5$, on a $96\!\times\!96$ grid, yield a leaf size of 9).
}
% Comparison of architectures across the NLSE \& Burgers 1D, 2D Darcy flow, and 2D NLSE benchmarks. Relative $\ell^{2}$ denotes the non-squared per-sample relative error used for evaluation. Inference time is measured per sample on the test set.
\label{tab:multi_benchmark_comparison}

\end{minipage}
\end{table}

\subsection{Comparison of FDSNet with existing neural operator architectures}\label{sec:ablation}
We compare FDSNet with four existing architectures: Multiscale Neural Network 
(MNN) of Fan et al.~\cite{fan2019multiscale}, Fourier Neural Operator 
(FNO) of Li et al.~\cite{li2021fourier}, DeepONet of Lu et 
al.~\cite{lu2021deeponet}, and a parameter-matched feedforward MLP with same width across all hidden layers. 
% The comparison is conducted across three PDE benchmarks: NLSE 1D, Burgers' 1D, and NLSE 2D. All models are trained using a training dataset of size $20,000$. 
The comparison is conducted across four PDE benchmarks: NLSE 1D, Burgers' 1D, NLSE 2D, and Darcy flow 2D. 
All models are trained with NAdam using a training dataset of size $20,000$ at batch size $128$ and early stopping. The learning rate is $10^{-3}$ for all problems except Darcy, where we use $2\times10^{-4}$. 
Each baseline network uses a configuration adopted from the corresponding original article and serves as a reasonable reference configuration rather than a tuned setting.
All models are trained on unnormalized input–output pairs.
% Training is driven by the mean-squared relative $\ell^{2}$ loss, which is itself scale-invariant in the target.
Test errors and inference times are averaged across multiple random seeds. We report the mean per-sample inference time measured on the test set.
% at batch size $32$
 
Table~\ref{tab:multi_benchmark_comparison} summarizes the results. FNO attains the lowest test error on every benchmark.
% , which is consistent with the smooth and, on three of the four benchmarks, translation-invariant structure of the problems, for which spectral parameterizations are known to be effective. 
On the 1D benchmarks, FNO uses roughly $7$--$15\times$ more parameters than FDSNet and incurs nearly equal per-sample inference times.
On the 2D benchmarks, FNO uses roughly $6$--$7\times$ more parameters than FDSNet and incurs higher per-sample inference times, by approximately $4\times$ on Darcy flow and $5.4\times$ on NLSE 2D.
% On NLSE 2D, FNO uses roughly $6$--$7\times$ more parameters than FDSNet and incurs higher per-sample inference time, by approximately $5.4\times$.
% FDSNet achieves intermediate accuracy across all benchmarks while maintaining a substantially smaller parameter count and inference-time cost.

On the 1D benchmarks, FDSNet and MNN perform comparably at similar 
parameter counts. On Burgers' 1D, both achieve test errors within a small margin of FNO ($9.57\times10^{-4}$ and 
$9.94\times10^{-4}$, respectively, versus $8.23\times10^{-4}$ for FNO) 
despite using roughly $8\times$ fewer parameters. On NLSE 1D, MNN 
outperforms FDSNet on test error ($3.27\times10^{-3}$ versus $1.06\times10^{-2}$) at a comparable parameter count. Both architectures exploit the translation 
invariance of the NLSE 1D ground-state map via shared convolutional kernels, 
so the difference between them reflects how each parameterizes the operator.
MNN learns the forward map of a hierarchical matrix, while FDSNet learns its 
inverse factorization. At this problem size and parameter budget, the forward 
parameterization appears easier to fit.

On the 2D benchmarks, FDSNet is the second most accurate architecture 
after FNO, while using a substantially smaller parameter budget than both FNO 
and DeepONet. On Darcy's flow, FDSNet achieves a test error of 
$4.89\times10^{-2}$, which is lower than $6.29\times10^{-2}$ obtained by DeepONet
(at roughly $37\times$ the parameter count) and is also lower than $1.50\times10^{-1}$ obtained by MLP. MNN does not reach competitive accuracy on Darcy under the baseline configuration used here. On NLSE 2D, FDSNet achieves a test error of $2.71\times10^{-2}$ at $163$k parameters, with DeepONet, MNN, and the MLP all reporting higher test errors. 
% On NLSE 2D, FDSNet is the second most accurate architecture after FNO,
% while using a substantially smaller parameter budget than both FNO and
% DeepONet. FDSNet achieves a test error of $2.71\times10^{-2}$ at $163$k 
% parameters, with DeepONet, MNN, and the MLP all reporting higher test errors.

Overall, the results show FDSNet is parameter-efficient and achieves competitive accuracy while maintaining low inference time.

\section{Conclusion}\label{sec:conclusion}
We presented linear FDSNet, a neural network for performing the inverse operation of a HODLR matrix, inspired by a fast direct solver. 
The architecture uses locally connected (LC) layers with linear activations. It has a parameter complexity of $\mathcal{O}(pN \log N)$, the same complexity as that of the solve phase of its classical fast direct solver counterpart. 
When the operator being learnt is translation-invariant, LC layers are replaced with CNN layers, further reducing complexity to $\mathcal{O}(pN)$.
Further, for learning non-linear solution operators associated with PDEs, we presented the non-linear FDSNet. It extends the linear FDSNet by replacing selected linear layers with stacks of layers where non-linear activations are used in all layers except the top layer, thereby enabling the network to learn the operator's nonlinearity.
% For 2D case, we linearize the two-dimensional grid map using K-d trees and use Algorithm~\ref{alg:linear} \& ~\ref{alg:non_linear}. 

We performed a comprehensive set of numerical experiments to validate the proposed architecture. Various inferences drawn from the experiments include: 
(1) FDSNet efficiently solved the Fredholm integral equation in 2D, nonlinear Schr\"odinger equation in 1D \& 2D, Burgers' equation in 1D, and Darcy's flow problem in 2D.
(2) In the NLSE 1D problem, FDSNet demonstrated strong generalisation across a parameter sweep. When trained on a dataset spanning multiple values of $\beta$, the network accurately predicted solutions for values of $\beta$ not encountered during training.
(3) Furthermore, on the NLSE 1D problem, the inference time of FDSNet is observed to be lower than the run time of a classical numerical solver. Classical solvers must be re-run independently for a new parameter value, whereas a single FDSNet trained over a parameter range can infer solutions at a fraction of the computational cost, highlighting its practical advantage over classical approaches.
% These results demonstrate that FDSNet maintains the computational efficiency of classical hierarchical methods while providing flexibility to learn solution operators for both linear and nonlinear problems. 
(4) The comparison of FDSNet with existing neural operator learning architectures shows that FDSNet achieves competitive accuracy using only a fraction of the parameters required by existing neural operator learning architectures. This highlights the benefit of incorporating algebraic structure inspired by fast direct solvers directly into the network design. 

In future, we aim to enhance the FDSNet architecture to improve its performance. Possible enhancements include: 
(1) In higher dimensions, the off-diagonal block rank of HODLR matrix grows in powers of $N$. As a result, the parameter complexity of FDSNet grows in powers of $N$ in higher dimensions. To build architectures with better parameter complexity and thereby better inference times, networks inspired by hierarchical matrices built on the strong admissibility condition are to be considered. In particular, we would like to extend it to the fast direct solver for FMM-able (Fast Multipole Method) matrices~\cite{gujjula2024algebraic}.
(2) When approximating oscillatory kernels, hierarchical matrices perform well in low to moderate frequency regimes, but they possess challenges in high frequency regime. To address this, we would like to extend the proposed architecture to the directional FMM framework~\cite{gujjula2023dafmm}. 

\bibliographystyle{siamplain}
\bibliography{references}
\end{document}